\documentclass[a4paper]{amsart}

\usepackage{amssymb,mathrsfs}
\usepackage[mathic=true]{mathtools}
\usepackage{enumitem}
\usepackage{hyperref}

\theoremstyle{plain}%
\newtheorem{theorem}{Theorem}
\newtheorem{proposition}[theorem]{Proposition}%
\newtheorem{lemma}[theorem]{Lemma}%
\newtheorem{corollary}[theorem]{Corollary}%

\theoremstyle{remark}%
\newtheorem{example}{Example}%

\theoremstyle{definition}%
\newtheorem{definition}{Definition}%

\let\emptyset\varnothing

\newcommand{\noarg}{\bullet}

\newcommand{\reals}{\mathbb{R}}
\newcommand{\extreals}{\overline{\mathbb{R}}}
\newcommand{\nnegreals}{\mathbb{R}_{\geq 0}}
\newcommand{\posreals}{\mathbb{R}_{>0}}
\newcommand{\nats}{\mathbb{N}}
\newcommand{\nnegints}{\mathbb{Z}_{\geq 0}}

\DeclarePairedDelimiter{\pr}{(}{)}
\DeclarePairedDelimiter{\st}{\{}{\}}
\DeclarePairedDelimiter{\ev}{\{}{\}}

\DeclarePairedDelimiter{\abs}{\vert}{\vert}
\DeclarePairedDelimiter{\card}{\vert}{\vert}
\DeclarePairedDelimiterX{\ooi}[2]{]}{[}{#1, #2}
\DeclarePairedDelimiterX{\oci}[2]{]}{]}{#1, #2}
\DeclarePairedDelimiterX{\coi}[2]{[}{[}{#1, #2}
\DeclarePairedDelimiterX{\cci}[2]{[}{]}{#1, #2}

\newcommand{\compl}{\mathrm{c}}

\newcommand{\stsp}{\mathcal{X}}
\newcommand{\timeaxis}{\mathbb{T}}
\newcommand{\timepoints}{\mathcal{T}}
\newcommand{\altpoints}{\mathcal{S}}

\DeclareMathOperator{\isos}{is}
\DeclareMathOperator{\llims}{llims}
\DeclareMathOperator{\rlims}{rlims}
\DeclareMathOperator{\lims}{lims}
\newcommand{\densepoints}{\mathcal{D}}
\newcommand{\countpoints}{\mathcal{C}}
\newcommand{\setoftseq}{\mathcal{U}}
\newcommand{\pth}{\omega}
\newcommand{\allpth}{\dot{\omega}}
\newcommand{\cadpth}{\mathring{\omega}}
\newcommand{\vpth}{\varpi}
\newcommand{\rpth}{\psi}
\newcommand{\paths}{\Omega}
\newcommand{\allpaths}{\dot{\Omega}}
\newcommand{\genpaths}{\Omega}
\newcommand{\cadpaths}{\mathring{\Omega}}

\newcommand{\indic}{\mathbb{I}}
\newcommand{\indica}[1]{\indic_{#1}}

\newcommand{\genX}{X}
\newcommand{\gencylevts}{\mathscr{C}}
\newcommand{\genalgebra}{\mathscr{A}}
\newcommand{\allX}{\dot{X}}

\newcommand{\allcylevts}{\dot{\mathscr{C}}}
\newcommand{\cadX}{\mathring{X}\vphantom{X}}
\newcommand{\cadcylevts}{\mathring{\mathscr{C}}}

\newcommand{\posssp}{\mathfrak{S}}
\newcommand{\sample}{\mathfrak{s}}
\newcommand{\algebra}{\mathfrak{A}}
\newcommand{\borel}{\mathscr{B}}

\newcommand{\dist}{\mu}
\newcommand{\prob}{P}
\newcommand{\prev}{E}

\newcommand{\pthjumps}{\eta}
\newcommand{\tupjumps}{\hat{\eta}}

\newcommand{\allprob}{\dot{\prob}}

\begin{document}

\title[Countable-state stochastic processes with càdlàg sample paths]{Countable-state stochastic processes \\ with càdlàg sample paths}

\author{Alexander Erreygers}
\author{Jasper De~Bock}
\address{Foundations Lab, Ghent University, Belgium}

\thanks{This work was supported by the Research Foundation -- Flanders (FWO) (project number 3G028919).}

\begin{abstract}
	The Daniell--Kolmogorov Extension Theorem is a fundamental result in the theory of stochastic processes, as it allows one to construct a stochastic process with prescribed finite-dimensional distributions.
	However, it is well-known that the domain of the constructed probability measure -- the product sigma-algebra in the set of all paths -- is not sufficiently rich.
	This problem is usually dealt with through a modification of the stochastic process, essentially changing the sample paths so that they become càdlàg.
	Assuming a countable state space, we provide an alternative version of the Daniell--Kolmogorov Extension Theorem that does not suffer from this problem, in that the domain is sufficiently rich and we do not need a subsequent modification step: we assume a rather weak regularity condition on the finite-dimensional distributions, and directly obtain a probability measure on the product sigma-algebra in the set of all càdlàg paths.
\end{abstract}


\maketitle

\section{Introduction}
\label{sec:introduction}

A stochastic process is a model for a system whose state changes over time in an uncertain manner.
More formally, a stochastic process is a joint uncertainty model for a sequence~\(\pr{X_t}_{t\in\timeaxis}\) of \(\stsp\)-valued maps on some sample space~\(\posssp\) indexed by~\(\timeaxis\), where \(\stsp\) is the state space of the system and \(\timeaxis\) is an infinite subset of the real numbers that is interpreted as the time axis.
Stochastic processes have been and still are an active field of research, as is clear from the large number of monographs that have appeared on the subject \cite{1933Kolmogorov-Grund,1953Doob-Stochastic,1995Billingsley-Probability,1992Breiman-Probability,1967Cramer-Stationary,1969GikhmanSkorokhod-Introduction,1994Rogers-Diffusions}.
One usually derives the joint uncertainty model for~\(\pr{X_t}_{t\in\timeaxis}\) from the finite-dimensional distributions: the joint uncertainty models for~\(\pr{X_t}_{t\in T}\), where \(T\) ranges over all finite subsets of~\(\timeaxis\).
For example, Kolmogorov~\cite{1933Kolmogorov-Grund,1950Kolmogorov-Foundations} already considered a stochastic processes with \(\stsp=\reals\) in his seminal `Grundbegriffe', formulating what is now known as the Daniell--Kolmogorov Extension Theorem \cite[Chapter~III, §~4]{1950Kolmogorov-Foundations}.
However, it is well-known \cite[Section~5.1.2]{2006ShaferVovk} that if \(\timeaxis\) is not countable, then the domain of the resulting stochastic process -- the product \(\sigma\)-algebra generated by the cylinder events -- is not rich enough.

The first one to get around this problem was Doob~\cite{1953Doob-Stochastic}, who in his influential work discusses stochastic processes in a very general setting.
He explains why the product \(\sigma\)-algebra of events is not sufficiently rich \cite[Chapter~II, Section~2]{1953Doob-Stochastic}, but also comes up with a solution: the notion of `separability of a stochastic process'~\cite[pp.~51 and 52]{1953Doob-Stochastic}, which ensures that the product \(\sigma\)-algebra of events \emph{is} sufficiently rich.
In the particular case of `Markov processes with infinitely many states and continuous parameter', his Theorem~2.4~\cite[p.~266]{1953Doob-Stochastic} implies that the sample paths of any `separable Markov process' are almost surely `step functions'~\cite[p.~245 and 246]{1953Doob-Stochastic}.
Billingsley~\cite[Sections~36 and 37]{1995Billingsley-Probability} does something similar, but only for stochastic processes with the real numbers as state space: he also argues that the Daniell--Kolmogorov Extension Theorem does not suffice \cite[pp.~492--494]{1995Billingsley-Probability}, and solves this through his notion of `separability' \cite[pp.~526--527]{1995Billingsley-Probability}.

Like Doob~\cite{1953Doob-Stochastic} and Billingsley~\cite{1995Billingsley-Probability}, Breiman~\cite{1992Breiman-Probability} treats stochastic processes in the particular setting where the state space is the set of real numbers, and explains that the standard construction is insufficient.
He solves this problem a bit differently, though: he imposes a form of `continuity in probability'~\cite[Definition~12.15]{1992Breiman-Probability} on the stochastic process, and shows that in combination with `almost-sure absolute continuity' this allows the construction of a modification with continuous -- often also called regular -- sample paths \cite[Theorem~12.6]{1992Breiman-Probability}.
More recently, Borovkov~\cite[Section~18.2]{2013Borovkov-Probability} -- following Cramér \& Leadbetter~\cite{1967Cramer-Stationary} -- does not rely on `almost-sure absolute continuity', but argues that in order to get a modification with regular sample paths, one needs a quantitative bound on the continuity.
See also \cite{2015Cohen-Stochastic,2002Kallenberg-Foundations}.

Fewer authors consider stochastic processes with state spaces other than the real numbers.
For example, Gikhman and Skorokhod~\cite{1969GikhmanSkorokhod-Introduction} consider the very general setting of metric spaces for time domain and state space, but for most of their results they need compactness \cite[Chapter~4, Section~4]{1969GikhmanSkorokhod-Introduction}.
Rogers and Williams~\cite[Chapter~II]{1994Rogers-Diffusions} also give a very broad account, but they only solve the issue of the product \(\sigma\)-algebra in some particular cases, for example that of Markov (or Feller-Dynkyn) processes \cite[Chapter~III]{1994Rogers-Diffusions}.
Fristedt and Gray~\cite[Chapter~31]{1997Fristedt-Modern} define what they call a `(pure-jump) Markov process' on the set of càdlàg paths (for a general Polish state space), but they never really move past the (finitary) cylinder events.
Finally, König~\cite{2006Konig} gives a solution to this problem that does not involve separability or modifications: his Theorem~3.2 extends the probability measure on the product \(\sigma\)-algebra to a `maximal non-sequential inner extension' on a larger domain, but never really shows that this domain contains `all' the interesting events.

This brings us to our main contribution: we give a version of the Daniell--Kolmogorov Extension Theorem with a \(\sigma\)-algebra that is sufficiently rich.
We will do so for a non-empty and countable state space~\(\stsp\), and an arbitrary time domain~\(\timeaxis\subseteq\reals\).
As is customary, we interpret the elements of~\(\timeaxis\) as time points; two important cases are \(\timeaxis=\nats\) and \(\timeaxis=\nnegreals=\coi{0}{+\infty}\).
We intend to construct a joint uncertainty model for the state of the system in a (not necessarily proper) subset of the index set~\(\timeaxis\); henceforth, we will denote this subset by \(\timepoints\).
The usual choice is \(\timepoints=\timeaxis\), but other choices are sometimes useful as well. For example, in the setting of `Markovian imprecise jump processes' (also known as imprecise continuous-time Markov chains) \cite{2017KrakDeBock,2022Erreygers}, we encounter the case \(\timeaxis=\nnegreals\) and \(\timepoints=\st{t_1, \dots t_{n-1}}\cup\coi{t_n}{+\infty}\), with \(t_1, \dots, t_n\in\nnegreals\) such that \(t_1<\cdots<t_n\).

Our main result is Theorem~\ref{the:regular then daniell-kolmogorov for cadlag}, which is similar to the Daniell--Kolmogorov Extension Theorem in that it starts from a consistent family of finite-dimensional distributions, but also different in the following ways.
A first difference is that we specify the finite-dimensional distributions for finite subsets of~\(\timepoints\subseteq\timeaxis\) instead of~\(\timeaxis\).
A second and more important difference is that we not only require consistency of these finite-dimensional distributions, but also regularity (Definition~\ref{def:regularity}), and it is this that allows us to obtain a probability measure on the product \(\sigma\)-algebra for the regular (in this case càdlàg) paths instead of, as is customary, on the product \(\sigma\)-algebra for all paths.
This way, we get rid of the additional step of constructing a regular modification.

\subsection{Some notation regarding tuples}
For any non-empty subset~\(\altpoints\) of~\(\timeaxis\), we let \(\setoftseq_\altpoints\) be the set of all (non-empty) tuples~\(\pr{t_1, \dots, t_n}\in\bigcup_{k\in\nats} \altpoints^k\) that are increasing, so with \(t_1<\cdots<t_n\); if \(\altpoints=\timepoints\), we simply write \(\setoftseq\).
We will usually denote a generic tuple in the set~\(\setoftseq_{\timeaxis}\) by~\(u\), but sometimes also by \(v\) or \(w\).
For any two tuples of time points~\(u=\pr{r_1, \dots, r_n}\) and \(v=\pr{s_1, \dots, s_m}\) in~\(\setoftseq_{\timeaxis}\), we write \(u\sqsubseteq v\) if all time points in~\(u\) are included in~\(v\), in the sense that \(\st{r_1, \dots, r_n}\subseteq \st{s_1, \dots, s_m}\).

For any \(u=\pr{t_1, \dots, t_n}\in\setoftseq_{\timeaxis}\), we let~\(\stsp_u\) denote the set of all \(n\)-tuples of states~\(\pr{x_{t_1}, \dots, x_{t_n}}\) indexed by the time points in~\(u\), and we usually denote a generic \(n\)-tuple of states in~\(\stsp_u\) by \(x_u\), \(y_u\) or \(z_u\).
Furthermore, given \(u=\pr{t_1, \dots, t_n}\in \setoftseq_{\timeaxis}\) and \(x_u=\pr{x_{t_1}, \dots, x_{t_n}}\in \stsp_u\) and for all \(v=\pr{s_1, \dots, s_m}\in \setoftseq_{\timeaxis}\) such that \(v\sqsubseteq u\), we let \(x_v\coloneqq\pr{x_{s_1}, \dots, x_{s_m}}\) be the \(m\)-tuple in~\(\stsp_v\) that consists of those components of~\(x_u\) with index in~\(v\).
Rogers \& Williams~\cite[Chapter~II, Section 25 onwards]{1994Rogers-Diffusions} use slightly different notation: they consider finite subsets~\(U\) of~\(\timepoints\) instead of increasing finite tuples~\(u\) in~\(\setoftseq\subseteq\bigcup_{k\in\nats}\timepoints^k\), and focus on the set~\(\stsp^U\) of \(\stsp\)-valued maps from \(U\) to~\(\stsp\) instead of on~\(\stsp_u\).
It should be clear that these two approaches are essentially equivalent; nevertheless, we choose to use increasing tuples because these will come in handy in several places, for example in Sections~\ref{ssec:expected number of jumps} and \ref{sssec:regularity} and Appendix~\ref{assec:expected number of jumps again}.

\section{Constructing a stochastic process}
\label{sec:constructing a stochastic process}
We want to model a system whose state, which takes values in~\(\stsp\), changes along the time axis~\(\timeaxis\) in an uncertain manner.
For this, we turn to (measure-theoretical) probability theory, as for example outlined in \cite{1994Rogers-Diffusions,1995Billingsley-Probability,1997Fristedt-Modern,1992Breiman-Probability,2002Kallenberg-Foundations}.
This means that we set out to determine a suitable probability space~\(\pr{\posssp, \algebra, \prob}\), where the non-empty set~\(\posssp\) is called the \emph{sample space}, \(\algebra\) is a \(\sigma\)-algebra of \emph{events} -- subsets of~\(\posssp\) --  that are of interest to us, and \(\prob\) is a probability measure on~\(\algebra\).
As we will presently see, we will only consider couples~\(\pr{\posssp, \algebra}\) of a specific form; for example, the sample space will be a set of paths.

\subsection{Paths and cylinder events}
Since our system changes state over~\(\timeaxis\), it makes sense to think of an element~\(\sample\) of the sample space~\(\posssp\) as a map from \(\timeaxis\) to~\(\stsp\); we call such an \(\stsp\)-valued map on~\(\timeaxis\) a \emph{path}, and we denote the set of all paths by~\(\stsp^{\timeaxis}\).
Hence, the obvious choice for the sample space~\(\posssp\) is some (non-empty) set of paths~\(\genpaths\subseteq\stsp^{\timeaxis}\).
However, in many cases -- for example when using coherent conditional probabilities to deal with conditioning, see \cite[Eqn.~(12)]{2017KrakDeBock} -- it makes sense to require that
\begin{equation}
\label{eqn:condition on set of paths}
	\pr{\forall u\in\setoftseq}
	\pr{\forall x_u\in\stsp_u}
	\pr{\exists \pth\in\genpaths}~
	\pth\pr{u}=x_u,
\end{equation}
where here and in the remainder, \(\pth\pr{u}\coloneqq\pr[\big]{\pth\pr{t_1}, \dots, \pth\pr{t_n}}\) for all \(\pth\in\stsp^{\timeaxis}\) and \(u=\pr{t_1, \dots, t_n}\in\setoftseq\).
We think this is a reasonable (and harmless) requirement on~\(\genpaths\), so we henceforth assume it.

Now that we have established that \(\genpaths\subseteq\stsp^{\timeaxis}\) is a sensible sample space, the question remains which events we are interested in.
At the very least, we are interested in events regarding the state of the system in a finite number of time points in~\(\timepoints\).
Let us formalise these events.

For all \(t\in\timeaxis\), we define the projector variable or evaluation map \cite[Chapter~II, Eqn.~(25.1)]{1994Rogers-Diffusions}
\begin{equation*}
	\genX_t\colon\genpaths\to\stsp
	\colon \pth\mapsto \pth\pr{t};
\end{equation*}
we extend this notation in the obvious way to tuples of time points~\(u\in\setoftseq\):
\begin{equation*}
	\genX_u
	\colon \genpaths\to\stsp_u
	\colon \pth\mapsto
	\pth\pr{u}.
\end{equation*}
Then an event regarding the state of the system in a finite number of time points in~\(\timepoints\) is a \emph{cylinder event} \cite[Chapter~II, Definition~25.4]{1994Rogers-Diffusions}: an event of the form
\begin{equation*}
	\ev{\genX_u\in A}
	\coloneqq \st{\pth\in\genpaths\colon \genX_u\pr{\pth}\in A}
	= \st{\pth\in\genpaths\colon \pth\pr{u}\in A}
\end{equation*}
with \(u\in\setoftseq\) and \(A\in\wp\pr{\stsp_u}\).\footnote{$\wp\pr{\stsp_u}$ denotes the powerset of $\stsp_u$.}
The condition in Eqn.~\eqref{eqn:condition on set of paths} ensures that for all \(u\in \setoftseq\) and \(A\in\wp\pr{\stsp_u}\), the corresponding cylinder event~\(\ev{\genX_u\in A}\) is the empty set if and only if \(A=\emptyset\).
Furthermore, it is not difficult to verify that the collection
\begin{equation*}
	\gencylevts
	\coloneqq\st[\big]{\ev{\genX_u\in A}\colon u\in\setoftseq, A\in\wp\pr{\stsp_u}}
\end{equation*}
of cylinder events with time points in~\(\timepoints\) is an algebra of events.

In order to obtain the probability triple that we are after, it remains for us to determine (i) the subset~\(\genpaths\) of~\(\stsp^{\timeaxis}\) that we want to use as possibility space, and (ii) a probability measure~\(\prob\) on a \(\sigma\)-algebra~\(\genalgebra\) in~\(\genpaths\) that includes~\(\gencylevts\).
Such a probability triple is what we will call a stochastic process.
\begin{definition}
\label{def:stochastic process}
	A \emph{stochastic process} is a probability triple~\(\pr{\genpaths, \genalgebra, \prob}\) such that \(\genpaths\subseteq\stsp^{\timeaxis}\) and \(\gencylevts\subseteq\genalgebra\).
\end{definition}
Our definition of a stochastic process differs somewhat from the `usual' one, which considers a generic probability space~\(\pr{\posssp, \algebra,\prob}\) and a family~\(\pr{Y_t}_{t\in\timeaxis}\) of \(\algebra/\wp\pr{\stsp}\)-measurable maps from~\(\posssp\) to~\(\stsp\) -- see for example \cite[Chapter~II]{1953Doob-Stochastic}, \cite[Chapter~IV]{1969GikhmanSkorokhod-Introduction}, \cite[Definition~12.1]{1992Breiman-Probability}, \cite[Chapter~II, Definition~27.1]{1994Rogers-Diffusions}, \cite[Section~36]{1995Billingsley-Probability}, \cite[Section~1.2]{1997Fristedt-Modern} or \cite[Definition~18.1.2]{2013Borovkov-Probability}.
The first difference is that we choose to restrict the sample space~\(\posssp\) to a set of paths a priori.
The second difference is that we are only interested in modelling our uncertainty about the state in the time points in~\(\timepoints\subseteq\timeaxis\); in that sense, our definition corresponds to the usual one with \(\posssp=\genpaths\) and \(\algebra=\genalgebra\supseteq\gencylevts\), at least if we only consider the family~\(\pr{\genX_t}_{t\in\timepoints}\) instead of the one indexed by~\(\timeaxis\).

To construct a stochastic process, one usually turns to the Daniell--Kolmogorov Extension Theorem -- see \cite[Chapter~II, Theorem~31.1]{1994Rogers-Diffusions} or, if \(\stsp\) is finite, \cite[Section~1.5.2]{1980Iosifescu-Finite}.
The idea behind this theorem is simple: given a model for our uncertainty about~\(\genX_u\) for all \(u\in\setoftseq\), it constructs a (joint) uncertainty model for~\(\pr{\genX_t}_{t\in\timepoints}\).
Since we are dealing with a countable state space, the uncertainty models for~\(\genX_u\) are finite-dimensional charges or distributions.

\subsection{Finite-dimensional charges and distributions}\label{ssec:finite charges and distributions}
A \emph{collection of finite-dimensional charges} is a collection~\(\dist_{\noarg}\coloneqq\pr{\dist_u}_{u\in\setoftseq}\) such that for all \(u\in \setoftseq\), \(\dist_u\) is a probability charge\footnote{
	\label{foot:charge}
	We use the terminology introduced in \cite[Definition~2.1.1]{1983RaoRao-Charges}: given an algebra~\(\algebra\) of events in~\(\posssp\), a probability charge~\(\prob\) is a real-valued map on~\(\algebra\) that is non-negative and finitely additive, with \(\prob\pr{\posssp}=1\).
	See also \cite[Definition~A.5]{1992Breiman-Probability}, \cite[Chapter~II, Section~4]{1994Rogers-Diffusions}, \cite[Section~7.2]{1997Fristedt-Modern}, \cite[Definition~1.15]{2014Troffaes-Lower} or \cite[Definition~A.1.2]{2015Cohen-Stochastic}.
} on~\(\wp\pr{\stsp_u}\).
Usually, we will assume that for all \(u\in\setoftseq\), the corresponding charge~\(\dist_u\) is countably additive, making it a \emph{distribution}\footnote{
	We call a probability measure on the powerset of a countable possibility space a distribution.
	Of course, a distribution~\(\dist_u\) on~\(\wp\pr{\stsp_u}\) is completely determined by the values it assumes on the atoms, so by the unique corresponding probability mass function
	\begin{equation*}
		\stsp_u\to\reals
		\colon x\mapsto \dist_u\pr{\st{x}}.
	\end{equation*}
}; whenever this is the case, we speak of a \emph{collection of finite-dimensional distributions} instead of a collection of finite-dimensional charges.

A collection of finite-dimensional charges -- or distributions -- \(\dist_{\noarg}\) is said to be \emph{consistent} -- or alternatively, `satisfies the compatibility condition' or `has the projective property', see \cite[Chapter~II, Eqn.~(29.6)]{1994Rogers-Diffusions} or \cite[Eqn.~(1.9)]{1980Iosifescu-Finite} -- if for all \(u,v\in\setoftseq\) such that \(u\sqsubseteq v\),
\begin{equation*}
	\dist_u\pr{A}
	= \dist_v\pr[\big]{\st{x_v\in\stsp_v\colon x_u\in A}}
	\quad\text{for all } A\in\wp\pr{\stsp_u}.
\end{equation*}

For the Poisson process, with \(\stsp=\nnegints\) and \(\timeaxis=\nnegreals=\timepoints\), the finite-dimensional distributions are derived from the Poisson distribution.
As explained in, for example, \cite[Chapter~II, Section~33]{1994Rogers-Diffusions}, \cite[Section~23]{1995Billingsley-Probability} or \cite[Chapter~12]{2002Kallenberg-Foundations}, one fixes a rate~\(\lambda\in\nnegreals\), and for all \(u=\pr{t_1, \dots, t_n}\in\setoftseq\) and \(x_u\in\stsp_u\), lets
\begin{equation*}
	\dist_u\pr{x_u}
	\coloneqq \begin{dcases}
		\prod_{k=2}^n e^{-\lambda\pr{t_k-t_{k-1}}} \frac{\pr[\big]{\lambda\pr{t_k-t_{k-1}}}^{\pr{x_{t_k}-x_{t_{k-1}}}}}{\pr{x_{t_k}-x_{t_{k-1}}}!} &\text{if } x_{t_1}\leq \cdots\leq x_{t_n}, \\
		0 &\text{otherwise.}
	\end{dcases}
\end{equation*}
In finite-state and countable-state Markov processes, where one usually also takes~\(\timeaxis=\nnegreals=\timepoints\), the collection of finite-dimensional distributions is constructed with the help of an initial distribution and a semi-group of transition matrices -- corresponding to the matrix exponential of a (transition) rate matrix -- as explained in \cite[Part~II]{1960Chung-Markov}, \cite[Chapter~8]{1980Iosifescu-Finite}, \cite[Chapter~III]{1994Rogers-Diffusions} or \cite[Chapters~2 and 3]{1997Norris-Markov}.

\subsection{From finite-dimensional charges to charges on the cylinder events}
The consistency condition allows us to construct a probability charge on the cylinder events~\(\gencylevts\) from a collection of finite-dimensional charges, which is the first step in establishing the Daniell-Kolomogorov Extension Theorem.
\begin{proposition}
\label{prop:Precursor of Kolmogorov}
	Consider a collection of finite-dimensional charges~\(\dist_{\noarg}\).
	Then there is a unique probability charge~\(\prob\) on the algebra of cylinder events~\(\gencylevts\) such that
	\begin{equation*}
		\prob\pr{\ev{\genX_u\in A}}
		= \dist_u\pr{A}
		\quad\text{for all } u\in\setoftseq, A\in\wp\pr{\stsp_u}
	\end{equation*}
	if and only if~\(\dist_{\noarg}\) is consistent.
\end{proposition}
\begin{proof}
	To see that the consistency of \(\dist_{\noarg}\) is necessary, assume \emph{ex absurdo} that there is a unique probability charge~\(\prob\) on~\(\gencylevts\) with the property in the statement and that \(\dist_{\noarg}\) is not consistent.
	Then there are \(u,v\in\setoftseq\) such that \(u\sqsubseteq v\) and \(A\in\wp\pr{\stsp_u}\) such that
	\begin{equation*}
		\dist_u\pr{A}
		\neq \dist_v\pr{A'}
		\quad\text{with }
		A'
		\coloneqq \st{x_v\in\stsp_v\colon x_u\in A}.
	\end{equation*}
	Now \(\ev{\genX_u\in A}=\ev{\genX_v\in A'}\), and therefore
	\begin{equation*}
		\dist_u\pr{A}
		= \prob\pr{\ev{\genX_u\in A}}
		= \prob\pr{\ev{\genX_v\in A'}}
		= \dist_v\pr{A'},
	\end{equation*}
	which is a clear contradiction.

	Next, we assume that \(\dist_\noarg\) is consistent, and show that there is a unique probability charge~\(\prob\) with the required property.
	We repeat the argument in the `Start of Proof' of Theorem~30.1 in \cite[Chapter~II]{1994Rogers-Diffusions} in our slightly different setting.

	First, fix some \(u,v\in\setoftseq\), \(A\in\wp\pr{\stsp_u}\) and \(B\in\wp\pr{\stsp_v}\).
	Let \(w\in\setoftseq\) be the unique tuple of time points that consists of all time points in~\(u\) and~\(v\).
	Then by construction, \(u\sqsubseteq w\) and \(v\sqsubseteq w\).
	Hence, we can let
	\begin{equation*}
		A^\star
		\coloneqq\st{x_w\in \stsp_w\colon x_u\in A}
		\quad\text{and}\quad
		B^\star
		\coloneqq\st{x_w\in \stsp_w\colon x_u\in B}.
	\end{equation*}
	If \(A^\star=B^\star\), then it follows from the consistency of \(\dist_{\noarg}\) that
	\begin{equation*}
		\dist_u\pr{A}
		= \dist_w\pr{A^\star}
		= \dist_w\pr{B^\star}
		= \dist_v\pr{B}.
	\end{equation*}
	Alternatively, if \(A^\star\cap B^\star=\emptyset\), then it follows from the consistency condition that
	\begin{equation}\label{eq:prop:Precursor of Kolmogorov}
		\dist_w\pr{A^\star\cup B^\star}
		= \dist_w\pr{A^\star}+\dist_w\pr{B^\star}
		= \dist_u\pr{A}+\dist_v\pr{B}.
	\end{equation}

	Now consider two events~\(\tilde{A}, \tilde{B}\) in~\(\gencylevts\).
	Then there are some \(u,v\in\setoftseq\), \(A\in\wp\pr{\stsp_u}\) and \(B\in\wp\pr{\stsp_v}\) such that \(\tilde{A}=\ev{\genX_u\in A}\) and \(\tilde{B}=\ev{\genX_v\in B}\).
	Let \(w\), \(A^\star\) and \(B^\star\) be as before.
	Then \(\ev{\genX_u\in A}=\ev{\genX_w\in A^\star}\) and \(\ev{\genX_v\in B}=\ev{\genX_w\in B^\star}\) by construction.
	If \(\tilde{A}=\tilde{B}\), then \(\ev{\genX_w\in A^\star}=\ev{\genX_w\in B^\star}\), and it then follows from this and Eqn.~\eqref{eqn:condition on set of paths} that \(A^\star=B^\star\).
	If \(\tilde{A}\cap\tilde{B}=\emptyset\), then \(\tilde{A}\cup\tilde{B}=\ev{\genX_w\in A^\star\cup B^\star}\) and \(\tilde{A}\cap\tilde{B}=\ev{\genX_w\in A^\star\cap B^\star}=\emptyset\), and it follows from the latter and Eqn.~\eqref{eqn:condition on set of paths} that \(A^\star\cap B^\star=\emptyset\).

	Due to the preceding two observations and because \(\dist_{{\noarg}}\) is a consistent collection of finite-dimensional charges, we can define the real-valued map~\(\prob\) on~\(\gencylevts\) for all \(\tilde{A}=\ev{\genX_u\in A}\in\gencylevts\) by
	\begin{equation*}
		\prob\pr{\tilde{A}}
		= \prob\pr{\ev{\genX_u\in A}}
		\coloneqq \dist_u\pr{A}.
	\end{equation*}
	Indeed, this map is well-defined because it does not depend on the particular choice of $u$ and $A$. To see this, it suffices to consider the case $\tilde{A}=\tilde{B}$ above, for which we had found that $A^\star=B^\star$ and therefore, that $\dist_u(A)=\dist_v(B)$.
	Furthermore, we trivially have that \(\prob\) is non-negative and that \(\prob\pr{\genpaths}=\prob\pr{\ev{\genX_u\in\stsp_u}}=\dist_u\pr{\stsp_u}=1\).
	To see that $P$ is finitely additive, we consider the case above with $\tilde{A}\cap\tilde{B}=\emptyset$, for which we found that \(\tilde{A}\cup\tilde{B}=\ev{\genX_w\in A^\star\cup B^\star}\) and $A^\star\cap B^\star=\emptyset$.
	It therefore follows from Eqn.~\eqref{eq:prop:Precursor of Kolmogorov} that
	\begin{equation*}
P(\tilde{A}\cup\tilde{B})=P(\ev{\genX_w\in A^\star\cup B^\star})=\dist_w\pr{A^\star\cup B^\star}
		= \dist_u\pr{A}+\dist_v\pr{B}=P(\tilde{A})+P(\tilde{B}).
	\end{equation*}
Hence, \(\prob\) is a probability charge.
	The uniqueness of~\(\prob\) is obvious.
\end{proof}

Proposition~\ref{prop:Precursor of Kolmogorov} allows us to go from a consistent collection of finite-dimensional charges to a probability charge, but we can also go the other way around: any (consistent) triple~\(\pr{\genpaths, \genalgebra,\prob}\)  induces a collection of finite-dimensional charges that is consistent.
\begin{definition}
\label{def:fidis from charge}
	Consider a subset~\(\genpaths\) of~\(\stsp^{\timeaxis}\), an algebra~\(\genalgebra\supseteq\gencylevts\) of events in~\(\genpaths\) and a probability charge~\(\prob\) on~\(\genalgebra\).
	Then for all \(u\in\setoftseq\), the map
	\begin{equation*}
		\dist_u
		\colon \wp\pr{\stsp_u} \to\cci{0}{1}
		\colon A\mapsto \dist_u\pr{A}\coloneqq\prob\pr{\ev{\genX_u\in A}}
	\end{equation*}
	is a probability charge.
	Moreover, the collection~\(\pr{\dist_u}_{u\in\setoftseq}\) is consistent, which is why we call it the \emph{(collection of) finite-dimensional charges of~\(\prob\)}.
	If \(\prob\) is countably additive, then \(\dist_u\) is a distribution for all \(u\in\setoftseq\); in this case, we call \(\pr{\dist_u}_{u\in\setoftseq}\) the \emph{(collection of) finite-dimensional distributions of~\(\prob\)}.
\end{definition}
\begin{proof}
	Follows immediately from the properties of (countably additive) probability charges.
\end{proof}

\subsection{From the finite-dimensional distributions to a stochastic process}
Proposition~\ref{prop:Precursor of Kolmogorov} and Definition~\ref{def:fidis from charge} show that for a fixed set of paths~\(\genpaths\) that satisfies Eqn.~\eqref{eqn:condition on set of paths}, there is a one to one correspondence between consistent collections of finite-dimensional charges and probability charges on~\(\gencylevts\).
It also follows from Definition~\ref{def:fidis from charge} that any stochastic process~\(\pr{\genpaths, \genalgebra, \prob}\) induces a consistent collection of finite-dimensional distributions.

This raises the following question.
Given a set of paths~\(\genpaths\subseteq\stsp^{\timeaxis}\) and a consistent collection of finite-dimensional distributions~\(\dist_\noarg\), is there a probability measure~\(\prob\) on some \(\sigma\)-algebra~\(\genalgebra\supseteq\gencylevts\) in~\(\genpaths\) such that the finite-dimensional distributions of~\(\prob\) are~\(\dist_{\noarg}\)?
The key to answering this question lies in Caratheodory's Extension Theorem -- see for example \cite[Theorem~A.9]{1992Breiman-Probability}, \cite[Chapter~II, Theorem~5.1]{1994Rogers-Diffusions}, \cite[Chapter~7, Theorem~14]{1997Fristedt-Modern} or \cite[Theorem~A.1.17]{2015Cohen-Stochastic} -- which says that any probability charge~\(\prob\) on some algebra~\(\genalgebra\) can be extended to a probability measure~\(\prob_{\sigma}\) on the generated \(\sigma\)-algebra~\(\sigma\pr{\genalgebra}\) if and only if \(\prob\) is countably additive, and that this extension is then unique.
Indeed, since the consistent collection of finite-dimensional distributions~\(\dist_{\noarg}\) corresponds to a unique probability charge~\(\prob\) on~\(\gencylevts\), it then follows that there is a probability measure on some \(\sigma\)-algebra~\(\genalgebra\supseteq\gencylevts\) if and only if \(\prob\) is countably additive.

In the remainder, we investigate the countable additivity of the induced probability charge~\(\prob\) for two important sets of paths: the set of all paths in Section~\ref{sec:set of all paths} and the set of all `càdlàg' paths in Section~\ref{sec:cadlag paths}.

\section{The set of all paths}
\label{sec:set of all paths}
We denote the set of all paths~\(\stsp^{\timeaxis}\) by \(\allpaths\); clearly, \(\allpaths\) satisfies Eqn.~\eqref{eqn:condition on set of paths}.
For all~\(t\in \timeaxis\) and \(u\in\setoftseq\), we denote the corresponding projector variables for~\(\allpaths\) by~\(\allX_t\) and \(\allX_u\).
We do something similar for the cylinder events: for any tuple of time points~\(u\in \setoftseq\) and any subset~\(A\) of~\(\stsp_u\), we denote the corresponding cylinder event by
\begin{equation*}
	\ev{\allX_u\in A}
	\coloneqq\st[\big]{\allpth\in\allpaths\colon \allpth\pr{u}\in A}.
\end{equation*}
Hence, the algebra of cylinder events for the set of all paths~\(\allpaths\) is
\begin{equation*}
	\allcylevts
	\coloneqq \st[\big]{\ev{\allX_u\in A}\colon u\in\setoftseq, A\in\wp\pr{\stsp_u}}.
\end{equation*}

\subsection{Establishing countable additivity}

Crucially, the set of all paths~\(\allpaths\) allows us to establish the countable additivity of any probability charge~\(\prob\) induced by a consistent collection of finite-dimensional distributions.
Since this result is essentially well-known, at least in case \(\timepoints=\timeaxis\), we have relegated our proof to Appendix~\ref{asec:a proof for Kolmogorov extension theorem} further on.
\begin{theorem}
\label{the:precursor for all paths}
	For any consistent collection~\(\dist_{\noarg}\) of finite-dimensional distributions, the corresponding probability charge \(\prob\) on~\(\allcylevts\) of Proposition~\ref{prop:Precursor of Kolmogorov} is countably additive.
\end{theorem}

As we explained at the end of Section~\ref{sec:constructing a stochastic process}, Theorem~\ref{the:precursor for all paths} admits us to invoke Caratheodory's Extension Theorem to obtain a stochastic process~\(\pr{\allpaths,\sigma\pr{\allcylevts},\prob}\) such that the finite-dimensional distributions of~\(\prob\) are the prescribed ones~\(\dist_{\noarg}\).
This result is known as the Daniell--Kolmogorov Extension Theorem, and is similar to -- and essentially implied by -- Theorem~31.1 in \cite[Chapter~II]{1994Rogers-Diffusions}, or in case \(\stsp\) is finite, the result outlined in \cite[Section~1.5.2]{1980Iosifescu-Finite}.
\begin{theorem}
\label{the:Kolomogorov extension theorem}
	For any consistent collection of finite-dimensional distributions~\(\dist_{\noarg}\), there is a unique probability measure~\(\prob\) on~\(\sigma\pr{\allcylevts}\) such that
	\begin{equation*}
		\prob\pr{\ev{\allX_u\in A}}
		= \dist_u\pr{A}
		\quad\text{for all } u\in\setoftseq, A\in\wp\pr{\stsp_u}.
	\end{equation*}
\end{theorem}
\begin{proof}
	Follows immediately from Proposition~\ref{prop:Precursor of Kolmogorov}, Theorem~\ref{the:precursor for all paths} and Caratheodory's Extension Theorem.
\end{proof}

\subsection{The insufficiency of the product sigma-algebra}
\label{ssec:insufficiency}
There is one major issue with Theorem~\ref{the:Kolomogorov extension theorem}: the generated \(\sigma\)-algebra \(\sigma\pr{\allcylevts}\) is not sufficiently rich, in the sense that many practically-relevant events do not belong to it.
Simply put it comes down to this: the events in~\(\sigma\pr{\allcylevts}\) only depend on the state of the paths in the time points in a countable subset of~\(\timepoints\).
We are by no means the first to signal this issue: Breiman~\cite[Proposition~12.8]{1992Breiman-Probability}, Billingsley~\cite[Theorem~36.3]{1995Billingsley-Probability} and Cohen \& Elliot~\cite[Lemma~A.2.2]{2015Cohen-Stochastic} all mention it, to give but three examples.
We give a (slightly modified version of) Lemma~25.9 in \cite[Chapter~II]{1994Rogers-Diffusions}, applicable to the general case in Section~\ref{sec:constructing a stochastic process}.
In it, we use the following hitherto undefined notation: for all subsets~\(\altpoints, \mathcal{R}\) of~\(\timeaxis\) such that \(\altpoints\subseteq\mathcal{R}\), we denote the restriction of \(\vpth\colon\mathcal{R}\to\stsp\) to~\(\altpoints\) by~\(\vpth\vert_{\altpoints}\), and we let
\begin{equation*}
	\tilde{\genalgebra}_{\altpoints}
	\coloneqq \sigma\pr[\big]{\st[\big]{\st{\vpth\in\stsp^\altpoints\colon \vpth\pr{s}\in A}\colon s\in\altpoints, A\in\wp\pr{\stsp}}}.
\end{equation*}

\begin{lemma}
\label{lem:generated sigma field is that of the sigma cylinders}
	Consider a non-empty subset \(\genpaths\) of~\(\stsp^{\timeaxis}\) and let \(\mathfrak{C}\) be the set of all non-empty countable subsets of~\(\timepoints\).
	Then
	\begin{equation*}
		\sigma\pr{\gencylevts}
		= \bigcup_{\countpoints\in\mathfrak{C}} \st[\big]{\st{\pth\in\genpaths\colon \pth\vert_{\countpoints}\in \tilde{A}_{\countpoints}}\colon \tilde{A}_{\countpoints}\in\tilde{\genalgebra}_{\countpoints}}.
	\end{equation*}
\end{lemma}
In our proof, we make use of the following claim (without formal proof) of Rogers \& Williams~\cite[Chapter~II, Eqn.~(25.3)]{1994Rogers-Diffusions}.
\begin{lemma}
\label{lem:restriction is measurable}
	For any two subsets~\(\altpoints,\mathcal{R}\) of~\(\timeaxis\) such that \(\altpoints\subseteq\mathcal{R}\), the restriction operator~\(\noarg\vert_{\altpoints}\colon\stsp^{\mathcal{R}}\to\stsp^{\altpoints}\) is \(\tilde{\genalgebra}_{\mathcal{R}}/\tilde{\genalgebra}_{\altpoints}\)-measurable.
\end{lemma}
\begin{proof}
	By a standard result in measure theory -- see for example Proposition~2.3 in~\cite[Chapter~II]{1994Rogers-Diffusions} -- we need to look only at the sets~\(\st{\rpth\in\stsp^{\altpoints}\colon \rpth\pr{s}\in A}\) that generate~\(\tilde{\genalgebra}_{\altpoints}\); more formally, this is the case if (and only if) for all \(s\in\altpoints\) and \(A\in\wp\pr{\stsp}\),
	\begin{equation*}
		\st[\big]{\vpth\in\stsp^{\mathcal{R}}\colon \vpth\vert_{\altpoints}\in \st{\rpth\in\stsp^{\altpoints}\colon \rpth\pr{s}\in A}}
		\in \tilde{\genalgebra}_{\mathcal{R}},
	\end{equation*}
	and this condition is trivially satisfied because
	\begin{equation*}
		\st[\big]{\vpth\in\stsp^{\mathcal{R}}\colon \vpth\vert_{\altpoints}\in \st{\rpth\in\stsp^{\altpoints}\colon \rpth\pr{s}\in A}}
		= \st[\big]{\vpth\in\stsp^{\mathcal{R}}\colon \vpth\pr{s}\in A}
		\in \tilde{\genalgebra}_{\mathcal{R}}.
	\end{equation*}
\end{proof}
\begin{proof}[Proof of Lemma~\ref{lem:generated sigma field is that of the sigma cylinders}]
	We adapt the argument in \cite[Proof of Lemma~25.9]{1994Rogers-Diffusions} to our slightly different setting.
	Since \(\gencylevts=\st{\dot{A}\cap\genpaths\colon \dot{A}\in\allcylevts}\) by definition of~\(\gencylevts\) and \(\allcylevts\), it follows from a standard result in measure theory -- see for example Theorem~10.1 in~\cite{1995Billingsley-Probability} -- that
	\begin{equation*}
		\sigma\pr{\gencylevts}
		= \st{\dot{A}\cap\genpaths\colon \dot{A}\in\sigma\pr{\allcylevts}}.
	\end{equation*}
	On the other hand, it is clear that for all \(\countpoints\in\mathfrak{C}\),
	\begin{equation*}
		\st[\big]{\st{\pth\in\genpaths\colon \pth\vert_{\countpoints}\in \tilde{A}_{\countpoints}}\colon \tilde{A}_{\countpoints}\in\tilde{\genalgebra}_{\countpoints}}
		= \st[\big]{\st{\allpth\in\allpaths\colon \allpth\vert_{\countpoints}\in \tilde{A}_{\countpoints}}\cap\genpaths\colon \tilde{A}_{\countpoints}\in\tilde{\genalgebra}_{\countpoints}}.
	\end{equation*}
	Hence, it suffices to prove that
	\begin{equation*}
		\sigma\pr{\allcylevts}
		= \bigcup_{\countpoints\in\mathfrak{C}} \st[\big]{\st{\allpth\in\allpaths\colon \allpth\vert_{\countpoints}\in \tilde{A}_{\countpoints}}\colon \tilde{A}_{\countpoints}\in\tilde{\genalgebra}_{\countpoints}},
	\end{equation*}
	which is precisely the equality in the statement in case \(\genpaths=\allpaths\).
	To simplify our notation, we let
	\begin{equation*}
		\tilde{\genalgebra}^{\uparrow\allpaths}_{\countpoints}
		\coloneqq \st[\big]{\st{\allpth\in\allpaths\colon \allpth\vert_{\countpoints}\in\tilde{A}_{\countpoints}}\colon \tilde{A}_{\countpoints}\in\tilde{\genalgebra}_{\countpoints}}
		\quad\text{for all } \countpoints\in\mathfrak{C}.
	\end{equation*}

	First, we verify that \(\genalgebra\coloneqq\bigcup_{\countpoints\in\mathfrak{C}}\tilde{\genalgebra}^{\uparrow\allpaths}_{\countpoints}\) is a \(\sigma\)-algebra.
	It is clear that \(\genalgebra\) includes the empty set and is closed under taking complements, so we really only need to show that \(\genalgebra\) is closed under countable unions.
	To this end, we fix any sequence~\(\pr{\dot{A}_n}_{n\in\nats}\) in~\(\genalgebra\).
	Then by definition, for all \(n\in\nats\) there are some \(\countpoints_n\in\mathfrak{C}\) and \(\tilde{A}_{\countpoints_n}\in\tilde{\genalgebra}_{\countpoints_n}\) such that
	\begin{equation*}
		\dot{A}_n
		= \st{\allpth\in\allpaths\colon \allpth\vert_{\countpoints_n}\in\tilde{A}_{\countpoints_n}}.
	\end{equation*}
	Then clearly, the set~\(\countpoints\coloneqq\bigcup_{n\in\nats} \countpoints_n\) is a countable subset of~\(\timepoints\), and therefore~\(\countpoints\in\mathfrak{C}\).
	Furthermore, for all \(n\in\nats\),
	\begin{equation*}
		\tilde{A}^n_{\countpoints}
		\coloneqq \st{\vpth\in\stsp^{\countpoints}\colon \vpth\vert_{\countpoints_n}\in\tilde{A}_{\countpoints_n}}
		\in \tilde{\genalgebra}_{\countpoints}
	\end{equation*}
	because the restriction operator~\(\noarg\vert_{\countpoints_n}\colon\stsp^{\countpoints}\to\stsp^{\countpoints_n}\) is \(\tilde{\genalgebra}_{\countpoints}/\tilde{\genalgebra}_{\countpoints_n}\)-measurable due to Lemma~\ref{lem:restriction is measurable}, and then
	\begin{equation*}
		\dot{A}_n
		= \st{\allpth\in\allpaths\colon\allpth\vert_{\countpoints_n}\in\tilde{A}_{\countpoints_n}}
		= \st{\allpth\in\allpaths\colon\allpth\vert_{\countpoints}\in\tilde{A}^n_{\countpoints}}.
	\end{equation*}
	Since \(\bigcup_{n\in\nats}\tilde{A}^n_{\countpoints}\in\tilde{\genalgebra}_{\countpoints}\) because \(\tilde{\genalgebra}_{\countpoints}\) is a \(\sigma\)-algebra, it follows from this that
	\begin{equation*}
		\bigcup_{n\in\nats} \dot{A}_n
		= \bigcup_{n\in\nats} \st{\allpth\in\allpaths\colon\allpth\vert_{\countpoints}\in\tilde{A}^n_{\countpoints}}
		= \st*{\allpth\in\allpaths\colon \allpth\vert_{\countpoints}\in \bigcup_{n\in\nats}\tilde{A}^n_{\countpoints}}
		\in\genalgebra,
	\end{equation*}
	as required.

	Next, we show that \(\sigma\pr{\allcylevts}\subseteq\genalgebra\).
	Since \(\genalgebra\) is a \(\sigma\)-algebra, it suffices to show that \(\allcylevts\subseteq\genalgebra\).
	To this end, we fix any \(\dot{A}=\ev{\allX_u\in A}\in\allcylevts\).
	If we enumerate the time points in~\(u\) as \(\pr{t_1, \dots, t_n}\) and let \(\countpoints\coloneqq\st{t_1, \dots, t_n}\), then clearly
	\begin{equation*}
		\tilde{A}_{\countpoints}
		\coloneqq \st{\vpth\in\stsp^{\countpoints}\colon \pr{\vpth\pr{t_1}, \dots, \vpth\pr{t_n}}\in A}
		= \bigcup_{x_u\in A} \bigcap_{k=1}^n \st{\vpth\in\stsp^{\countpoints}\colon \vpth\pr{t_k}=x_{t_k}}
		\in \tilde{\genalgebra}_{\countpoints}
	\end{equation*}
	because \(A\subseteq\stsp_u\) is countable.
	Hence,
	\begin{equation*}
		\dot{A}
		= \st{\allX_u\in A}
		= \st{\allpth\in\allpaths\colon \allpth\vert_{\countpoints}\in\tilde{A}_{\countpoints}}
		\in \genalgebra,
	\end{equation*}
	as required.

	To obtain the equality that we are after, it remains for us to verify that \(\genalgebra\subseteq\sigma\pr{\allcylevts}\), or equivalently, that $\st{\allpth\in\allpaths\colon \allpth\vert_{\countpoints}\in \tilde{A}_{\countpoints}}
		\in \sigma\pr{\allcylevts}$ for all \(\countpoints\in\mathfrak{C}\) and $\tilde{A}_{\countpoints}\in\tilde{\genalgebra}_{\countpoints}$.
	So let us fix any \(\countpoints\in\mathfrak{C}\) and $\tilde{A}_{\countpoints}\in\tilde{\genalgebra}_{\countpoints}$. Then
\begin{equation*}
\st{\allpth\in\allpaths\colon \allpth\vert_{\countpoints}\in \tilde{A}_{\countpoints}}=\st{\allpth\in\stsp^{\timeaxis}\colon \allpth\vert_{\countpoints}\in \tilde{A}_{\countpoints}}\in\tilde{\genalgebra}_{\timeaxis}
\end{equation*}
because the restriction operator~\(\noarg\vert_{\countpoints}\colon\stsp^{\timeaxis}\to\stsp^{\countpoints}\) is \(\tilde{\genalgebra}_{\timeaxis}/\tilde{\genalgebra}_{\countpoints}\)-measurable due to Lemma~\ref{lem:restriction is measurable}. Since
\begin{align*}
\tilde{\genalgebra}_{\timeaxis}
	= \sigma\pr[\big]{\st[\big]{\st{\allpth\in\stsp^\timeaxis\colon \allpth\pr{s}\in A}\colon s\in\timeaxis, A\in\wp\pr{\stsp}}}
\end{align*}
and, for all $s\in\timeaxis$ and $A\in\wp\pr{\stsp}$,
\begin{equation*}
\{\allpth\in{\stsp}^{\timeaxis}\colon\allpth(s)\in A\}
=\{\allpth\in\allpaths\colon\allpth(s)\in A\}=\{\allX_{(s)}\in A\},
\end{equation*}
we furthermore have that
\begin{align*}
\tilde{\genalgebra}_{\timeaxis}
	\subseteq \sigma\pr[\big]{\st[\big]{\ev{\allX_u\in A}\colon u\in\setoftseq, A\in\wp\pr{\stsp_u}}}=\sigma\pr{\allcylevts}.
\end{align*}
It therefore follows that $\st{\allpth\in\allpaths\colon \allpth\vert_{\countpoints}\in \tilde{A}_{\countpoints}}
		\in\tilde{\genalgebra}_{\timeaxis}\subseteq\sigma\pr{\allcylevts}$, as required.
\end{proof}

In the particular case where $\paths$ is the set $\allpaths$ of all paths, and hence $\gencylevts=\allcylevts$, the main take-away from Lemma~\ref{lem:generated sigma field is that of the sigma cylinders} is that if \(\timepoints\) is not countable, then \(\sigma\pr{\allcylevts}\) may \emph{not} contain `all' events that are of interested to us. The following example illustrates this.
\begin{example}
\label{ex:level set of hitting time for all}
	Let \(\timeaxis\coloneqq\nnegreals\eqqcolon\timepoints\).
	Suppose we are interested in the event that our system is in some given state~\(x\in\stsp\) before some time~\(T\in\posreals\):
	\begin{equation*}
		\dot{A}
		\coloneqq \st{\allpth\in\allpaths\colon \pr{\exists t\in\cci{0}{T}}~\allpth\pr{t}=x}.
	\end{equation*}
	Events like this are common in applications, for example in model checking \cite{2008Baier-Princliples}.
	Obviously, we can write this event as an uncountable union of events in~\(\allcylevts\):
	\begin{equation*}
		\dot{A}
		= \bigcup_{t\in\cci{0}{T}} \ev{\allX_t=x}.
	\end{equation*}
	However, the generated \(\sigma\)-algebra \(\sigma\pr{\allcylevts}\) is only guaranteed to be closed under countable unions.
	In fact, we can use Lemma~\ref{lem:generated sigma field is that of the sigma cylinders} to prove that \(\dot{A}\) does not belong to~\(\sigma\pr{\allcylevts}\).
	Assume \emph{ex absurdo} that \(\dot{A}\in\sigma\pr{\allcylevts}\).
	Then by Lemma~\ref{lem:generated sigma field is that of the sigma cylinders}, there are some non-empty countable subset~\(\countpoints\) of~\(\timepoints\) and subset~\(\tilde{A}\) of~\(\stsp^{\countpoints}\) such that
	\begin{equation*}
		\dot{A}
		= \st{\allpth\in\allpaths\colon \allpth\vert_{\countpoints}\in \tilde{A}}.
	\end{equation*}
	Take any path~\(\allpth\in \dot{A}^{\compl}\).
	Then clearly, \(\allpth\pr{t}\neq x\) for all \(t\in\cci{0}{T}\).
	Now take any \(t^\star\in\cci{0}{T}\setminus\countpoints\), and let \(\tilde{\pth}\in\allpaths\) be the path defined for all \(t\in\timepoints\) by \(\tilde{\pth}\pr{t}\coloneqq x\) if \(t=t^\star\) and \(\tilde{\pth}\pr{t}=\allpth\pr{t}\) otherwise.
	Then on the one hand, \(\tilde{\pth}\in \dot{A}\) by construction.
	On the other hand, \(\tilde{\pth}\vert_{\countpoints}=\allpth\vert_{\countpoints}\) by construction, and since \(\allpth\notin \dot{A}\) this implies that \(\tilde{\pth}\notin \dot{A}\), a clear contradiction.
\end{example}

The reason why this problem occurs is because, whenever \(\timepoints\) is uncountable, (the restriction to~\(\timepoints\) of) a path~\(\pth\) in~\(\allpaths\) is not fully defined by the states it assumes on a countable subset of~\(\timepoints\).
The obvious solution to this problem is therefore to focus on a subset~\(\genpaths\) of the set of all paths~\(\allpaths\) that does not contain such paths, and it is standard to choose the set of càdlàg paths for this; see Section~\ref{sec:cadlag paths} further on.

The usual manner to proceed is then to construct a `modification'~\(\tilde{X}_{\noarg}\) of \(\allX_{\noarg}\) that has `càdlàg sample paths', which can only be done under some conditions on the finite-dimensional distributions; we are only aware of work where this is done for Markov processes -- that is, only for finite-dimensional distributions of a specific form.
Most authors then proceed with this modified process~\(\tilde{X}_{\noarg}\) as if it was the original process~\(\allX_{\noarg}\), but this modification is a drastically changed version of~\(\allX_{\noarg}\).
Why one would be allowed to even change the outcomes (or sample paths) of a stochastic process is often not really given much thought or motivation, expect perhaps for the (implicit) pragmatic justification that it `just works'.
That said, one could use the modified process~\(\tilde{X}_{\noarg}\) to obtain a probability measure on a (sufficiently rich) \(\sigma\)-algebra of events in the set of càdlàg paths, but most -- if not almost all -- works that we are aware of skip this crucial step.
We have found two notable exceptions: Cramér \& Leadbetter~\cite[Sections~3.2, 3.3 and 3.6]{1967Cramer-Stationary} more or less explain how this can be done, as does Borovkov~\cite[Section~18.1]{2013Borovkov-Probability} in less detail.

In contrast, we now set out to provide a necessary and sufficient regularity condition on any (so not necessarily Markov) collection of finite-dimensional distributions for there to be a corresponding probability measure on the \(\sigma\)-algebra of events generated by the cylinder events in the set of \emph{càdlàg} paths -- instead of the usual product \(\sigma\)-algebra~\(\sigma\pr{\allcylevts}\) of events in the set of \emph{all} paths~\(\allpaths\).
The crucial benefit of our approach is that we do \emph{not} need to modify our stochastic process to avoid the issue illustrated in Example~\ref{ex:level set of hitting time for all} -- that is, the lack of richness of the generated \(\sigma\)-algebra

\section{The set of càdlàg paths}
\label{sec:cadlag paths}
Before we can introduce the càdlàg paths, we need to deal with some topological considerations.
Consider a subset~\(\altpoints\) of~\(\reals\).
A real number~\(s\in \reals\) is an \emph{isolated point} of~\(\altpoints\) if there is some~\(\delta\in \posreals\) such that \(\ooi{s-\delta}{s+\delta}\cap\altpoints=\st{s}\), a \emph{right-sided limit point} of~\(\altpoints\) if \(\ooi{s}{s+\delta}\cap\altpoints\neq\emptyset\) for all \(\delta\in \posreals\) and a \emph{left-sided limit point} of~\(\altpoints\) if \(\ooi{s-\delta}{s}\cap\altpoints\neq\emptyset\) for all \(\delta\in \posreals\).
We collect all isolated points, right-sided limit points and left-sided limit points of~\(\altpoints\) in \(\isos\pr{\altpoints}\), \(\rlims\pr{\altpoints}\) and \(\llims\pr{\altpoints}\), respectively, and let \(\lims\pr{\altpoints}\coloneqq\rlims\pr{\altpoints}\cup\llims\pr{\altpoints}\); obviously, \(\isos\pr{\altpoints}\cap\lims\pr{\altpoints}=\emptyset\).
The union of the isolated points and the right-sided limits points of~\(\altpoints\) is its closure with respect to the \emph{lower limit topology} -- also known as the `right half-open topology' -- on~\(\reals\) \cite[Section~16.33.e]{1997Schechter-Handbook}, the union of the isolated points and the left-sided limits points of~\(\altpoints\) is its closure with respect to the \emph{upper limit topology} -- also know as the `left half-open interval topology' -- on~\(\reals\) \cite[Chapter~III, Section~3, Ex.~4]{1966Dugundji-Topology}, while the union of the isolated points and the limit points of~\(\altpoints\) is the closure with respect to the usual topology on~\(\reals\).

Next, we take directional limits along~\(\altpoints\).
The following is a specialization of the general (topological) definition of a limit \cite[Appendix~C, Definition~8 and Promblem~18]{1997Fristedt-Modern} to the lower and upper limit topologies, see also \cite[Section~15.21]{1997Schechter-Handbook}.
Consider a map~\(\phi\) from \(\altpoints\) to a metric space~\(\pr{\mathbb{M}, \rho}\); we will need two particular cases: (i) the set of real numbers~\(\reals\) with the usual metric induced by the absolute value (in Definition~\ref{def:regularity} further on), and (ii) the set of states~\(\stsp\) with the Kronecker metric (in Definition~\ref{def:cadlag} further on).
Then for any right-sided limit point~\(t\) of~\(\altpoints\), we say that \emph{the right-sided limit of~\(\phi\) in~\(t\) along~\(\altpoints\) exists} if there is some -- necessarily unique -- \(\ell\in \mathbb{M}\) such that
\begin{equation*}
	\pr{\forall \epsilon\in\posreals}
	\pr[\big]{\exists\delta\in\posreals}
	\pr[\big]{\forall r\in\ooi{t}{t+\delta}\cap\altpoints}~
	\rho\pr{\phi\pr{r}, \ell}
	< \epsilon;
\end{equation*}
whenever this is the case, we denote this limit~\(\ell\) by~\(\lim_{\altpoints\ni r\searrow t}\phi\pr{r}\).
Similarly, for any left-sided limit point~\(t\) of~\(\altpoints\), we say that \emph{the left-sided limit of~\(\phi\) in~\(t\) along~\(\altpoints\) exists} if there is some -- necessarily unique -- \(\ell\in\mathbb{M}\) such that
\begin{equation*}
	\pr[\big]{\forall\epsilon\in\posreals}
	\pr[\big]{\exists\delta\in\posreals}
	\pr[\big]{\forall s\in\ooi{t-\delta}{t}\cap\altpoints}~
	\rho\pr{\phi\pr{s},  \ell}
	< \epsilon,
\end{equation*}
and then we denote this limit~\(\ell\) by~\(\lim_{\altpoints\ni s\nearrow t}\phi\pr{s}\).

With this terminology and notation, we can define càdlàg paths as those paths that have left-sided and right-sided limits and are continuous from the right in all time points~\(t\) in the closure of~\(\timeaxis\) where these concepts make sense.
\begin{definition}
\label{def:cadlag}
Consider a non-empty subset $\mathcal{S}$ of $\reals$.
	A map~\(\vpth\colon\mathcal{S}\to\stsp\) is \emph{càdlàg} if it has a left-sided limit (along~\(\mathcal{S}\)) in all~\(t\in\llims\pr{\mathcal{S}}\), meaning that
	\begin{equation*}
		\pr[\big]{\exists \delta\in\posreals}
		\pr[\big]{\exists x\in\stsp}
		\pr[\big]{\forall s\in\mathcal{S}\cap\ooi{t-\delta}{t}}~
		\vpth\pr{s}
		= x,
	\end{equation*}
	has a right-sided limit (along~\(\mathcal{S}\)) in all~\(t\in\rlims\pr{\mathcal{S}}\), meaning that
	\begin{equation*}
		\pr[\big]{\exists \delta\in\posreals}
		\pr[\big]{\exists x\in\stsp}
		\pr[\big]{\forall r\in\mathcal{S}\cap\ooi{t}{t+\delta}}~
		\vpth\pr{r}
		= x,
	\end{equation*}
	and is right-continuous (along~\(\mathcal{S}\)) in all~\(t\in\mathcal{S}\cap\rlims\pr{\mathcal{S}}\), meaning that \(\lim_{\mathcal{S}\ni r\searrow t} \vpth\pr{r}=\vpth\pr{t}\). Càdlàg paths correspond to the case $\mathcal{S}=\timeaxis$; that is, càdlàg paths are càdlàg maps on $\timeaxis$.
	We collect all such càdlàg paths in~\(\cadpaths\).
\end{definition}

It is straightforward to verify that \(\cadpaths\) satisfies Eqn.~\eqref{eqn:condition on set of paths}.
Another important property of càdlàg paths is that they can always be extended to a càdlàg map on~\(\reals\).
Because this is an intermediary technical result, we have relegated its proof to Appendix~\ref{asec:Properties of cadlag paths}.
\begin{lemma}
\label{lem:extension of cadlag path}
	Any càdlàg map~\(\vpth\colon\mathcal{S}\to\stsp\), and therefore also any càdlàg path $\cadpth\in\cadpaths$, can be extended to a map~\(\rpth\colon\reals\to\stsp\) that is right-continuous and has left-sided limits everywhere.
\end{lemma}

For all~\(t\in\timeaxis\) and \(u\in\setoftseq\), we denote the corresponding projector variable for the set of all càdlàg paths~\(\cadpaths\) by \(\cadX_t\) and \(\cadX_u\), and we use similar notation for the cylinder events: for any tuple of time points~\(u\in \setoftseq\) and any subset~\(A\) of~\(\stsp_u\), we denote the corresponding cylinder event by
\begin{equation*}
	\ev{\cadX_u\in A}
	\coloneqq\st[\big]{\cadpth\in\cadpaths\colon \cadpth\pr{u}\in A}
	= \ev{\allX_u\in A}\cap\cadpaths.
\end{equation*}
Hence, the algebra of cylinder events for the set of all càdlàg paths~\(\cadpaths\) is
\begin{equation*}
	\cadcylevts
	\coloneqq \st[\big]{\ev{\cadX_u\in A}\colon u\in\setoftseq, A\in\wp\pr{\stsp_u}}
	= \st[\big]{\dot{A}\cap\cadpaths\colon \dot{A}\in\allcylevts}.
\end{equation*}

\subsection{The sufficiency of the generated sigma-algebra}

So are càdlàg paths fully defined by the states they assume on a countable subset of~\(\timeaxis\)?
The following result makes clear that this is the case, at least if \(\timeaxis\) has a countable subset~\(\densepoints\) that is dense in~\(\timeaxis\) for the lower limit topology, or equivalently, such that \(\timeaxis\subseteq\densepoints\cup\rlims\pr{\densepoints}\); its simple proof can be found in Appendix~\ref{asec:Properties of cadlag paths}.
\begin{lemma}
\label{lem:cadlag defined on subset}
	Consider a subset~\(\densepoints\) of~\(\timeaxis\) such that \(\timeaxis\subseteq\densepoints\cup\rlims\pr{\densepoints}\).
	Then for all~\(\cadpth_1,\cadpth_2\in\cadpaths\),
	\begin{equation*}
		\cadpth_1=\cadpth_2
		\Leftrightarrow
		\pr{\forall d\in\densepoints}~
		\cadpth_1\pr{d}
		= \cadpth_2\pr{d}.
	\end{equation*}
\end{lemma}
Fortunately, the requirement that \(\timeaxis\) has such a countable dense subset is always met.
\begin{lemma}
\label{lem:subset of sorgenfrey is seperable}
	Any subset~\(\altpoints\) of~\(\reals\) has a countable subset~\(\densepoints\) that is dense with respect to the lower limit topology, meaning that
	\begin{equation*}
		\altpoints
		= \st[\big]{s\in\altpoints\colon \pr{\forall \delta\in\posreals}~\coi{s}{s+\delta}\cap\densepoints\neq\emptyset}
		\subseteq \densepoints\cup\rlims\pr{\densepoints}.
	\end{equation*}
\end{lemma}
\begin{proof}
	According to \cite[Chapter~VIII, Section~7, Ex.~5]{1966Dugundji-Topology}, any subspace~\(\altpoints'\) of the set of real numbers~\(\reals\) equipped with the upper limit topology is separable, in the sense that it has a countable dense subset~\(\densepoints'\) \cite[Chapter~VIII, Section~7, Definition~7.1]{1966Dugundji-Topology}.
	Since the intervals~\(\oci{a}{b}\) form a basis for the upper limit topology \cite[Chapter~III, Section~3, Ex.~4]{1966Dugundji-Topology}, it follows from Theorem~4.13 in~\cite[Chapter~III]{1966Dugundji-Topology} that
	\begin{equation*}
		\pr{\oci{a}{b}\cap\altpoints'}\cap\densepoints'
		\neq \emptyset
		\text{ for all } a,b\in\reals \text{ such that } a<b,
	\end{equation*}
	and therefore
	\begin{equation*}
		\oci{s-\delta}{s}\cap\densepoints'
		\neq \emptyset
		\text{ for all } s\in\altpoints', \delta\in\posreals.
	\end{equation*}
	Due to symmetry, this implies that any subspace $\altpoints$ of the set of real numbers equipped with the lower limit topology has a countable dense subset $\densepoints$ as well, with
	\begin{equation*}
		\coi{s}{s+\delta}\cap\densepoints
		\neq \emptyset
		\text{ for all } s\in\altpoints, \delta\in\posreals.
	\end{equation*}
	This also implies that for any $s\in\altpoints\setminus\densepoints$, $\ooi{s}{s+\delta}\cap\densepoints\neq\emptyset$ for all $\delta\in\posreals$, and therefore $s\in\rlims\pr{\densepoints}$.
\end{proof}

This is already promising, but we have left one big question unanswered: is the \(\sigma\)-algebra~\(\sigma\pr{\cadcylevts}\) generated by the cylinder events sufficiently rich, meaning that it contains `all' events that we are interested in?
Due to Lemmas~\ref{lem:generated sigma field is that of the sigma cylinders}, \ref{lem:cadlag defined on subset} and \ref{lem:subset of sorgenfrey is seperable}, we would expect that the answer to this question is yes, but since it is hard to formalise what precisely `all events of interest' are, we cannot answer this question definitively.
We can, however, easily show that the issue we observed in Example~\ref{ex:level set of hitting time for all} is now resolved.
\begin{example}
	Let us return to the setting of Example~\ref{ex:level set of hitting time for all}.
	Note that the set of positive rationals~\(\mathbb{Q}_{\geq 0}\) is a countable subset of~\(\timeaxis=\reals_{\geq0}=\timepoints\) that is dense with respect to the lower limit topology.
	For the set of càdlàg paths, the event of interest is
	\begin{equation*}
		\mathring{A}
		\coloneqq \bigcup_{t\in\cci{0}{T}} \ev{\cadX_t=x}.
	\end{equation*}
	Fix any càdlàg path~\(\cadpth\in\cadpaths\).
	Then \(\cadpth\in\mathring{A}\) if and only if there is some~\(t\in\cci{0}{T}\) such that \(\cadpth\pr{t}=x\); since \(\cadpth\) is right-continuous in~\(\timeaxis\cap\rlims\pr{\timeaxis}=\reals_{\geq0}\), either this is the case for~\(t=T\) or for some rational time point~\(t\in\coi{0}{T}\cap\mathbb{Q}_{\geq 0}\).
	Hence,
	\begin{equation*}
		\mathring{A}
		= \bigcup_{t\in \countpoints} \ev{\cadX_t=x}
		\quad\text{with }
		\countpoints
		\coloneqq \st{T} \cup \pr*{\coi{0}{T}\cap\mathbb{Q}_{\geq 0}}.
	\end{equation*}
	Since \(\countpoints\) is clearly countable and the \(\sigma\)-algebra generated by~\(\cadcylevts\) is closed under countable unions, we infer from this equality that \(\mathring{A}\in\sigma\pr{\cadcylevts}\).
\end{example}
So it seems like $\sigma\pr{\cadcylevts}$ would be a suitable domain, but can we also obtain a probability measure on $\sigma\pr{\cadcylevts}$?
To that end, we need to show that the probability charge on $\cadcylevts$ induced by a consistent collection of finite-dimensional distributions is countably additive.
We will see that this is the case if and only if an extra regularity condition is satisfied.
Before we can introduce this condition in Section~\ref{sssec:regularity} further on, we need to introduce the notion of (the expected number of) jumps.

\subsection{The number of jumps}
\label{ssec:expected number of jumps}
A crucial property of càdlàg paths is that for all \(s,r\in\timeaxis\) such that \(s\leq r\), the number of `discontinuities' or `jumps' in~\(\cci{s}{r}\cap\timeaxis\) is finite.
For specific choices of $\timeaxis$ -- such as $\reals$ or $\reals_{\geq0}$ -- this property is well known -- see for example Lemma~5.20 in \cite{2021Erreygers-Phd} or Lemma~1 in \cite[Section~12]{1999Billingsley-Convergence} for related results -- and a rather straightforward consequence of the definition of càdlàg paths and the Heine-Borel Theorem. We will establish this result for general $\timeaxis$ though. But let us begin by explaining what we mean with `jumps'.

Fix \(s,r\in\reals\) such that \(s\leq r\) and a càdlàg path~\(\cadpth\in\cadpaths\).
First, we assume that \(\cci{s}{r}\subseteq\timeaxis\).
Then the number of jumps in~\(\cci{s}{r}\cap\timeaxis\) of~\(\cadpth\) is
\begin{equation*}
	\card[\Big]{\st[\Big]{t\in\oci{s}{r}\colon \lim_{\Delta\searrow 0} \cadpth\pr{t-\Delta}\neq \cadpth\pr{t}}},
\end{equation*}
and it follows fairly easily from the properties of càdlàg paths and the Heine-Borel Theorem -- or Billingsley's~\cite[Section~12, Lemma~1]{1999Billingsley-Convergence} argument -- that this number is finite.
Therefore, it is not difficult to verify that
\begin{equation}
\label{eqn:pthjumps for cadlag as sup}
	\card[\Big]{\st[\Big]{t\in\oci{s}{r}\colon \lim_{\Delta\searrow 0} \cadpth\pr{t-\Delta}\neq \cadpth\pr{t}}}
	= \sup\st[\big]{\tupjumps_u\pr{\cadpth\pr{u}}\colon u\in\setoftseq_{\cci{s}{r}\cap\timeaxis}},
\end{equation}
where for any tuple of time points~\(u=\pr{t_1, \dots, t_n}\in \setoftseq_{\timeaxis}\), we let
\begin{equation*}
	\tupjumps_u
	\colon \stsp_u\to\st{0, \dots, n-1}
	\colon x_u\mapsto
	\card{k\in\st{2, \dots, n}\colon x_{t_{k-1}}\neq x_{t_k}}.
\end{equation*}
If \(\cci{s}{r}\) is not included in~\(\timeaxis\), we cannot define the number of jumps in~\(\cci{s}{r}\) in the same way.
Instead, we generalise the right-hand side of Eqn.~\eqref{eqn:pthjumps for cadlag as sup}, but with the sufficient generality we need in the remainder: for all subsets~\(\altpoints\) of~\(\timeaxis\),
\begin{equation*}
	\pthjumps_{\altpoints}\pr{\allpth}
	\coloneqq \sup\st[\big]{\tupjumps_u\pr{\allpth\pr{u}}\colon u\in\setoftseq_{\altpoints}}\text{ for all }\allpth\in\allpaths.
\end{equation*}
In any case, the number of jumps of a \emph{càdlàg} path in~\(\cci{s}{r}\) is always finite; because the proof of this (essentially well-known) result is not very instructive, we have relegated it to Appendix~\ref{asec:Properties of cadlag paths}.
\begin{proposition}
\label{prop:number of jumps is finite}
	For all~\(\cadpth\in\cadpaths\) and~\(s,r\in\reals\) such that \(s<r\) and \(\cci{s}{r}\cap\timeaxis\neq\emptyset\), \(\pthjumps_{\cci{s}{r}\cap\timeaxis}\pr{\cadpth}<+\infty\).
\end{proposition}

Finally, we will also consider the expected number of jumps.
To this end, we observe that for all \(u=\pr{t_1, \dots, t_m}\in\setoftseq\), \(\tupjumps_u\) is a \(\wp\pr{\stsp_u}\)-simple variable\footnote{
	We adhere to the definition given by Troffaes \& De Cooman~\cite[Definition~1.16]{2014Troffaes-Lower}, but see also Definition~4.2.12 in \cite{1983RaoRao-Charges}.
} because
\begin{equation*}
	\tupjumps_u
	= \sum_{k=1}^{m-1} k \indica{\st{x_u\in\stsp_u\colon \tupjumps_u\pr{\stsp_u}=k}}
	= \sum_{k=2}^m \indica{\st{x_u\in\stsp_u\colon x_{t_k}\neq x_{t_{k-1}}}};
\end{equation*}
hence, for any probability charge~\(\dist\) on~\(\wp\pr{\stsp_u}\), its expectation with respect to~\(\dist\) is well-defined through the Dunford integral:\footnote{
	We adhere to the definition in \cite[Definition~8.13]{2014Troffaes-Lower}, see also Definition~4.4.1 in \cite{1983RaoRao-Charges}.
}
\begin{equation}
\label{eqn:prevdist jump}
	\prev_{\dist}\pr{\tupjumps_u}
	= \sum_{k=1}^{m-1} k\dist\pr{\st{x_u\in\stsp_u\colon \tupjumps_u\pr{x_u}=k}}
	= \sum_{k=2}^m \dist_u\pr[\big]{\st{x_u\in\stsp_u\colon x_{t_k}\neq x_{t_{k-1}}}}.
\end{equation}
If \(\dist_{\noarg}\) is a consistent collection of finite-dimensional charges, then for all \(u=\pr{t_1, \dots, t_m}\in\setoftseq\), the final expression in Eqn.~\eqref{eqn:prevdist jump} simplifies as follows:
\begin{equation}
\label{eqn:prevdist jump as sum over differences}
	\prev_{\dist_u}\pr{\tupjumps_u}
	= \sum_{k=2}^m \dist_u\pr[\big]{\st{x_u\in\stsp_u\colon x_{t_k}\neq x_{t_{k-1}}}}
	= \sum_{k=2}^m \dist_{\pr{t_{k-1}, t_k}}\pr[\big]{\stsp^2_{\neq}},
\end{equation}
where here and in the remainder, we let \(\stsp^2_{\neq}\coloneqq\st[\big]{\pr{x,y}\in\stsp^2\colon x\neq y}\) and \(\stsp^2_{=}\coloneqq\st[\big]{\pr{x,y}\in\stsp^2\colon x= y}\), and for all \(s,r\in\timepoints\) such that \(s<r\), with some slight abuse of notation, we also write \(\stsp^2_{\neq}\) and \(\stsp^2_{=}\) when we actually mean \(\st{\pr{x_s, x_r}\in\stsp_{\pr{x_s, x_r}}\colon x_s\neq x_r}\) and \(\st{\pr{x_s, x_r}\in\stsp_{\pr{x_s, x_r}}\colon x_s=x_r}\), respectively.





\subsection{Regularity}
\label{sssec:regularity}
Finally, we can get down to proving that the following two regularity conditions on a consistent collection of finite-dimensional distributions are necessary and sufficient for the corresponding probability charge on~\(\cadcylevts\) to be countably additive.
\begin{definition}
\label{def:regularity}
	A collection of finite-dimensional charges~\(\dist_{\noarg}\) is \emph{regular} if
	\begin{enumerate}[label=R\arabic*., ref=(R\arabic*)]
		\item\label{def:dist regular:cont} for all \(t\in\rlims\pr{\timepoints}\cap\timepoints\),
		\begin{equation*}
			\lim_{\timepoints\ni r\searrow t} \dist_{\pr{t,r}}\pr[\big]{\stsp^2_=}
			= 1;
		\end{equation*}
		\item\label{def:dist regular:bound on exp of jumps} for all \(n\in\nats\) such that \(\cci{-n}{n}\cap\timepoints\neq\emptyset\),
		\begin{equation*}
			\lim_{k\to+\infty} \sup\st[\big]{\dist_u\pr[\big]{\st{x_u\in\stsp_u\colon \tupjumps_u\pr{x_u}\geq k}}\colon u\in\setoftseq_{\cci{-n}{n}\cap\timepoints}}
			= 0.
		\end{equation*}
	\end{enumerate}
\end{definition}
The two conditions for regularity basically ensure that the collection of finite-dimensional charges is compatible with the properties of càdlàg paths.
The first condition \ref{def:dist regular:cont} mirrors the continuity from the right: it demands that the finite-dimensional charges are (stochastically) continuous from the right -- similar to the notions in, for example, \cite[Definition~12.15]{1992Breiman-Probability} or \cite[Definition~18.2.2]{2013Borovkov-Probability}.
It may be less obvious on first inspection, but the second condition \ref{def:dist regular:bound on exp of jumps} mirrors the fact that càdlàg paths have a finite number of jumps in any closed and bounded interval.
While Proposition~\ref{prop:number of jumps is finite} says that this should be the case for all intervals~\(\cci{s}{r}\), it clearly suffices to limit ourselves to intervals of the form~\(\cci{s}{r}=\cci{-n}{n}\) with \(n\in\nats\) because the number of jumps is monotone.

That regularity is indeed necessary and sufficient for a consistent collection of finite-dimensional distributions to have a corresponding probability charge~\(\prob\) on~\(\cadcylevts\) that is countably additive, is established by our next result.
Our proof is a bit long and needs quite a bit more technical machinery; for this reason, we have relegated it to Appendix~\ref{asec:proof for the:precursor for cadlag}.
\begin{theorem}
\label{the:precursor for cadlag paths}
	Consider a consistent collection~\(\dist_{\noarg}\) of finite-dimensional distributions, and let \(\prob\) be the corresponding probability charge on~\(\cadcylevts\) of Proposition~\ref{prop:Precursor of Kolmogorov}.
	Then \(\prob\) is countably additive if an only if \(\dist_{\noarg}\) is regular.
\end{theorem}

We now use Theorem~\ref{the:precursor for cadlag paths} to prove our main result: a version of the Daniell--Kolmogorov Extension Theorem where the set of all paths is replaced by the set of all càdlàg paths.
\begin{theorem}
\label{the:regular then daniell-kolmogorov for cadlag}
	Consider a consistent collection of finitary distributions~\(\dist_{\noarg}\).
	Then there is a probability measure~\(\prob\) on~\(\sigma\pr{\cadcylevts}\) such that
	\begin{equation*}
		\prob\pr{\ev{\cadX_u\in A}}
		= \dist_u\pr{A}
		\quad\text{for all } u\in\setoftseq, A\in\wp\pr{\stsp_u}
	\end{equation*}
	if and only if \(\dist_\noarg\) is regular; if this is the case, then this probability measure is unique.
\end{theorem}
\begin{proof}
	Follows immediately from Proposition~\ref{prop:Precursor of Kolmogorov}, Theorem~\ref{the:precursor for cadlag paths} and Caratheodory's Extension Theorem.
\end{proof}

That we need the two conditions \ref{def:dist regular:cont} and \ref{def:dist regular:bound on exp of jumps} for Theorem~\ref{the:regular then daniell-kolmogorov for cadlag} is in line with Borovkov's~\cite{2013Borovkov-Probability} findings for \(\reals\)-valued processes: right after his Definition~18.2.3, he explains that stochastic continuity (so his version of~\ref{def:dist regular:cont}) alone does not suffice, and that `in order to characterise the properties of trajectories, one needs quantitative bounds for [the magnitude of the jumps]'.
Note that in our setting of countable state spaces, we do \emph{not} need a quantitative bound on the `magnitude' of the jumps but on the \emph{number} of jumps.

The two conditions \ref{def:dist regular:cont} and \ref{def:dist regular:bound on exp of jumps} for regularity are perhaps not the most easy ones to check.
Fortunately, there are plenty of sufficient conditions for regularity that can be more easily verified.
For example, the one in the following result comes in handy in the setting of Poisson processes and Markovian (imprecise) jump processes.
\begin{proposition}
\label{prop:only expected number of jumps for regularity}
	Consider a consistent collection~\(\dist_{\noarg}\) of finite-dimensional distributions.
	If for all \(n\in\nats\) with \(\cci{-n}{n}\cap\timepoints\neq\emptyset\) there is some \(\lambda_n\in\nats\) such that
	\begin{equation*}
		\prev_{\dist_u}\pr{\tupjumps_u}
		\leq \lambda_n\pr{t_m-t_1}
		\quad\text{for all } u=\pr{t_1, \dots, t_m}\in\setoftseq_{\cci{-n}{n}\cap\timepoints},
	\end{equation*}
	then \(\dist_{\noarg}\) is regular.
\end{proposition}
\begin{proof}
	To verify \ref{def:dist regular:cont}, we fix any \(t\in\timepoints\cap\rlims\pr{\timepoints}\).
	Let \(n\) be any natural number such that \(\abs{t}<n\).
	Then \(\cci{-n}{n}\cap\timepoints\neq\emptyset\) by construction and \(\oci{t}{n}\cap\timepoints\neq\emptyset\) because \(t\) is a right-sided limit point of~\(\timepoints\).
	For all \(r\in\oci{t}{n}\cap\timepoints\) and with \(u\coloneqq\pr{t,r}\), \(\st{x_u\in\stsp_u\colon\tupjumps_u\pr{x_u}=1}=\stsp^2_{\neq}\), and it follows from this and the condition in the statement that
	\begin{equation*}
		\dist_{u}\pr[\big]{\stsp^2_{\neq}}
		= \dist_u\pr{\st{x_u\in\stsp_u\colon\tupjumps_u\pr{x_u}=1}}
		= \prev_{\dist_u}\pr{\tupjumps_u}
		\leq \lambda_n\pr{r-t}.
	\end{equation*}
	From this inequality, we infer that
	\begin{equation*}
		\lim_{\timepoints\ni r\searrow t} \dist_{\pr{t,r}}\pr[\big]{\stsp^2_{=}}
		= \lim_{\timepoints\ni r\searrow t} 1-\dist_{\pr{t,r}}\pr[\big]{\stsp^2_{\neq}}
		= 1,
	\end{equation*}
	as required for \ref{def:dist regular:cont}.

	Checking \ref{def:dist regular:bound on exp of jumps} is straightforward.
	Fix some \(n\in\nats\) with \(\cci{-n}{n}\cap\timepoints\neq\emptyset\).
	For all \(k\in\nats\) and \(u=\pr{t_1, \dots, t_m}\in\setoftseq_{\cci{-n}{n}\cap\timepoints}\), it follows from Markov's inequality and the condition in the statement that
	\begin{equation*}
		\dist_u\pr{\st{x_u\in\stsp_u\colon \tupjumps_u\pr{x_u}\geq k}}
		\leq \frac{\prev_{\dist_u}\pr{\tupjumps_u}}{k}
		\leq \frac{\lambda_n\pr{t_m-t_1}}{k}
		\leq \frac{2n\lambda_n}{k}.
	\end{equation*}
	From this, we infer that
	\begin{equation*}
		\lim_{k\to+\infty} \sup\st[\big]{\dist_u\pr{\st{x_u\in\stsp_u\colon \tupjumps_u\pr{x_u}\geq k}}\colon u\in\setoftseq_{\cci{-n}{n}\cap\timepoints}}
		= 0,
	\end{equation*}
	as required.
\end{proof}

The consistent collection of finite-dimensional distributions~\(\dist_{\noarg}\) for the Poisson process with rate~\(\lambda\in\nnegreals\), as introduced at the end of Section~\ref{ssec:finite charges and distributions}, satisfies the sufficient condition in Proposition~\ref{prop:only expected number of jumps for regularity}: it is not difficult to verify that for all \(u=\pr{t_1, \dots, t_m}\in\setoftseq\),
\begin{equation*}
	\prev_{\dist_u}\pr{\tupjumps_u}
	= \sum_{k=2}^m \lambda\pr{t_k-t_{k-1}}
	= \lambda\pr{t_m-t_1}.
\end{equation*}
The same is true for the collections of finite-dimensional distributions for finite-state Markov processes whose semi-group of transition matrices is generated by a (transition) rate matrix.
Even more, this sufficient condition also suffices in the more general setting of Markovian imprecise jump processes with a finite state space \cite{2017KrakDeBock,2022Erreygers}.
In that setting, one considers sets of consistent collections of finite-dimensional distributions that are `consistent' -- in the sense of \cite[Definition~6.1]{2017KrakDeBock} or \cite[Definition~3.50]{2021Erreygers-Phd} -- with a non-empty and bounded set~\(\mathcal{Q}\) of (transition) rate matrices, and Erreygers~\cite[Corollary~5.18 and Theorem~5.27]{2021Erreygers-Phd} proves that for each of these collections of finite-dimensional distributions,
\begin{equation*}
	\prev_{\dist_u}\pr{\tupjumps_u}
	\leq \frac12\Vert \mathcal{Q}\Vert\pr{t_m-t_1}
	< +\infty
	\quad\text{for all } u=\pr{t_1, \dots, t_m}\in\setoftseq,
\end{equation*}
where \(\Vert \mathcal{Q}\Vert\) is the supremum over the norms of the rate matrices in the bounded set~\(\mathcal{Q}\).

Our work on Markovian imprecise jump processes has also inspired us to establish two additional sufficient conditions for regularity.
The main difference with the condition in Proposition~\ref{prop:only expected number of jumps for regularity} is that these conditions only concern the finite-dimensional distributions for two time points; they can be thought of as placing bounds on the dynamics of the process.
Because these results are of a more technical nature, we have relegated them to Appendix~\ref{asec:two additional sufficient conditions for regularity}.
Here, we only give a simplified version of Proposition~\ref{prop:stronger II}.
\begin{corollary}
\label{cor:stronger II}
	Suppose that \(\timepoints=\nnegreals\), and consider a consistent collection~\(\dist_{\noarg}\) of finite-dimensional distributions.
	If for all \(n\in\nats\) there is some~\(\lambda_n\in\nnegreals\) such that
	\begin{align*}
		\limsup_{r\searrow t}\frac{\dist_{\pr{t,r}}\pr[\big]{\stsp^2_{\neq}}}{r-t}
		&\leq \lambda_n
		\quad\text{for all } t\in\coi{0}{n}
	\shortintertext{and}
		\limsup_{s\nearrow t}\frac{\dist_{\pr{s,t}}\pr[\big]{\stsp^2_{\neq}}}{t-s}
		&\leq \lambda_n
		\quad\text{for all } t\in\oci{0}{n},
	\end{align*}
	then \(\dist_{\noarg}\) is regular.
\end{corollary}
Note that Corollary~\ref{cor:stronger II} is also relevant when dealing with the Poisson process: with \(\dist_{\noarg}\) the consistent collection of finite-dimensional distributions for the Poisson process with rate~\(\lambda\),
\begin{equation*}
	\lim_{r\searrow t}\frac{\dist_{\pr{t,r}}\pr{\stsp^2_{\neq}}}{r-t}
	= \lim_{r\searrow t} \frac{1-e^{-\lambda\pr{r-t}}}{r-t}
	= \lambda
	\quad\text{for all } t\in\nnegreals
\end{equation*}
and
\begin{equation*}
	\lim_{s\nearrow t}\frac{\dist_{\pr{s,t}}\pr{\stsp^2_{\neq}}}{t-s}
	= \lim_{s\nearrow t} \frac{1-e^{-\lambda\pr{t-s}}}{t-s}
	= \lambda
	\quad\text{for all } t\in\posreals
\end{equation*}

\section{Conclusion}
Our version of the Daniell--Kolmogorov Extension Theorem (Theorem~\ref{the:regular then daniell-kolmogorov for cadlag}) does not suffer from the issue with the original result (Theorem~\ref{the:Kolomogorov extension theorem}): by using the càdlàg paths as sample space, we obtain a probability measure whose domain is indeed sufficiently rich, in that it contains most -- if not all -- of the events that one can be interested in.
Moreover, our results considers stochastic processes indexed by any subset~\(\timeaxis\) of the real numbers and looked at events that depend only on the time points in some subset~\(\timepoints\) of~\(\timeaxis\), while it is standard to assume~\(\timepoints=\timeaxis\); this distinction between \(\timeaxis\) and \(\timepoints\) may seem somewhat strange at first, but it was motivated by our work on Markovian imprecise jump processes, where this proved to be useful.

In contrast to the original version of the theorem, our new version is limited to the specific setting of countable-state stochastic processes due to the form of our regularity condition.
For this reason, a first line of follow-up research we envision is to generalise our work to uncountable state spaces.
Crucial is that one should find a meaningful alternative to the number of jumps; for example, with the set of real numbers, it is customary -- in the setting of modifications, at least -- to look at upcrossing numbers instead of the number of jumps.
A second line of follow-up research could be to apply our results.
We have already briefly touched on the relevance to the Poisson process, and we are convinced that our theorem can be useful for more general counting processes as well.
Furthermore, we know from our work in the setting of Markovian imprecise jump processes that our version of the Daniell--Kolmogorov Extension Theorem can be of use in the setting of jump processes that need not be time-homogeneous nor even Markovian.

\bibliographystyle{amsplain}
\bibliography{bibliography}

\appendix
\section{Proof for Theorem~\ref{the:precursor for all paths}}
\label{asec:a proof for Kolmogorov extension theorem}
In case \(\timepoints=\timeaxis\), Theorem~\ref{the:precursor for all paths} is implied by Theorem~31.1 in \cite[Chapter~II]{1994Rogers-Diffusions} -- we leave it to the reader to check that \(\pr{\stsp, \wp\pr{\stsp}}\) is a Lusin space.
We could also show that this in turn implies Theorem~\ref{the:precursor for all paths} in case \(\timepoints\subset\timeaxis\), but we choose to take a more direct route towards proving this result; our reasons for this are twofold: (i) it does not take a lot more work, and (ii) it would be a shame not to elucidate how the argument simplifies in our particular case where $\stsp$ is countable.

Our proof for Theorem~\ref{the:precursor for all paths} is essentially the argument used by Rogers \& Williams~\cite[Chapter~II]{1994Rogers-Diffusions} (who assume that \(\stsp\) is a complete metric space) to prove their Lemma 30.7, but translated to the setting of a countable state space~\(\stsp\) -- note that if \(\stsp\) is countably infinite, then \(\pr{\stsp, \wp\pr{\stsp}}\) is not a complete metric space. The crucial result that allows us to avoid the assumption that $\stsp$ is a commplete metric space, is the following intermediate lemma.
\begin{lemma}
\label{lem:finite intersection paths}
For all $n\in\nats$, consider a finite subset $B_n$ of $\stsp_{u_n}$, with $u_n\in\setoftseq$, and let $\dot{B}_n\coloneqq\ev{\allX_{u_n}\in B_n}$.
If $\bigcap_{n=1}^m\dot{B}_n\neq\emptyset$ for all $m\in\nats$, then also $\bigcap_{n\in\nats}\dot{B}_n\neq\emptyset$.
\end{lemma}
\begin{proof}
We endow $\stsp$ with the cofinite topology \cite[Section 5.15.c]{1997Schechter-Handbook}, the open sets of which are $\{\emptyset\}\cup\{A^{\compl}\colon A\subseteq\stsp,\abs{A}<+\infty\}$, and we endow $\allpaths=\stsp^{\timeaxis}$ with the corresponding product topology. Since any product of topological spaces, each equiped with the cofinite topology, is compact \cite[Section~17.20 (AC25)]{1997Schechter-Handbook}, this turns $\allpaths$ into a compact space.

For all $x\in\stsp$ and $t\in\timeaxis$, since $\{x\}^{\compl}$ belongs to the cofinite topology on $\stsp$, we know that $\ev{\allX_t=x}^{\compl}=\ev{\allX_t\in\st{x}^{\compl}}=\st{\allpth\in\allpaths\colon\allpth(t)\in\st{x}^{\compl}}$ belongs to the product topology on $\allpaths$.
For all $u=\pr{t_1, \dots, t_k}\in\setoftseq$ and finite $B\subseteq\stsp_u$, since
\begin{equation*}
\ev{\allX_u\in B}
=\bigcup_{x_u\in B}\ev{\allX_u=x_u}
=\bigcup_{x_u\in B}\bigcap_{\ell=1}^k\ev{\allX_{t_\ell}=x_{t_\ell}},
\end{equation*}
this implies that
\begin{equation*}
\{\allX_u\in B\}^{\compl}
=\bigcap_{x_u\in B}\bigcup_{\ell=1}^k\{\allX_{t_\ell}=x_{t_\ell}\}^{\compl}
\end{equation*}
belongs to the product topology as well.
Hence, for all $n\in\nats$, $\dot{B}_n=\ev{\allX_{u_n}\in B_n}$ is a closed subset of $\allpaths$.

Let us now assume that $\bigcap_{n=1}^m\dot{B}_n\neq\emptyset$ for all \(m\in\nats\). Then clearly, for any finite $N\subseteq\nats$, we have that $\bigcap_{n\in N}\dot{B}_n\neq\emptyset$. The compactness of~\(\allpaths\) therefore implies that $\bigcap_{n\in\nats}\dot{B}_n\neq\emptyset$ -- see \cite[Definition~17.2.(B)]{1997Schechter-Handbook}.
\end{proof}

Our proof for Theorem~\ref{the:precursor for all paths} now follows more or less the same argument as the one for Lemma~30.7 in \cite[Chapter~II]{1994Rogers-Diffusions}, but without the use of their Lemma 29.7.
\begin{proof}[Proof of Theorem~\ref{the:precursor for all paths}]
According to a classic result from measure theory -- see for example \cite[Chapter~7, Proposition~9]{1997Fristedt-Modern} or \cite[Lemma~4.3]{1994Rogers-Diffusions} -- $P$ is countably additive if and only if, for any decreasing sequence $\pr{\dot{A}_n}_{n\in\nats}$ in $\allcylevts$ such that $\bigcap_{n\in\nats}\dot{A}_n=\emptyset$, also $\lim_{n\to+\infty}\prob\pr{\dot{A}_n}=0$.
We will prove this alternative condition instead.

By definition of $\allcylevts$, there are, for all $n\in\nats$, some $u_n\in\setoftseq$ and $A_n\in\stsp_{u_n}$ such that $\dot{A}_n=\ev{\allX_{u_n}\in A_n}$.
Fix any $\epsilon\in\posreals$.
For all $n\in\nats$, since $\dist_{u_n}$ is countably additive and $\stsp_{u_n}$ -- and hence also $A_n$ -- is a countable set, we have that $\dist_{u_n}\pr{A_n}=\sum_{x_{u_n}\in A_n}\dist_{u_n}\pr{\ev{x_{u_n}}}$.
Since every $\dist_{u_n}\pr{\ev{x_{u_n}}}$ is positive, the convergence of this countable sum to $\dist_{u_n}\pr{A_n}$ implies that there is a finite subset $B_n$ of $A_n$ such that $\dist_{u_n}\pr{B_n}=\sum_{x_{u_n}\in B_n}\dist_{u_n}\pr{\ev{x_{u_n}}}>\dist_{u_n}(A_n)-2^{-n}\epsilon$ and hence also $\dist_{u_n}\pr{A_n\setminus B_n}<2^{-n}\epsilon$.
Consider any such $B_n$ and let $\dot{B}_n\coloneqq\ev{\allX_{u_n}\in B_n}\subseteq\dot{A}_n$.

Since $\bigcap_{n\in\nats}\dot{A}_n=\emptyset$, we clearly also have that $\bigcap_{n\in\nats}\dot{B}_n=\emptyset$.
It therefore follows from Lemma~\ref{lem:finite intersection paths} that there is some $m\in\nats$ such that $\bigcap_{n=1}^m\dot{B}_n=\emptyset$.
Then for any $\allpth\in\dot{A}_m$, since $\allpth\notin\bigcap_{n=1}^m\dot{B}_n$, there is some $1\leq n_{\allpth}\leq m$ such that $\allpth\notin\dot{B}_{n_{\allpth}}$ and therefore $\allpth\in\dot{A}_m\setminus\dot{B}_{n_{\allpth}}\subseteq\dot{A}_{n_{\allpth}}\setminus\dot{B}_{n_{\allpth}}\subseteq\bigcup_{n=1}^m\dot{A}_n\setminus\dot{B}_n$.
This implies that $\dot{A}_m\subseteq\bigcup_{n=1}^m\dot{A}_n\setminus\dot{B}_n$ and therefore, that
\begin{align*}
\prob(\dot{A}_m)
\leq\prob\pr*{\bigcup_{n=1}^m\dot{A}_n\setminus\dot{B}_n}
\leq\sum_{n=1}^m\prob(\dot{A}_n\setminus\dot{B}_n)
&=\sum_{n=1}^m\prob\pr{\ev{\allX_{u_n}\in A_n\setminus B_n}}\\
&=\sum_{n=1}^m\dist_{u_n}\pr{A_n\setminus B_n}
<\sum_{n=1}^m2^{-n}\epsilon<\epsilon.
\end{align*}
Since $\pr{\dot{A}_n}_{n\in\nats}$ is a decreasing sequence, it follows that $\lim_{n\to+\infty}\prob\pr{\dot{A}_n}<\epsilon$.
As the choice of $\epsilon>0$ was arbitrary, we conclude that $\lim_{n\to+\infty}\prob\pr{\dot{A}_n}=0$.
\end{proof}

\section{Properties of càdlàg pahts}
\label{asec:Properties of cadlag paths}
In this appendix, we give the proof for the three important properties of càdlàg paths.
The first property is that a càdlàg path always has a càdlàg extension to~\(\reals\).
\begin{proof}[Proof of Lemma~\ref{lem:extension of cadlag path}]
	Fix some state~\(x\in\stsp\), and for all \(t\in\reals\), let \(\mathcal{S}_{\leq t}\coloneqq\oci{-\infty}{t}\cap\mathcal{S}\) and \(s^{\star}_t\coloneqq\sup\mathcal{S}_{\leq t}\).
	Then
	\begin{equation*}
		\rpth
		\colon\reals\to\stsp
		\colon t\mapsto \begin{cases}
			\vpth\pr{t}	&\text{if } t\in\mathcal{S} \\
			\lim_{\mathcal{S}\ni r\searrow t}\vpth\pr{r} &\text{if } t\notin\mathcal{S}, t\in\rlims\pr{\mathcal{S}} \\
			\vpth\pr{s^{\star}_t} &\text{if } t\notin\mathcal{S}, t\notin\rlims\pr{\mathcal{S}}, \mathcal{S}_{\leq t}\neq\emptyset, s^\star_t\in\mathcal{S} \\
			\lim_{\mathcal{S}\ni s\nearrow s^{\star}_t} \vpth\pr{s} &\text{if } t\notin\mathcal{S}, t\notin\rlims\pr{\mathcal{S}}, \mathcal{S}_{\leq t}\neq\emptyset, s^\star_t\notin\mathcal{S} \\
			x &\text{if }t\notin\mathcal{S}, t\notin\rlims\pr{\mathcal{S}}, \mathcal{S}_{\leq t}=\emptyset
		\end{cases}
	\end{equation*}
	clearly coincides with~\(\vpth\) on~\(\mathcal{S}\).
	Less obvious is that \(\rpth\) (i) has a left-sided limit everywhere and (ii) is continuous from the right.

	For the first condition, we fix any \(t\in\reals\), and distinguish two cases.
	The first case is that \(t\) is not a left-sided limit point of~\(\mathcal{S}\).
	Then there is some \(\delta\in\posreals\) such that \(\ooi{t-\delta}{t}\cap\mathcal{S}=\emptyset\).
	Then for all \(s\in\ooi{t-\delta}{t}\), \(s\notin\mathcal{S}\cup\rlims\pr{\mathcal{S}}\) -- because with \(\delta'\coloneqq t-s>0\), \(\coi{s}{s+\delta'}\cap\mathcal{S}\subseteq\ooi{t-\delta}{t}\cap\mathcal{S}=\emptyset\) -- and \(\mathcal{S}_{\leq s}=\mathcal{S}_{\leq t-\delta}\) -- because \(\oci{t-\delta}{s}\cap\mathcal{S}\subseteq\ooi{t-\delta}{t}\cap\mathcal{S}=\emptyset\).
	It follows from this and the definition of~\(\rpth\) that for all \(s\in\ooi{t-\delta}{t}\),
	\begin{equation*}
		\rpth\pr{s}
		= \begin{cases}
			\vpth\pr{s^{\star}_{t-\delta}} &\text{if } \mathcal{S}_{\leq t-\delta}\neq\emptyset, s^{\star}_{t-\delta}\in\mathcal{S} \\
			\lim_{\mathcal{S}\ni r\nearrow s^\star_{t-\delta}}\vpth\pr{r} &\text{if } \mathcal{S}_{\leq t-\delta}\neq\emptyset, s^{\star}_{t-\delta}\notin\mathcal{S}, \\
			x &\text{otherwise.}
		\end{cases}
	\end{equation*}
	From this, we conclude that \(\rpth\) is constant on~\(\ooi{t-\delta}{t}\).

	The second case is that \(t\) is a left-sided limit point of~\(\mathcal{S}\).
	Because \(\vpth\) is càdlàg, there are some \(\delta'\in\posreals\) and \(x_t\in\stsp\) such that \(\vpth\pr{s}=x_t\) for all \(s\in\ooi{t-\delta'}{t}\cap\mathcal{S}\neq\emptyset\).
	Take any \(s'\in\ooi{t-\delta'}{t}\cap\mathcal{S}\).
	Then for all \(s\in\ooi{s'}{t}\cap\mathcal{S}\), \(\vpth\pr{s}=\vpth\pr{s'}=x_t\) and \(s'\leq s^\star_s<t\).
	From this, it follows that for all \(s\in\ooi{s'}{t}\), (i) \(\vpth\pr{s}=x_t\) if \(s\in\mathcal{S}\); (ii) \(\lim_{\mathcal{S}\ni r\searrow s}\vpth\pr{r}=x_t\) if \(s\in\rlims\pr{\mathcal{S}}\setminus\mathcal{S}\); (iii) \(\vpth\pr{s^\star_s}=x_t\) if \(s\notin\mathcal{S}\cup\rlims\pr{\mathcal{S}}\) and \(s^\star_s\in\mathcal{S}\); and (iv) \(\lim_{\mathcal{S}\ni r\nearrow s^\star_s}\vpth\pr{r}=x_t\) if $s\notin\mathcal{S}\cup\rlims\pr{\mathcal{S}}$ and $s^\star_s\notin\mathcal{S}$ -- the case \(\mathcal{S}_{\leq s}=\emptyset\) is clearly impossible because \(s'\in\mathcal{S}_{\leq s}\).
	We conclude from this and the definition of~\(\rpth\) that \(\rpth\pr{s}=x_t\) for all \(s\in\ooi{s'}{t}\).
	In both cases, this shows that \(\rpth\) has a left-sided limit in~\(t\).

	For the second condition, we fix any \(t\in\reals\).
	Here too, we distinguish two cases.
	The first case is that \(t\) is not a right-sided limit point of~\(\mathcal{S}\), meaning that there is some \(\delta\in\posreals\) such that \(\ooi{t}{t+\delta}\cap\mathcal{S}=\emptyset\).
	Then for all \(r\in\ooi{t}{t+\delta}\), \(r\notin\mathcal{S}\cup\rlims\pr{\mathcal{S}}\) -- because with \(\delta'\coloneqq t+\delta-r>0\), \(\coi{r}{r+\delta'}\cap\mathcal{S}\subseteq\ooi{t}{t+\delta}\cap\mathcal{S}=\emptyset\) -- and \(\mathcal{S}_{\leq r}=\mathcal{S}_{\leq t}\) -- because \(\oci{t}{r}\cap\mathcal{S}\subseteq\ooi{t}{t+\delta}\cap\mathcal{S}=\emptyset\).
	We distinguish two subcases.
	On the one hand, if \(t\in\mathcal{S}\) then it follows from this and our definition of~\(\rpth\) that (i) \(\rpth\pr{t}=\vpth\pr{t}\) and (ii) for all \(r\in\ooi{t}{t+\delta}\), \(s^\star_r=\sup\mathcal{S}_{\leq r}=\sup\mathcal{S}_{\leq t}=t\in\mathcal{S}\) and therefore \(\rpth\pr{r}=\vpth(s^\star_r)=\vpth\pr{t}\).
	If on the other hand \(t\notin\mathcal{S}\) and therefore $t\notin\mathcal{S}\cup\rlims(\mathcal{S})$, then it follows from this and the definition of~\(\rpth\) that for all \(r\in\ooi{t}{t+\delta}\), \(\rpth\pr{r}=\rpth\pr{t}\).

	This leaves the case that \(t\) is a right-sided limit point of~\(\mathcal{S}\).
	Because \(\vpth\) is càdlàg, there are some \(\delta\in\posreals\) and \(x_t\) such that \(\vpth\pr{r}=x_t\) for all \(r\in\coi{t}{t+\delta}\cap\mathcal{S}\neq\emptyset\).
	Here too, we distinguish two subcases.
	If \(t\in\mathcal{S}\), then we infer from this that (i) \(\rpth\pr{t}=\vpth\pr{t}=x_t\); (ii) for all \(r\in\ooi{t}{t+\delta}\cap\pr{\mathcal{S}\cup\rlims\pr{\mathcal{S}}}\), \(\rpth\pr{r}=x_t\); and (iii) for all \(r\in\ooi{t}{t+\delta}\setminus\pr{\mathcal{S}\cup\rlims\pr{\mathcal{S}}}\), \(\mathcal{S}_{\leq r}\neq\emptyset\) and \(s^\star_r\in\ooi{t}{t+\delta}\), and therefore \(\rpth\pr{r}=x_t\).
	If on the other hand \(t\notin\mathcal{S}\), then we infer from this that (i) \(\rpth\pr{t}=\lim_{\mathcal{S}\ni s\searrow t}\vpth\pr{s}=x_t\); (ii) for all \(r\in\ooi{t}{t+\delta}\cap\pr{\mathcal{S}\cup\rlims\pr{\mathcal{S}}}\), \(\rpth\pr{r}=x_t\); and (iii) for all \(r\in\ooi{t}{t+\delta}\setminus\pr{\mathcal{S}\cup\rlims\pr{\mathcal{S}}}\), \(\ooi{t}{r}\cap\mathcal{S}\neq\emptyset\) -- as \(t\) is a right-sided limit point of~\(\mathcal{S}\) -- and therefore $\mathcal{S}_{\leq r}\neq\emptyset$ and \(s^\star_r\in\ooi{t}{t+\delta}\), whence
	\begin{equation*}
		\rpth\pr{r}
		= \begin{cases}
			\vpth\pr{s^\star_r}=x_t &\text{if } s^\star_r\in\mathcal{S}, \\
			\lim_{\mathcal{S}\ni s\nearrow s^\star_r} \vpth\pr{s}=x_t &\text{if } s^\star_r\notin\mathcal{S}.
		\end{cases}
	\end{equation*}
	It is clear in both cases that \(\rpth\) is continuous from the right in~\(t\), as required.
\end{proof}

The second property we need to prove is Lemma~\ref{lem:cadlag defined on subset}.
\begin{proof}[Proof of Lemma~\ref{lem:cadlag defined on subset}]
	The direct implication is trivial because \(\densepoints\subseteq\timeaxis\), so it remains for us to show the converse implication.
	To this end, we assume that \(\cadpth_1\pr{d}=\cadpth_2\pr{d}\) for all \(d\in\densepoints\), and set out to show that then \(\cadpth_1\pr{t}=\cadpth_2\pr{t}\) for all \(t\in\timeaxis\).
	Fix any time point~\(t\) in~\(\timeaxis\).
	If \(t\in\densepoints\), then \(\cadpth_1\pr{t}=\cadpth_2\pr{t}\) by assumption, as required.
	If \(t\notin\densepoints\), then \(t\in\rlims\pr{\densepoints}\) because \(\timeaxis\subseteq\densepoints\cup\rlims\pr{\densepoints}\) by assumption.
	Because \(\cadpth_1\) and \(\cadpth_2\) are both càdlàg, there are some \(\delta_1, \delta_2\in\posreals\) such that
	\begin{equation*}
		\pr[\big]{\forall r\in\timeaxis\cap\ooi{t}{t+\delta_1}}~
		\cadpth_1\pr{r} = \cadpth_1\pr{t}
		\quad\text{and}\quad
		\pr[\big]{\forall r\in\timeaxis\cap\ooi{t}{t+\delta_2}}~
		\cadpth_2\pr{r} = \cadpth_2\pr{t}.
	\end{equation*}
	Let \(\delta\coloneqq\min\st{\delta_1, \delta_2}\).
	Fix any \(d\in\ooi{t}{t+\delta}\cap\densepoints\) -- this is always possible because \(t\) is a right-sided limit point of~\(\densepoints\), so \(\ooi{t}{t+\delta}\cap\densepoints\neq\emptyset\).
	On the one hand, it follows from the preceding that \(\cadpth_1\pr{t}=\cadpth_1\pr{d}\) and \(\cadpth_2\pr{t}=\cadpth_2\pr{d}\).
	On the other hand, we know that \(\cadpth_1\pr{d}=\cadpth_2\pr{d}\) because \(d\in\densepoints\).
	From this, we infer that \(\cadpth_1\pr{t}=\cadpth_2\pr{t}\).
	Since this is the case for all \(t\in\timeaxis\), we have shown that \(\cadpth_1=\cadpth_2\), as required.
\end{proof}

Finally, we prove that a càdlàg path can only jump a finite number of times in any closed and bounded interval.
\begin{proof}[Proof of Proposition~\ref{prop:number of jumps is finite}]
	By Lemma~\ref{lem:extension of cadlag path}, there is some càdlàg map~\(\rpth\colon\reals\to\stsp\) that coincides with~\(\cadpth\) on~\(\timeaxis\).
	As we have explained in Section~\ref{ssec:expected number of jumps}, it is not difficult to show that
	\begin{equation*}
		\sup\st[\big]{\tupjumps_u\pr{\rpth\pr{u}}\colon u\in\setoftseq_{\cci{s}{r}}}
		= \card[\Big]{\st[\Big]{t\in\oci{s}{r}\colon \lim_{\Delta\searrow 0}\rpth\pr{t-\Delta}\neq\rpth\pr{t}}}
		< +\infty,
	\end{equation*}
	where for all \(u=\pr{t_1, \dots, t_n}\in\setoftseq_{\cci{s}{r}}\), we let \(\rpth\pr{u}\coloneqq\pr{\rpth\pr{t_1}, \dots, \rpth\pr{t_n}}\).
	Since clearly
	\begin{equation*}
		\pthjumps_{\cci{s}{r}\cap\timeaxis}\pr{\cadpth}
		= \sup\st[\big]{\tupjumps_u\pr{\rpth\pr{u}}\colon u\in\setoftseq_{\cci{s}{r}\cap\timeaxis}}
		\leq \sup\st[\big]{\tupjumps_u\pr{\rpth\pr{u}}\colon u\in\setoftseq_{\cci{s}{r}}},
	\end{equation*}
	this proves the statement.
\end{proof}

\section{Proof for Theorem~\ref{the:precursor for cadlag paths}}
\label{asec:proof for the:precursor for cadlag}
In this appendix, we prove Theorem~\ref{the:precursor for cadlag paths}.
We will get to this in Appendix~\ref{asec:actual proof}, but first we need to introduce some additional technical machinery.
In Appendix~\ref{assec:expected number of jumps again}, we take a second look at the expected number of jumps, but this time in the setting of probability measures and stochastic processes.
This theme continues in Appendix~\ref{assec:modification}, where we prove that whenever a (not necessarily càdlàg) path~\(\pth\) has a finite number of jumps (for some countable subset of~\(\timepoints\)), we can always `modify' this path in such a way that it becomes càdlàg.
Appendix~\ref{assec:two additional intermediary results} introduces the two remaining intermediary results that we will need to prove Theorem~\ref{the:precursor for cadlag paths}.

\subsection{The expected number of jumps, again}
\label{assec:expected number of jumps again}
Consider again a subset~\(\genpaths\) of~\(\allpaths\) that satisfies Eqn.~\eqref{eqn:condition on set of paths}.
Then for any tuple of time points~\(u=\pr{t_1, \dots, t_n}\in\setoftseq\), the `number of jumps of~\(\genX_{\noarg}\) along~\(u\)' is the functional composition of~\(\tupjumps_u\) after the projection~\(\genX_u\):
\begin{equation*}
	\tupjumps_u\pr{\genX_u}
	\colon \genpaths\to\st{0, \dots, n-1}
	\colon \pth\mapsto \tupjumps_u\pr{\pth\pr{u}}.
\end{equation*}
It is easy to see that
\begin{equation}
\label{eqn:jumps as sum alt}
	\tupjumps_u\pr{\genX_u}
	= \sum_{k=2}^n \indica{\ev{\genX_{t_{k-1}}\neq \genX_{t_k}}},
\end{equation}
where for all \(s,r\in\timepoints\) such that \(s<r\), we let
\begin{equation*}
	\ev{\genX_s\neq \genX_r}
	\coloneqq \ev{\genX_{\pr{s,r}}\in\stsp^2_{\neq}}
	\in\gencylevts
\quad\text{and}\quad
	\ev{\genX_s=\genX_r}
	\coloneqq \ev{\genX_{\pr{s,r}}\in\stsp^2_=}
	\in\gencylevts.
\end{equation*}
This means that \(\tupjumps_u\pr{\genX_u}\) is a \(\gencylevts\)-simple and therefore trivially \(\sigma\pr{\gencylevts}/\borel\pr{\extreals}\)-measurable variable.
Consequently, we also have that
\begin{equation*}
	\ev{\tupjumps_u\pr{\genX_u}\leq \alpha}
	\coloneqq \st{\pth\in\genpaths\colon \tupjumps_u\pr{\pth\pr{u}}\leq \alpha}
	\in\gencylevts
	\quad\text{for all } \alpha\in\reals.
\end{equation*}

Crucial to our proof of Theorem~\ref{the:precursor for cadlag paths} is the number of jumps of~\(\genX_{\noarg}\) in~\(\cci{s}{r}\cap\countpoints\), with \(\countpoints\) a countable subset of~\(\timepoints\).
For any countable subset~\(\countpoints\) of~\(\timepoints\) and \(s,r\in\reals\) such that \(s<r\) and \(\cci{s}{r}\cap\densepoints\neq\emptyset\), the number of jumps of~\(\genX_{\noarg}\) in~\(\cci{s}{r}\cap\densepoints\) is the variable
\begin{equation*}
	\pthjumps_{\cci{s}{r}\cap\densepoints}\pr{\genX_{\noarg}}
	\colon \genpaths\to\nnegints\cup\st{+\infty}
	\colon \pth\mapsto
	\pthjumps_{\cci{s}{r}\cap\densepoints}\pr{\pth};
\end{equation*}
in other words, \(\pthjumps_{\cci{s}{r}\cap\densepoints}\pr{\genX_{\noarg}}\) is the functional composition of~\(\pthjumps_{\cci{s}{r}\cap\densepoints}\) and the `projection'~\(\genX_{\noarg}\) which maps any~\(\pth\in\genpaths\) to itself.
This variable is \(\sigma\pr{\gencylevts}/\borel\pr{\extreals}\)-measurable.
\begin{lemma}
\label{lem:jumps pointwise convergence}
	Consider a countable subset~\(\countpoints\) of~\(\timepoints\) and some~\(s,r\in \reals\) such that \(s<r\) and \(\cci{s}{r}\cap\countpoints\neq\emptyset\).
	Then there is a sequence \(\pr{u_n}_{n\in\nats}\in \setoftseq_{\cci{s}{r}\cap\countpoints}\) such that \(\pr{\tupjumps_{u_n}\pr{\genX_{u_n}}}_{n\in\nats}\) is an increasing sequence of positive \(\gencylevts\)-simple -- and therefore \(\sigma\pr{\gencylevts}/\borel\pr{\extreals}\)-measurable -- variables that converges point-wise to~\(\pthjumps_{\cci{s}{r}\cap\countpoints}\pr{\genX_\noarg}\).
	Hence, \(\pthjumps_{\cci{s}{r}\cap\countpoints}\pr{\genX_{\noarg}}\) is \(\sigma\pr{\gencylevts}/\borel\pr{\extreals}\)-measurable.
\end{lemma}
In our proof of Lemma~\ref{lem:jumps pointwise convergence}, we make use of the straightforward observation that if~\(u\sqsubseteq v\), then the number of jumps of~\(\genX_{\noarg}\) along~\(v\) is greater than or equal to the number of jumps of~\(\genX_{\noarg}\) along~\(u\); clearly, it suffices to state this property for paths.
\begin{lemma}
\label{lem:jumps is monotone}
	For all \(\allpth\in\allpaths\) and \(u,v\in \setoftseq\) such that \(u\sqsubseteq v\), \(\tupjumps_u\pr{\allpth\pr{u}}\leq \tupjumps_v\pr{\allpth\pr{v}}\).
\end{lemma}
\begin{proof}
	Let us write \(u=\pr{r_1, \dots, r_n}\) and \(v=\pr{s_1, \dots, s_m}\).
	If \(n=1\), then the inequality in the statement holds trivially because \(\tupjumps_u=0\) and \(\tupjumps_v\geq 0\).
	Hence, we assume that \(n>1\).
	Because \(u\sqsubseteq v\) by the conditions of the statement, for all \(k\in \st{1, \dots, n}\) there is a unique natural number~\(\ell_k\in \st{1, \dots, m}\) such that \(r_k=s_{\ell_k}\).
	Then
	\begin{equation*}
		\tupjumps_u\pr{\allpth\pr{u}}
		= \card[\big]{\st[\big]{k\in\st{2, \dots, n}\colon \allpth\pr{s_{\ell_{k-1}}}
		\neq \allpth\pr{s_{\ell_k}}}}.
	\end{equation*}
	Because the tuples of time points~\(u\) and \(v\) are increasing, \(\ell_{k-1}<\ell_k\) for all \(k\in \st{2, \dots, n}\).
	Furthermore, it is easy to verify that for all \(k\in\st{2, \dots, n}\), \(\allpth\pr{s_{\ell_{k-1}}}\neq \allpth\pr{s_{\ell_k}}\) implies that \(\allpth\pr{s_{\ell-1}}\neq \allpth\pr{s_\ell}\) for at least one~\(\ell\in\st{\ell_{k-1}+1, \dots, \ell_k}\) -- if this is not the case, then \(\allpth\pr{s_{\ell_{k-1}}}=\allpth\pr{s_{\ell_{k-1}+1}}=\cdots=\allpth\pr{s_{\ell_k-1}}=\allpth\pr{s_{\ell_k}}\), which is a clear contradiction.
	Hence,
	\begin{equation*}
		\card[\big]{\st[\big]{k\in\st{2, \dots, n}\colon \allpth\pr{s_{\ell_{k-1}}}
		\neq \allpth\pr{s_{\ell_k}}}}
		\leq \card[\big]{\st[\big]{\ell\in\st{2, \dots, m}\colon \allpth\pr{s_{\ell-1}}
		\neq \allpth\pr{s_\ell}}},
	\end{equation*}
	and therefore
	\begin{equation*}
		\tupjumps_u\pr{\allpth\pr{u}}
		\leq \tupjumps_v\pr{\allpth\pr{v}}.
	\end{equation*}
\end{proof}
\begin{proof}[Proof of Lemma~\ref{lem:jumps pointwise convergence}]
	Because \(\countpoints\) is countable by assumption, so is~\(\cci{s}{r}\cap\countpoints\).
	Consequently, there is a sequence~\(\pr{u_n}_{n\in\nats}\in \setoftseq_{\cci{s}{r}\cap\countpoints}\) such that (i) \(u_n\sqsubseteq u_{n+1}\) for all \(n\in \nats\) and (ii) for all \(d\in \cci{s}{r}\cap\countpoints\), there is some \(n\in \nats\) such that \(d\) belongs to~\(u_n\); for example, simply add all the time points in~\(\cci{s}{r}\cap\countpoints\) one by one.

	For all \(n\) in \(\nats\), we know from before that \(\tupjumps_{u_n}\pr{\genX_{u_n}}\) is a positive \(\gencylevts\)-simple and therefore trivially \(\sigma\pr{\gencylevts}/\borel\pr{\extreals}\)-measurable variable, and because \(u_n\sqsubseteq u_{n+1}\) by construction, \(\tupjumps_{u_n}\pr{\genX_{u_n}}\leq \tupjumps_{u_{n+1}}\pr{\genX_{u_{n+1}}}\) due to Lemma~\ref{lem:jumps is monotone}.
	Hence, we have verified that \(\pr{\tupjumps_{u_n}\pr{\genX_{u_n}}}_{n\in\nats}\) is a sequence of positive \(\sigma\pr{\gencylevts}/\borel\pr{\extreals}\)-measurable variables that is monotonously increasing.
	Any monotonously increasing sequence converges point-wise, so it remains for us to show that this point-wise limit is~\(\pthjumps_{\cci{s}{r}\cap\countpoints}\pr{\genX_{\noarg}}\).
	On the one hand, for all \(n\in \nats\), \(u_n\) belongs to~\(\setoftseq_{\cci{s}{r}\cap\countpoints}\) so \(\tupjumps_{u_n}\pr{\genX_{u_n}}\leq \pthjumps_{\cci{s}{r}\cap\countpoints}\pr{\genX_\noarg}\); hence,
	\begin{equation*}
		\lim_{n\to+\infty} \tupjumps_{u_n}{\pr{\genX_{u_n}}}
		\leq \pthjumps_{\cci{s}{r}\cap\countpoints}\pr{\genX_\noarg}.
	\end{equation*}
	On the other hand, for any tuple of time points~\(u\in \setoftseq_{\countpoints}^{\cci{s}{r}}\) there is some~\(n\in \nats\) such that \(u\sqsubseteq u_n\), and therefore
	\begin{equation*}
		\pthjumps_{\cci{s}{r}\cap\countpoints}\pr{\genX_{\noarg}}
		\leq \lim_{n\to+\infty} \tupjumps_{u_n}\pr{\genX_{u_n}}.
	\end{equation*}
	From these two inequalities, we infer that \(\pr{\tupjumps_{u_n}\pr{\genX_{u_n}}}_{n\in\nats}\) converges point-wise to~\(\pthjumps_{\cci{s}{r}\cap\countpoints}\pr{\genX_{\noarg}}\).
	Since \(\pr{\tupjumps_{u_n}\pr{\genX_{u_n}}}_{n\in\nats}\) is an increasing sequence of positive \(\sigma\pr{\gencylevts}/\borel\pr{\extreals}\)-measurable variables, this implies that \(\pthjumps_{\cci{s}{r}\cap\countpoints}\pr{\genX_{\noarg}}\) is \(\sigma\pr{\gencylevts}/\borel\pr{\extreals}\)-measurable \cite[Chapter~2, Proposition~11 or Lemma~13]{1997Fristedt-Modern}.
\end{proof}
We end this intermezzo with an alternative expression for the probability that the number of jumps in~\(\cci{s}{r}\cap\countpoints\), with \(\countpoints\) a countable subset of~\(\timepoints\), is infinite.
\begin{lemma}
\label{lem:prob vs sup}
	Consider a stochastic process~\(\pr{\genpaths, \sigma\pr{\gencylevts}, \prob}\).
	Fix a countable subset~\(\countpoints\) of~\(\timepoints\).
	Then for all~\(s,r\in\reals\) such that \(s<r\) and \(\cci{s}{r}\cap\countpoints\neq\emptyset\),
	\begin{equation*}
		\prob\pr{\ev{\pthjumps_{\cci{s}{r}\cap\countpoints}\pr{\genX_{\noarg}}=+\infty}}
		= \lim_{k\to+\infty} \sup\st[\big]{\prob\pr{\ev{\tupjumps_u\pr{\genX_u}\geq k}}\colon u\in\setoftseq_{\cci{s}{r}\cap\countpoints}}.
	\end{equation*}
\end{lemma}
\begin{proof}
	Clearly, the sequence of level sets~\(\pr{\ev{\pthjumps_{\cci{s}{r}\cap\countpoints}\pr{\genX_{\noarg}}\geq k}}_{k\in\nats}\) decreases to the level set~\(\ev{\pthjumps_{\cci{s}{r}\cap\countpoints}\pr{\genX_{\noarg}}=+\infty}\).
	All these sets furthermore belong to $\sigma\pr{\gencylevts}$ because $\pthjumps_{\cci{s}{r}\cap\countpoints}\pr{\genX_{\noarg}}$ is \(\sigma\pr{\gencylevts}/\borel\pr{\extreals}\)-measurable due to Lemma~\ref{lem:jumps pointwise convergence}.
	Since the probability measure~\(\prob\) is continuous with respect to monotone sequences, it follows that
	\begin{equation*}
		\prob\pr{\ev{\pthjumps_{\cci{s}{r}\cap\countpoints}\pr{\genX_{\noarg}}=+\infty}}
		= \lim_{k\to+\infty} \prob\pr{\ev{\pthjumps_{\cci{s}{r}\cap\countpoints}\pr{\genX_{\noarg}}\geq k}}.
	\end{equation*}
	The equality in the statement follows if we can prove that for all \(k\in\nats\),
	\begin{equation}
	\label{eqn:proof of prob vs supprobs}
		\prob\pr{\ev{\pthjumps_{\cci{s}{r}\cap\countpoints}\pr{\genX_{\noarg}}\geq k}}
		= \sup\st[\big]{\prob\pr{\ev{\tupjumps_u\pr{\genX_u}\geq k}}\colon u\in\setoftseq_{\cci{s}{r}\cap\countpoints}}.
	\end{equation}
	To this end, we fix any \(k\in\nats\).
	By definition of~\(\pthjumps_{\cci{s}{r}\cap\countpoints}\pr{\genX_{\noarg}}\), \(\tupjumps_u\pr{\genX_u}\leq \pthjumps_{\cci{s}{r}\cap\countpoints}\pr{\genX_{\noarg}}\) for all \(u\in\setoftseq_{\cci{s}{r}\cap\countpoints}\).
	Hence, it follows from the monotonicity of the probability measure~\(\prob\) that
	\begin{equation*}
		\prob\pr{\ev{\pthjumps_{\cci{s}{r}\cap\countpoints}\pr{\genX_{\noarg}}\geq k}}
		\geq \prob\pr{\ev{\tupjumps_u\pr{\genX_u}\geq k}}
		\quad\text{for all } u\in\setoftseq_{\cci{s}{r} \cap\countpoints},
	\end{equation*}
	and therefore
	\begin{equation*}
		\prob\pr{\ev{\pthjumps_{\cci{s}{r}\cap\countpoints}\pr{\genX_{\noarg}}\geq k}}
		\geq
		\sup\st[\big]{\prob\pr{\ev{\tupjumps_u\pr{\genX_u}\geq k}}\colon u\in\setoftseq_{\cci{s}{r}\cap\countpoints}}.
	\end{equation*}
	Recall from Lemma~\ref{lem:jumps pointwise convergence} that there is a sequence~\(\pr{u_n}_{n\in\nats}\) in~\(\setoftseq_{\cci{s}{r}\cap\countpoints}\) such that \(\pr{\tupjumps_{u_n}\pr{\genX_{u_n}}}_{n\in\nats}\) is an increasing sequence sequence of positive \(\sigma\pr{\gencylevts}/\borel\pr{\extreals}\)-measurable variables that converges point-wise to~\(\pthjumps_{\cci{s}{r}\cap\countpoints}\pr{\genX_{\noarg}}\).
	Since all of the involved variables take values in $\mathbb{Z}_{\geq0}\cup\{+\infty\}$, this implies that the corresponding sequence of level sets~\(\pr{\ev{\tupjumps_{u_n}\pr{\genX_{u_n}}\geq k}}_{n\in\nats}\) increases to the level set~\(\ev{\pthjumps_{\cci{s}{r}\cap\countpoints}\pr{\genX_{\noarg}}\geq k}\).
	Again, it therefore follows from the continuity of the probability measure~\(\prob\) with respect to monotone sequences that
	\begin{align*}
		\prob\pr{\ev{\pthjumps_{\cci{s}{r}\cap\countpoints}\pr{\genX_{\noarg}}\geq k}}
		&= \lim_{n\to+\infty} \prob\pr{\ev{\tupjumps_{u_n}\pr{\genX_{u_n}}\geq k}} \\
		&\leq \sup\st[\big]{\prob\pr{\ev{\tupjumps_u\pr{\genX_u}\geq k}\colon u\in\setoftseq_{\cci{s}{r}\cap\countpoints}}}.
	\end{align*}
	Eqn.~\eqref{eqn:proof of prob vs supprobs} now follows immediately from the preceding two inequalities.
\end{proof}

\subsection{Modifying a path with a finite number of jumps along a countable dense subset}
\label{assec:modification}
We are interested in the number of jumps of~\(\allX_{\noarg}\) in \(\cci{s}{r}\cap\densepoints\) because if this is finite for some~\(\allpth\in \allpaths\), then the right-sided limit of~\(\allpth\) along~\(\densepoints\) exists in any right-sided limit point~\(t\) of~\(\densepoints\cap\coi{s}{r}\), and similarly for the left-sided limit in the left-sided limit points.
The following intermediary result -- which can be seen as a version of Theorem~62.7 in \cite[Chapter~II]{1994Rogers-Diffusions} in our setting -- essentially shows this, and also establishes a sufficient condition under which we can `modify' the path~\(\allpth\) in such a way that it becomes càdlàg.
\begin{lemma}
\label{lem:limit along densetpoints}
	Consider some \(\allpth\in\allpaths\), and fix a subset~\(\densepoints\) of~\(\timepoints\) such that \(\timepoints\subseteq\densepoints\cup\rlims\pr{\densepoints}\).
	Suppose that for all \(n\in\nats\) such that \(\cci{-n}{n}\cap\densepoints\neq\emptyset\),
	\begin{equation*}
		\pthjumps_{\cci{-n}{n}\cap\densepoints}\pr{\allpth}
		<+\infty.
	\end{equation*}
	Then for all \(t\in \rlims\pr{\densepoints}\), \(\lim_{\densepoints\ni r\searrow t} \allpth\pr{r}\) exists and for all \(t\in\llims\pr{\densepoints}\), \(\lim_{\densepoints\ni s\nearrow t} \allpth\pr{s}\) exists.
	Furthermore, there is a càdlàg path~\(\cadpth\in\cadpaths\) such that
	\begin{equation*}
		\cadpth\pr{t}
		= \begin{cases}
			\lim_{\densepoints\ni r\searrow t} \allpth\pr{r} &\text{if  } t\in\rlims\pr{\densepoints} \\
			\allpth\pr{t} &\text{if } t\in\densepoints\setminus\rlims\pr{\densepoints}
		\end{cases}
		\quad\text{for all } t\in\timepoints.
	\end{equation*}
\end{lemma}
\begin{proof}
	In the first part of this proof, we prove a convenient intermediary result.
	Fix any time point~\(t\in\reals\), and let \(n\) be a natural number such that \(-n<t<n\) and \(\cci{-n}{n}\cap\densepoints\neq\emptyset\) -- this is always possible because \(\densepoints\) is dense in the non-empty subset~\(\timepoints\) of~\(\timeaxis\).
	By the condition in the statement, there is some tuple of time points~\(u=\pr{t_1, \dots, t_m}\in \setoftseq_{\cci{-n}{n}\cap\densepoints}\) such that
	\begin{equation*}
		\tupjumps_u\pr{\allpth\pr{u}}
		= \pthjumps_{\cci{-n}{n}\cap\densepoints}\pr{\allpth}
		= \sup\st[\Big]{\tupjumps_v\pr{\allpth\pr{v}}\colon v\in\setoftseq_{\cci{-n}{n}\cap\densepoints}}
		< + \infty
	\end{equation*}
	and, for any tuple of time points~\(v\in \setoftseq_{\cci{-n}{n}\cap\densepoints}\) such that \(u\sqsubseteq v\),
	\begin{equation*}
		\tupjumps_u\pr{\allpth\pr{u}}
		\leq \tupjumps_v\pr{\allpth\pr{v}}
		\leq
		\pthjumps_{\cci{-n}{n}\cap\densepoints}\pr{\allpth}
		= \sup\st[\Big]{\tupjumps_w\pr{\allpth\pr{w}}\colon w\in\setoftseq_{\cci{-n}{n}\cap\densepoints}}
		= \tupjumps_u\pr{\allpth\pr{u}},
	\end{equation*}
	where for the first inequality we used Lemma~\ref{lem:jumps is monotone}.
	In other words, no matter how many time points from~\(\cci{-n}{n}\cap\densepoints\) we add to~\(u\), the number of jumps of~\(\allpth\) along the sequence will remain the same.
	From this, we infer that \(\allpth\pr{s}=\allpth\pr{t_1}\) for all \(s\in \cci{-n}{t_1}\cap\densepoints\), \(\allpth\pr{r}=\allpth\pr{t_m}\) for all \(r\in \cci{t_m}{n}\cap\densepoints\) and that for all \(k\in \st{1, \dots, m-1}\), there is some \(t^\star_k\in \cci{t_k}{t_{k+1}}\) such that \(\allpth\pr{s}=\allpth\pr{t_k}\) for all \(s\in \coi{t_k}{t^\star_k}\cap\densepoints\) and \(\allpth\pr{r}=\allpth\pr{t_{k+1}}\) for all \(r\in \oci{t^\star_k}{t_{k+1}}\cap\densepoints\).
	Because \(-n<t<n\), we infer from this that there are strictly positive real numbers~\(\delta_{t,+},\delta_{t,-}\) and states~\(x_{t,+}, x_{t,-}\in \stsp\) such that
	\begin{equation}
	\label{eqn:proof of limit along densepoints:two sides}
		\pr[\big]{\forall r\in\ooi{t}{t+\delta_{t,+}}\cap\densepoints}~\allpth\pr{r}=x_{t,+}
		\quad\text{and}\quad
		\pr[\big]{\forall s\in\ooi{t-\delta_{t,-}}{t}\cap\densepoints}~\allpth\pr{s}=x_{t,-}.
	\end{equation}

	In the second part of this proof, we use Eqn.~\eqref{eqn:proof of limit along densepoints:two sides} to prove the first part of the statement.
	Fix a right-sided limit point~\(t\) of~\(\densepoints\).
	From the first part of this proof -- that is, from right before Eqn.~\eqref{eqn:proof of limit along densepoints:two sides} -- we know that there are some positive real number~\(\delta_{t,+}\) and state~\(x_{t,+}\) such that for all \(r\in \ooi{t}{t+\delta_{t,+}}\cap\densepoints\), \(\allpth\pr{r}=x_{t,+}\).
	Then \(\lim_{\densepoints\ni r\searrow t}\allpth\pr{r}\) exists, and is equal to~\(x_{t,+}\).
	Similarly, for any left-sided limit point~\(t\) of~\(\densepoints\), there are some positive real number~\(\delta_{t,-}\) and state~\(x_{t,-}\) such that for all \(s\in \ooi{t-\delta_{t,-}}{t}\cap\densepoints\), \(\allpth\pr{r}=x_{t,-}\); so \(\lim_{\densepoints\ni s\nearrow t}\allpth\pr{s}\) exists and is equal to~\(x_{t,-}\).

	In the third part of this proof, we show the existence of the càdlàg path~\(\cadpth\in\cadpaths\).
	We start by constructing a map $\vpth\colon\timepoints\to\stsp$, defined by
\begin{equation*}
		\vpth\pr{t}
		\coloneqq \begin{cases}
			\lim_{\densepoints\ni r\searrow t} \allpth\pr{r} &\text{if  } t\in\rlims\pr{\densepoints} \\
			\allpth\pr{t} &\text{if } t\in\densepoints\setminus\rlims\pr{\densepoints}
		\end{cases}
		\quad\text{for all } t\in\timepoints.
	\end{equation*}
In the remainder of this proof, we will show that $\vpth$ is càdlàg. It then follows from Lemma~\ref{lem:extension of cadlag path} that $\vpth$ can be extended to a map~\(\rpth\colon\reals\to\stsp\) that is càdlàg. Since restrictions of càdlàg maps are càdlàg too, this implies that $\cadpth\coloneqq\rpth\vert_{\timeaxis}$ is a càdlàg path in $\cadpaths$ that, since it extends $\vpth$, clearly satisfies the condition in the statement.

To show that $\vpth$ is càdlàg, we consider any $t\in\reals$. We've shown earlier in this proof that there are strictly positive real numbers~\(\delta_{t,+},\delta_{t,-}\) and states~\(x_{t,+}, x_{t,-}\in \stsp\) that satisfy Eqn.~\eqref{eqn:proof of limit along densepoints:two sides}. Taking into account the definition of $\vpth$, this implies that
\begin{equation*}
		\pr[\big]{\forall r\in\ooi{t}{t+\delta_{t,+}}\cap\timepoints}~\vpth\pr{r}=x_{t,+}
		\quad\text{and}\quad
		\pr[\big]{\forall s\in\ooi{t-\delta_{t,-}}{t}\cap\timepoints}~\vpth\pr{s}=x_{t,-}.
	\end{equation*}
Since this is true for any $t\in\reals$, and therefore definitely for any $t\in\llims(\timepoints)$ or $t\in\rlims(\timepoints)$ this already establishes the first two conditions in Definition~\ref{def:cadlag}.
To establish the third condition, we need to show that if $t\in\timepoints\cap\rlims(\timepoints)$, then also $\vpth(t)=x_{t,+}$. So consider the case $t\in\timepoints\cap\rlims(\timepoints)$. Since $t\in\rlims(\timepoints)$, we also have that $t\in\rlims(\densepoints)$.
To see why that is the case, consider any $\delta\in\posreals$. Since $t\in\rlims(\timepoints)$, we know that \(\ooi{t}{t+\delta}\cap\timepoints\neq\emptyset\).
Choose any $t_\delta\in\ooi{t}{t+\delta}\cap\timepoints$.
Since $t_\delta\in\timepoints\subseteq\densepoints\cup\rlims\pr{\densepoints}$, we then have that $\coi{t_\delta}{t+\delta}\cap\densepoints\neq\emptyset$, which implies that also $\ooi{t}{t+\delta}\cap\densepoints\neq\emptyset$.
Since $\delta\in\posreals$ was arbitrary, we indeed have that $t\in\rlims(\densepoints)$.
It therefore follows from Eqn.~\eqref{eqn:proof of limit along densepoints:two sides} that $\vpth(t)=x_{t,+}$.
\end{proof}

\subsection{Two additional intermediary results}
\label{assec:two additional intermediary results}
In our proof for Theorem~\ref{the:precursor for cadlag paths}, we need the following continuity property of probability measures -- for a proof, see for example Theorem~2 in~\cite[Section~6.1]{1997Fristedt-Modern}.
\begin{lemma}
\label{lem:limit of events}
	Consider a probability space~\(\pr{\posssp, \algebra, \prob}\) and a sequence~\(\pr{A_n}_{n\in\nats}\) in~\(\algebra\).
	If \(\lim_{n\to+\infty} \indica{A_n}\pr{\sample}\) exists for all \(\sample\in\posssp\), then
	\begin{equation*}
		\lim_{n\to+\infty} A_n
		\coloneqq \st[\Big]{\sample\in\posssp\colon \lim_{n\to+\infty}\indica{A_n}\pr{\sample}=1}\in\algebra
	\end{equation*}
	and
	\begin{equation*}
		\prob\pr[\Big]{\lim_{n\to+\infty} A_n}
		= \lim_{n\to+\infty} \prob\pr{A_n}.
	\end{equation*}
\end{lemma}

Another example where Lemma~\ref{lem:limit of events} comes in handy is the following intermediary result, which is the final lemma we need to prove Theorem~\ref{the:precursor for cadlag paths}.
To state it, we observe that for all \(n\in\nats\) and \(v,w\in\setoftseq\) such that \(v=\pr{v_1, \dots, v_m}\) and \(w=\pr{w_1, \dots, w_m}\),
\begin{equation*}
	\ev{\genX_v=\genX_w}
	\coloneqq \bigcap_{j=1}^m \ev{\genX_{v_j}=\genX_{w_j}}
	\in\gencylevts.
\end{equation*}
\begin{lemma}
\label{lem:approximating sequence of sequence of time points}
	Consider a probability charge~\(\prob\) on a domain~\(\genalgebra\) that includes~\(\gencylevts\).
	Fix some \(v=\pr{t_1, \dots, t_m}\in\setoftseq\) and some countable subset~\(\densepoints\) of~\(\timepoints\) such that \(\st{t_1, \dots, t_m}\subseteq\densepoints\cup\rlims\pr{\densepoints}\).
	Suppose the collection~\(\dist_{\noarg}\) of finite-dimensional charges of~\(\prob\) satisfies~\ref{def:dist regular:cont}.
	Then for all \(\epsilon\in\posreals\), there is some~\(w=\pr{t'_1, \dots, t'_m}\in\setoftseq_{\densepoints}\) such that
	\begin{equation*}
		\abs{\prob\pr{A}- \prob\pr{A\cap\ev{\genX_v=\genX_w}}}
		< \epsilon
		\quad\text{for all } A\in\genalgebra.
	\end{equation*}
	Furthermore, for any $\delta\in\posreals$, we can guarantee that $t'_j\in\ooi{t_j}{t_j+\delta}$ if $t_j\in\rlims(\densepoints)$ and $t'_j=t_j$ otherwise.
	Hence, there is some sequence \(\pr{w_\ell}_{\ell\in\nats}\) in~\(\setoftseq_{\densepoints}\) such that (i) for all \(\ell\in\nats\), \(w_{\ell}=\pr{t^\ell_1, \dots, t^\ell_m}\), (ii) for all \(j\in\st{1, \dots, m}\), \(\pr{t^\ell_j}_{\ell\in\nats}\) is a strictly decreasing sequence in~\(\densepoints\) that converges to~\(t_j\) if \(t_j\in\rlims\pr{\densepoints}\) and a constant sequence that is equal to $t_j$ otherwise, and (iii)
	\begin{equation*}
		\prob\pr{A}
		= \lim_{\ell\to+\infty} \prob\pr{A\cap \ev{\genX_v=\genX_{w_\ell}}}
	\end{equation*}
\end{lemma}
\begin{proof}
	Fix some \(j\in\st{1, \dots, m}\).
	If \(t_j\in\densepoints\setminus\rlims\pr{\densepoints}\), then we let \(t'_j\coloneqq t_j\), and then trivially
	\begin{equation*}
		\prob\pr{\ev{\genX_{t_j}=\genX_{t'_j}}}
		= 1
		> 1-\frac{\epsilon}{m}.
	\end{equation*}
	If on the other hand \(t_j\in\rlims\pr{\densepoints}\), then it follows from the assumptions on~\(\prob\) and \(\densepoints\) in the statement that there is some \(t'_j\) in~\(\ooi{t_j}{t_{j+1}}\cap\ooi{t_j}{t_j+\delta}\cap\densepoints\) -- where we let \(t_{m+1}\coloneqq+\infty\) -- such that
	\begin{equation*}
		\prob\pr{\ev{\genX_{t_j}=\genX_{t'_j}}}
		= \dist_{\pr{t_j, t'_j}}\pr{\stsp^2_{=}}
		> 1-\frac{\epsilon}{m}.
	\end{equation*}

	Let \(w\coloneqq\pr{t'_1, \dots, t'_m}\in\setoftseq_{\densepoints}\), and note that
	\begin{equation*}
		\ev{\genX_v=\genX_w}
		= \bigcap_{j=1}^m \ev{\genX_{t_j}=\genX_{t'_j}}
		\quad\text{and}\quad
		\ev{\genX_v=\genX_w}^{\compl}
		= \bigcup_{j=1}^m \ev{\genX_{t_j}=\genX_{t'_j}}^{\compl}.
	\end{equation*}
	It follows from this, the properties of the probability charge~\(\prob\) and our construction of~\(w\) that
	\begin{align*}
		\prob\pr{\ev{\genX_v=\genX_w}^{\compl}}
		= \prob\pr*{\bigcup_{j=1}^m \ev{\genX_{t_j}=\genX_{t'_j}}^{\compl}}
		&\leq \sum_{j=1}^m \prob\pr{\ev{\genX_{t_j}=\genX_{t'_j}}^{\compl}} \\
		&= \sum_{j=1}^m \pr[\big]{1-\prob\pr{\ev{\genX_{t_j}=\genX_{t'_j}}}} \\
		&<\sum_{j=1}^m\frac{\epsilon}{m}=\epsilon.
	\end{align*}
	Now for all \(A\in\genalgebra\),
	\begin{equation*}
		\prob\pr{A}
		= \prob\pr{A\cap\ev{\genX_v=\genX_w}} + \prob\pr{A\cap\ev{\genX_v=\genX_w}^{\compl}},
	\end{equation*}
	and therefore
	\begin{align*}
		\abs{\prob\pr{A}-\prob\pr{A\cap\ev{\genX_v=\genX_w}}}
		= \prob\pr{A\cap\ev{\genX_v=\genX_w}^{\compl}}
		\leq \prob\pr{\ev{\genX_v=\genX_w}^{\compl}}< \epsilon,
	\end{align*}
	as required.

	The second part of the statement follows immediately from the first part, since it is easy to ensure that the sequences~\(\pr{t^\ell_j}_{\ell\in\nats}\) are decreasing if $t_j\in\rlims(\densepoints)$: when constructing $w_{\ell+1}$, simply choose $\delta$ small enough such that $t_j+\delta<t_j^\ell$ for all $j\in\{1,\ldots,m\}$ such that $t_j\in\rlims(\densepoints)$.
\end{proof}

\subsection{Proof for Theorem~\ref{the:precursor for cadlag paths}}
\label{asec:actual proof}
Finally, we can get around to proving the precursor to our main result.
\begin{proof}[Proof of Theorem~\ref{the:precursor for cadlag paths}]
	First, we prove the necessity; that is, we assume that the probability charge~\(\prob\) is countably additive, and show that this implies that \(\dist_{\noarg}\) is regular.
	Since \(\prob\) is countably additive, we know from Caratheodory's Extension Theorem that there is a unique probability measure~\(\prob_\sigma\) on~\(\sigma\pr{\cadcylevts}\) that extends~\(\prob\).
	Note that the finite-dimensional distributions of~\(\prob_\sigma\) are~\(\dist_{\noarg}\).
	We will use this to show that \(\dist_{\noarg}\) is regular.

	For \ref{def:dist regular:cont}, we fix some \(t\in\rlims\pr{\timepoints}\cap\timepoints\).
	For any decreasing sequence~\(\pr{r_n}_{n\in\nats}\) in~\(\ooi{t}{+\infty}\cap{\timepoints}\) with \(\lim_{n\to+\infty}r_n=t\), \(\lim_{n\to+\infty}\pth\pr{r_n}=\pth\pr{t}\) for all \(\pth\in\cadpaths\) due to the right-continuity of càdlàg paths, and therefore $\lim_{n\to +\infty}\ev{\cadX_{t}=\cadX_{r_n}}=\cadpaths$; hence, it follows from Lemma~\ref{lem:limit of events} that
	\begin{equation*}
		\lim_{n\to+\infty} \prob_{\sigma}\pr{\ev{\cadX_{t}=\cadX_{r_n}}}
		= \prob_{\sigma}\pr{\cadpaths}
		= 1.
	\end{equation*}
	Since this is true for any such sequence, it follows that
	\begin{equation*}
		\lim_{\timepoints\ni r\searrow t} \prob_{\sigma}\pr{\ev{\cadX_t=\cadX_r}}
		= 1.
	\end{equation*}
	Furthermore, for all $r\in\timepoints$ such that $r>t$, we have that
	\begin{align*}
	\prob_{\sigma}\pr{\ev{\cadX_t=\cadX_r}}
	=\prob\pr{\ev{\cadX_t=\cadX_r}}
	=\prob\pr{\ev{\cadX_{(t,r)}\in \stsp^2_=}}
	=\dist_{\pr{t,r}}\pr{\stsp^2_{=}}.
	\end{align*}
	From these two observations, we infer that
	$\lim_{\timepoints\ni r\searrow t} \dist_{\pr{t,r}}\pr{\stsp^2_{=}}	= 1$, as required for \ref{def:dist regular:cont}.

	For \ref{def:dist regular:bound on exp of jumps}, we fix some \(n\in\nats\) such that \(\cci{-n}{n}\cap\timepoints\neq\emptyset\), and some countable subset~\(\countpoints\) of~\(\cci{-n}{n}\cap\timepoints\) such that \(\cci{-n}{n}\cap\timepoints\subseteq\countpoints\cup\rlims\pr{\countpoints}\); this is always possible due to Lemma~\ref{lem:subset of sorgenfrey is seperable}.
	Since \(\countpoints\subseteq\cci{-n}{n}\cap\timeaxis\), we know from Proposition~\ref{prop:number of jumps is finite} that
	\begin{equation*}
		\pthjumps_{\countpoints}\pr{\cadpth}
		\leq \pthjumps_{\cci{-n}{n}\cap\timeaxis}\pr{\cadpth}
		< +\infty
		\quad\text{for all } \cadpth\in\cadpaths.
	\end{equation*}
	Hence, \(\ev{\pthjumps_{\countpoints}\pr{\cadX_{\noarg}}=+\infty}=\emptyset\), and therefore it must be that
	\begin{equation*}
		\prob_{\sigma}\pr{\ev{\pthjumps_{\countpoints}\pr{\cadX_{\noarg}}=+\infty}}
		= \prob_{\sigma}\pr{\emptyset}
		= 0.
	\end{equation*}
	From this and Lemma~\ref{lem:prob vs sup}, it follows that
	\begin{equation}
	\label{eqn:proof of main:R2 is necessary}
		\lim_{k\to+\infty} \sup\st[\big]{\prob_{\sigma}\pr{\ev{\tupjumps_u\pr{\cadX_u}\geq k}}\colon u\in\setoftseq_{\countpoints}}
		= 0.
	\end{equation}
	Now fix any \(k\in\nats\) and \(\epsilon\in\posreals\).
	Then for all \(v\in\setoftseq_{\cci{-n}{n}\cap\timepoints}\), we know from Lemma~\ref{lem:approximating sequence of sequence of time points} -- with \(\countpoints\) and \(\ev{\tupjumps_v\pr{\cadX_v}\geq k}\) here in the role of \(\densepoints\) and \(A\) there -- that there is some \(w\in\setoftseq_{\countpoints}\) such that
	\begin{equation*}
		\prob_{\sigma}\pr{\ev{\tupjumps_v\pr{\cadX_v}\geq k}}
		< \prob_{\sigma}\pr{\ev{\tupjumps_v\pr{\cadX_v}\geq k}\cap\ev{\cadX_v=\cadX_w}}+\epsilon
		\leq \prob_{\sigma}\pr{\ev{\tupjumps_w\pr{\cadX_w}\geq k}}+\epsilon,
	\end{equation*}
	where for the non-strict inequality we used the monotonicity of~\(\prob_{\sigma}\) and that
	\begin{equation*}
		\ev{\tupjumps_v\pr{\cadX_v}\geq k}\cap\ev{\cadX_v=\cadX_w}
		\subseteq \ev{\tupjumps_w\pr{\cadX_w}\geq k}.
	\end{equation*}
	From this inequality, we infer that
	\begin{equation*}
		\sup\st[\big]{\prob_{\sigma}\pr{\ev{\tupjumps_u\pr{\cadX_u}\geq k}}\colon u\in\setoftseq_{\cci{-n}{n}\cap\timepoints}}
		\leq
		\sup\st[\big]{\prob_{\sigma}\pr{\ev{\tupjumps_u\pr{\cadX_u}\geq k}}\colon u\in\setoftseq_{\countpoints}}+\epsilon.
	\end{equation*}
	Since \(\epsilon\in\posreals\) is arbitrary, and $\countpoints\subseteq\cci{-n}{n}\cap\timepoints$, it follows that
	\begin{equation*}
		\sup\st[\big]{\prob_{\sigma}\pr{\ev{\tupjumps_u\pr{\cadX_u}\geq k}}\colon u\in\setoftseq_{\cci{-n}{n}\cap\timepoints}}
		=
		\sup\st[\big]{\prob_{\sigma}\pr{\ev{\tupjumps_u\pr{\cadX_u}\geq k}}\colon u\in\setoftseq_{\countpoints}}.
	\end{equation*}
	This is true for all $k\in\nats$, so we infer from Equation~\eqref{eqn:proof of main:R2 is necessary} that
	\begin{equation*}
		\lim_{k\to+\infty} \sup\st[\big]{\prob_{\sigma}\pr{\ev{\tupjumps_u\pr{\cadX_u}\geq k}}\colon u\in\setoftseq_{\cci{-n}{n}\cap\timepoints}}
		= 0.
	\end{equation*}
	Because the finite-dimensional distributions of~\(\prob_{\sigma}\) are~\(\dist_\noarg\), we conclude from this that
	\begin{equation*}
		\lim_{k\to+\infty} \sup\st[\big]{\dist_u\pr{\st{x_u\in\stsp_u\colon\tupjumps_u\pr{x_u}\geq k}}\colon u\in\setoftseq_{\cci{-n}{n}\cap\timepoints}}
		= 0,
	\end{equation*}
	which is exactly what \ref{def:dist regular:bound on exp of jumps} demands.

	Second, we prove the sufficiency: we assume that \(\dist_{\noarg}\) is regular, and show that the corresponding probability charge~\(\prob\) is countably additive, or equivalently, that for any sequence~\(\pr{\mathring{A}_i}_{i\in\nats}\) in~\(\cadcylevts\) such that \(\mathring{A}_i\supseteq \mathring{A}_{i+1}\) for all \(i\in\nats\) and \(\bigcap_{i\in\nats} \mathring{A}_i=\emptyset\), \(\lim_{i\to+\infty} \prob\pr{\mathring{A}_i}=0\).

	To prove this, we fix any such sequence~\(\pr{\mathring{A}_i}_{i\in\nats}\).
	By definition of~\(\cadcylevts\), for all \(i\in\nats\) there are some \(u_i\in\setoftseq\) and \(A_i\in\wp\pr{\stsp_{u_i}}\) such that \(\mathring{A}_i=\ev{\cadX_{u_i}\in A_i}\).
	Without loss of generality, we may assume that \(u_i \sqsubseteq u_{i+1}\) -- if this is not the case, add the missing points to~\(u_{i+1}\) to obtain~\(u'_{i+1}\), and replace \(A_{i+1}\) by \(A'_{i+1}\coloneqq\st{x_{u'_{i+1}}\in\stsp_{u'_{i+1}}\colon x_{u_{i+1}}\in A_{i+1}}\).
	Then for all \(i\in\nats\), \(\st{x_{u_{i+1}}\in\stsp_{u_{i+1}}\colon x_{u_i}\in A_i}\supseteq A_{i+1}\) since \(\pr{\mathring{A}_i}_{i\in\nats}\) is decreasing.

	To prove that the probability of these events converges to~\(0\), we will rely on Theorem~\ref{the:Kolomogorov extension theorem}.
	Since \(\dist_{\noarg}\) is consistent, this result says that there is a (unique) probability measure~\(\allprob\) on~\(\sigma\pr{\allcylevts}\) such that
	\begin{equation}
	\label{eqn:proof of main:cond on allprob}
		\allprob\pr{\ev{\allX_u\in A}}
		= \dist_u\pr{A}
		= \prob\pr{\ev{\cadX_u\in A}}
		\quad\text{for all } u\in\setoftseq, A\in\wp\pr{\stsp_u}.
	\end{equation}
	It follows from these equalities that
	\begin{equation}
	\label{eqn:proof of main:intermed1}
		\lim_{i\to+\infty} \prob\pr{\mathring{A}_i}
		= \lim_{i\to+\infty} \prob\pr{\ev{\cadX_{u_i}\in A_i}}
		= \lim_{i\to+\infty} \allprob\pr{\ev{\allX_{u_i}\in A_i}}.
	\end{equation}
	Note that since \(\pr{\ev{\cadX_{u_i}\in A_i}}_{i\in\nats}\) is decreasing, the same must hold for \(\pr{\ev{\allX_{u_i}\in A_i}}_{i\in\nats}\).

	Fix some countable subset~\(\densepoints\) of~\(\timepoints\) such that \(\timepoints\subseteq\densepoints\cup\rlims\pr{\densepoints}\); this is always possible due to Lemma~\ref{lem:subset of sorgenfrey is seperable}.
	Then it follows from Eqn.~\eqref{eqn:proof of main:cond on allprob} and \ref{def:dist regular:cont} that for all \(t\in\timepoints\cap\rlims\pr{\timepoints}\)
	\begin{equation*}
		\lim_{\timepoints\ni r\searrow t} \allprob\pr{\ev{\allX_t=\allX_r}}
		= \lim_{\timepoints\ni r\searrow t} \dist_{\pr{t,r}}\pr{\stsp^2_{=}}
		= 1.
	\end{equation*}
	For the implication of Eqn.~\eqref{eqn:proof of main:cond on allprob} and \ref{def:dist regular:bound on exp of jumps}, we fix some \(n\in\nats\) such that \(\cci{-n}{n}\cap\densepoints\neq\emptyset\).
	Recall from Lemma~\ref{lem:prob vs sup} that
	\begin{equation*}
		\allprob\pr{\ev{\pthjumps_{\cci{-n}{n}\cap\densepoints}\pr{\allX_\noarg}=+\infty}}
		= \lim_{k\to+\infty} \sup\st[\big]{\allprob\pr{\ev{\tupjumps_u\pr{\allX_u}\geq k}}\colon u\in\setoftseq_{\cci{-n}{n}\cap\densepoints}}.
	\end{equation*}
	It therefore follows from Eqns.~\ref{def:dist regular:bound on exp of jumps} and~\eqref{eqn:proof of main:cond on allprob} that
	\begin{align*}
	0&=\lim_{k\to+\infty} \sup\st[\big]{\dist_u\pr{\st{x_u\in\stsp_u\colon\tupjumps_u\pr{x_u}\geq k}}\colon u\in\setoftseq_{\cci{-n}{n}\cap\timepoints}}\\
	&\geq\lim_{k\to+\infty} \sup\st[\big]{\dist_u\pr{\st{x_u\in\stsp_u\colon\tupjumps_u\pr{x_u}\geq k}}\colon u\in\setoftseq_{\cci{-n}{n}\cap\densepoints}}\\
		&= \lim_{k\to+\infty} \sup\st[\big]{\allprob\pr{\ev{\tupjumps_u\pr{\allX_u}\geq k}}\colon u\in\setoftseq_{\cci{-n}{n}\cap\densepoints}} \\
		&= \allprob\pr{\ev{\pthjumps_{\cci{-n}{n}\cap\densepoints}\pr{\allX_\noarg}=+\infty}}\geq0,
	\end{align*}
	which implies that $\allprob\pr{\ev{\pthjumps_{\cci{-n}{n}\cap\densepoints}\pr{\allX_\noarg}=+\infty}}=0$
	and therefore, that
	\begin{equation*}
		\allprob\pr{\ev{\pthjumps_{\cci{-n}{n}\cap\densepoints}\pr{\allX_\noarg}<+\infty}}
		= 1-\allprob\pr{\ev{\pthjumps_{\cci{-n}{n}\cap\densepoints}\pr{\allX_\noarg}=+\infty}}
		= 1.
	\end{equation*}

	Let \(\dot{A}\coloneqq\bigcap_{n=n^\star}^{+\infty}\ev{\pthjumps_{\cci{-n}{n}\cap\densepoints}\pr{\allX_{\noarg}}<+\infty}\), where \(n^\star\) is the smallest natural number such that \(\cci{-n}{n}\cap\densepoints\neq\emptyset\).
	Since each of the events~\(\ev{\pthjumps_{\cci{-n}{n}\cap\densepoints}\pr{\allX_{\noarg}}<+\infty}\) has probability~\(1\), so does their countable intersection~\(\dot{A}\).
	Hence,
	\begin{equation}
	\label{eqn:proof of main:intermed2}
		\lim_{i\to+\infty}\allprob\pr{\ev{\allX_{u_i}\in A_i}}
		= \lim_{i\to+\infty} \allprob\pr{\dot{A}\cap\ev{\allX_{u_i}\in A_i}}.
	\end{equation}
	From Lemma~\ref{lem:approximating sequence of sequence of time points}, we know that for all \(i\in\nats\) there is a sequence \(\pr{u_{i,\ell}}_{\ell\in\nats}\) in~\(\setoftseq_{\densepoints}\) such that (i) each of the \(u_{i,\ell}\)'s has the same number of time points as~\(u_i\), (ii) if the \(j\)-th component of~\(u_i\) is a right-sided limit point of~\(\densepoints\), then the \(j\)-th component of~\(u_{i,\ell}\) strictly decreases to it as~\(\ell\) recedes to~\(+\infty\), and if the \(j\)-th component of~\(u_i\) is not a right-sided limit point of~\(\densepoints\), the \(j\)-th component of~\(u_{i,\ell}\) is equal to it, and (iii)
	\begin{equation*}
		\allprob\pr{\dot{A}\cap\ev{\allX_{u_i}\in A_i}}
		= \lim_{\ell\to+\infty} \allprob\pr{\dot{A}\cap\ev{\allX_{u_i}\in A_i}\cap\ev{\allX_{u_i}=\allX_{u_{i,\ell}}}}.
	\end{equation*}
	Since for all $i,\ell\in\nats$,
\begin{align*}
\dot{A}\cap\ev{\allX_{u_i}\in A_i}\cap\ev{\allX_{u_i}=\allX_{u_{i,\ell}}}
&=
\dot{A}\cap\ev{\allX_{u_{i,\ell}}\in A_i}\cap\ev{\allX_{u_i}=\allX_{u_{i,\ell}}}\\
&\subseteq\dot{A}\cap\ev{\allX_{u_{i,\ell}}\in A_i},
\end{align*}
	we infer from this that for all \(i\in\nats\),
	\begin{align}
		\allprob\pr{\dot{A}\cap\ev{\allX_{u_i}\in A_i}}
		&\leq \limsup_{\ell\to+\infty} \allprob\pr{\dot{A}\cap\ev{\allX_{u_{i,\ell}}\in A_i}}.
	\label{eqn:proof of main:intermed3}
	\end{align}

	For all \(i\in\nats\) and \(\allpth\in\dot{A}\), it follows from the construction of~\(\pr{u_{i,\ell}}_{\ell\in\nats}\) and -- for the components of $u_i$ that are a right-sided limit point of $\densepoints$ -- Lemma~\ref{lem:limit along densetpoints} that \(\lim_{\ell\to+\infty}\allpth\pr{u_{i,\ell}}\) exists (with the limit taken component-wise); from this and Lemma~\ref{lem:limit of events}, it follows that for all $i\in\nats$
	\begin{equation*}
		\dot{A}_{i}
		\coloneqq
\{\allpth\in\dot{A}\colon\lim_{\ell\to+\infty}\allpth(u_{i,\ell})\in A_i\}
		=\lim_{\ell\to+\infty} \dot{A}\cap\ev{\allX_{u_{i,\ell}}\in A_i}
		\in\sigma\pr{\allcylevts}
	\end{equation*}
	and
	\begin{equation}
	\label{eqn:proof of main:intermed4}
		\allprob\pr{\dot{A}_i}
		=\lim_{\ell\to+\infty} \allprob\pr{\dot{A}\cap\ev{\allX_{u_{i,\ell}}\in A_i}}
		=\limsup_{\ell\to+\infty} \allprob\pr{\dot{A}\cap\ev{\allX_{u_{i,\ell}}\in A_i}}.
	\end{equation}
	We now establish two properties of these sets $A_i$.

	The first property is that \(\pr{\dot{A}_i}_{i\in\nats}\) is a decreasing sequence.
	To prove this, we consider any $i\in\nats$ and set out to show that $\dot{A}_{i+1}\subseteq\dot{A}_i$, or equivalently, that any $\allpth$ in $\dot{A}_{i+1}$ belongs to $\dot{A}_i$ as well. So consider any $\allpth\in\dot{A}_{i+1}$. Then $\allpth\in\dot{A}$ and $\lim_{\ell\to+\infty}\allpth(u_{i+1,\ell})\in A_{i+1}\subseteq\st{x_{u_{i+1}}\in\stsp_{u_{i+1}}\colon x_{u_i}\in A_i}$. So there is some $x_{u_i}\in A_i$ such that, for each $t\in u_i$, the corresponding component of $\lim_{\ell\to+\infty}\allpth(u_{i+1,\ell})$ is equal to $x_t$. Furthermore, if $t$ is a right-sided limit point of $\densepoints$, then the corresponding components of $\lim_{\ell\to+\infty}\allpth(u_{i+1,\ell})$ and $\lim_{\ell\to+\infty}\allpth(u_{i,\ell})$ are equal because $\lim_{\densepoints\ni r\searrow t} \allpth\pr{r}$ exists due to Lemma~\ref{lem:limit along densetpoints}, and if $t$ is not a right-sided limit point of $\densepoints$, then the corresponding components of $\lim_{\ell\to+\infty}\allpth(u_{i+1,\ell})$ and $\lim_{\ell\to+\infty}\allpth(u_{i,\ell})$ are equal because they are both equal to $\allpth(t)$. In all cases, we conclude that for each $t\in u_i$, the corresponding component of $\lim_{\ell\to+\infty}\allpth(u_{i,\ell})$ is equal to $x_t$. So $\lim_{\ell\to+\infty}\allpth(u_{i,\ell})=x_{u_i}\in A_i$ and therefore, since $\allpth\in\dot{A}$, we see that $\allpth\in\dot{A}_i$, as required.

The second property is that \(\bigcap_{i\in\nats} \dot{A}_i=\emptyset\).
	To see why this is true, assume \emph{ex absurdo} that \(\bigcap_{i\in\nats} \dot{A}_i\neq\emptyset\).
	Then there is some~\(\allpth\in\dot{A}\) such that for all \(i\in\nats\), \(\lim_{\ell\to+\infty}\allpth\pr{u_{i,\ell}}\in A_i\).
	By Lemma~\ref{lem:limit along densetpoints}, there is some~\(\cadpth\in\cadpaths\) such that for all \(i\in\nats\), \(\cadpth\pr{u_i}=\lim_{\ell\to+\infty}\allpth\pr{u_{i,\ell}}\in A_i\); for the components $t$ of $u_i$ that are right-sided limit points of $\densepoints$, the equality in this statement follows because $\cadpth(t)=\lim_{\densepoints\ni r\searrow t} \allpth\pr{r}$, and for the components $t$ of $u_i$ that are not right-sided limit points of $\densepoints$, this follows because $\cadpth(t)=\allpth\pr{t}$ and because the corresponding components of $u_{i,\ell}$ are all equal to $t$.
	But then \(\cadpth\in \bigcap_{i\in\nats} \ev{\cadX_{u_i}\in A_i}=\bigcap_{i\in\nats} \mathring{A}_i\), which is a contradiction because \(\bigcap_{i\in\nats} \mathring{A}_i=\emptyset\).

	Since \(\allprob\) is a probability measure -- and therefore definitely continuous for decreasing sequences -- it follows from these two properties that
	\begin{equation}
		\label{eqn:proof of main:intermed5}
		\lim_{i\to+\infty} \allprob\pr{\dot{A}_i}
		= \allprob\pr{\emptyset}
		= 0.
	\end{equation}
Finally, it follows from Eqns.~\eqref{eqn:proof of main:intermed1} to \eqref{eqn:proof of main:intermed5} -- and the non-negativity of~\(\prob\) -- that
	\begin{align*}
		0\leq\lim_{i\to+\infty} \prob\pr{\mathring{A}_i}
		&= \lim_{i\to+\infty} \allprob\pr{\ev{\allX_{u_i}\in A_i}}\\
		&= \lim_{i\to+\infty} \allprob\pr{\dot{A}\cap\ev{\allX_{u_i}\in A_i}}\\
		&\leq\limsup_{\ell\to+\infty} \allprob\pr{\dot{A}\cap\ev{\allX_{u_{i,\ell}}\in A_i}}
		= \lim_{i\to+\infty} \allprob\pr{\dot{A}_i}
		= 0.
	\end{align*}
	So we find that $\lim_{i\to+\infty} \prob\pr{\mathring{A}_i}$, which finalises our proof of the sufficiency.
\end{proof}

\section{Conditions for regularity}
\label{asec:two additional sufficient conditions for regularity}

The two regularity conditions \ref{def:dist regular:cont} and \ref{def:dist regular:bound on exp of jumps} are perhaps not the most easy to check.
We have already given a simpler sufficient condition in Proposition~\ref{prop:only expected number of jumps for regularity} in the main text.
In this appendix, we give two alternative sufficient conditions in the special case that \(\timepoints\cap\cci{-n}{n}\) is closed for all $n\in\nats$.
These only involve the finite-dimensional distributions for two time points, and can be interpreted as bounding the dynamics of the stochastic process.
\begin{proposition}
\label{prop:stronger I}
	Suppose that \(\timepoints\cap\cci{-n}{n}\) is closed for all $n\in\nats$.
	Consider a consistent collection~\(\dist_{\noarg}\) of finite-dimensional distributions.
	If for all \(n\in\nats\) there is some~\(\lambda_n\in\nnegreals\) such that
	\begin{multline*}
		\pr{\forall t\in\timepoints\cap\cci{-n}{n}}\pr{\forall\epsilon\in\posreals}\pr{\exists\delta\in\posreals}\pr{\forall s,r\in\timepoints\cap\ooi{t-\delta}{t+\delta}\colon s< r}\\
		\frac{\dist_{\pr{s,r}}\pr{\stsp^2_{\neq}}}{r-s}
		< \lambda_n+\epsilon,
	\end{multline*}
	then \(\dist_{\noarg}\) is regular.
\end{proposition}
In our proof, we will lean on the following two intermediary results, which will come in handy in our proof for Proposition~\ref{prop:stronger II} further on as well.
The first one is an immediate consequence of the Heine-Borel Theorem.
\begin{lemma}
\label{lem:Modified Heine Borel}
	Consider a closed and bounded subset~\(\altpoints\) of~\(\timepoints\).
	For all \(s\in \altpoints\), fix some positive real number~\(\delta_s\).
	Then there is a tuple of time points~\(\pr{s_1, \dots, s_m}\in \setoftseq_{\altpoints}\) such that, with \(\Delta_\ell\coloneqq\delta_{s_\ell}\) for all \(\ell\in \st{1, \dots, m}\),
	\begin{enumerate}[label=\upshape(\roman*), ref=\upshape(\roman*)]
		\item\label{lem:MHB:cover} \(\altpoints\subseteq \bigcup_{\ell=1}^m \ooi{s_\ell-\Delta_\ell}{s_\ell+\Delta_\ell}\); and
		\item\label{lem:MHB:no complete overlap} \(s_i-\Delta_i<s_j-\Delta_j\) and \(s_i+\Delta_i<s_j+\Delta_j\) for all \(i,j\in \st{1, \dots, m}\) such that \(i<j\).
	\end{enumerate}
	Hence, for any sequence~\(u=\pr{t_1, \dots, t_n}\in \setoftseq_{\altpoints}\) with \(n\geq 2\),
	\begin{equation}
	\label{lem:MHB:at most m}
		\card[\big]{\st[\big]{k\in\st{2, \dots, n}\colon \pr{\not\exists\ell\in\st{1, \dots, m}}~\st{t_{k-1}, t_k}\subseteq \ooi{s_\ell-\Delta_\ell}{s_\ell+\Delta_\ell}}}
		< m.
	\end{equation}
	If furthermore \(\altpoints\) is convex, then
	\begin{equation}
		\label{lem:MHB:overlap}
		s_{\ell-1}+\Delta_{\ell-1}>s_{\ell}-\Delta_\ell \text{ for all } \ell\in\st{2, \dots, m}.
	\end{equation}
\end{lemma}
\begin{proof}
	Statements~\ref{lem:MHB:cover} and \ref{lem:MHB:no complete overlap} and Eqn.~\eqref{lem:MHB:overlap} essentially follow from an argument in \cite[Proof of Lemma~F.1]{2017KrakDeBock}, which we modify to our slightly different context.
	For all \(s\in \altpoints\), we let \(C_s\coloneqq\ooi{s-\delta_s}{s+\delta_s}\).
	Note that \(\pr{C_s}_{s\in\altpoints}\) is an open cover of the bounded and closed set~\(\altpoints\).
	By the Heine-Borel Theorem, this open cover has a finite subcover: \(C_{s_1}\), \dots, \(C_{s_m}\).
	Without loss of generality, we may assume that \(s_1<\cdots<s_m\) and that this subcover is minimal, in the sense that removing one of its elements has as a consequence that we no longer have a subcover of~\(\altpoints\) any more.
	This proves \ref{lem:MHB:cover}.

	Our proof for \ref{lem:MHB:no complete overlap} is one by contradiction.
	Assume \emph{ex absurdo} that there are indices~\(i,j\in \st{1, \dots, m}\) such that \(i<j\) and either \(s_i-\Delta_i\geq s_j-\Delta_j\) or \(s_i+\Delta_i\geq s_j+\Delta_j\).
	Let us start with the first case.
	Since \(i<j\) by assumption, it follows that \(s_i<s_j\), and therefore, since \(s_i-\Delta_i\geq s_j-\Delta_j\) by assumption, that \(\Delta_j\geq s_j-s_i+\Delta_i>\Delta_i\).
	From this, it follows that \(s_i+\Delta_i<s_i+\Delta_j<s_j+\Delta_j\), where for the second inequality we again used that \(s_i<s_j\).
	Hence, \(C_{s_i}=\ooi{s_i-\Delta_i}{s_i+\Delta_i}\subseteq\ooi{s_j-\Delta_j}{s_j+\Delta_j}=C_{s_j}\), but this is a contradiction because \(C_{s_1}\), \dots \(C_{s_m}\) is minimal.
	For the second case that \(i<j\) and \(s_i+\Delta_i\geq s_j+\Delta_j\), a similar argument leads to a contradiction as well.

	Next, we fix any \(u=\pr{t_1, \dots, t_n}\in \setoftseq_{\altpoints}\) with \(n\geq 2\), and set out to prove Eqn.~\eqref{lem:MHB:at most m}.
	For all \(k\in \st{1, \dots, n}\), the time point~\(t_k\) belongs to~\(\altpoints\), so it follows from
	\ref{lem:MHB:cover} that the index set
	\begin{equation*}
		\mathcal{L}_{k}
		\coloneqq \st[\big]{\ell\in\st{1, \dots, m}\colon t_{k}\in C_{s_\ell}}
	\end{equation*}
	is non-empty and finite.
	Furthermore, for all \(k\in \st{1, \dots, n}\),
	\begin{equation}
	\label{eqn:proof of MHB:intermed1}
		\mathcal{L}_k
		= \st{\ell\in\nats\colon \min\mathcal{L}_k\leq\ell\leq\max\mathcal{L}_k}.
	\end{equation}
	Indeed, if this were not the case, there would be indices~\(i,\ell, j\in \st{1, \dots, m}\) such that \(i<\ell<j\), \(t_k\in C_{s_i}\cap C_{s_j}\) and \(t_k\notin C_{s_\ell}\).
	Then either \(t_k\leq s_\ell-\Delta_\ell\) or \(t_k\geq s_\ell+\Delta_\ell\), but by \ref{lem:MHB:no complete overlap} this implies that \(t_k<s_j-\Delta_j\) -- so \(t_k\notin C_{s_j}\) -- or \(t_k>s_i+\Delta_i\) -- so \(t_k\notin C_{s_i}\) -- which is a contradiction because \(t_k\in C_{s_i}\cap C_{s_j}\) by assumption.
	Finally, for all indices \(k_1,k_2\in \st{1, \dots, n}\) such that \(k_1<k_2\),
	\begin{equation}
	\label{eqn:proof of MHB:intermed2}
		\min \mathcal{L}_{k_1}
		\leq \min \mathcal{L}_{k_2}
		\text{ and }
		\max\mathcal{L}_{k_1}
		\leq\max\mathcal{L}_{k_2};
	\end{equation}
	again, this follows more or less immediately from \ref{lem:MHB:no complete overlap}.
	Indeed, assume \emph{ex absurdo} that \(j\coloneqq\min \mathcal{L}_{k_1}>\min \mathcal{L}_{k_2}\eqqcolon i\).
	Then \(t_{k_1}\in \ooi{s_j-\Delta_j}{s_j+\Delta_j}\) and \(t_{k_1}\notin \ooi{s_i-\Delta_i}{s_i+\Delta_i}\ni t_{k_2}\), and by \ref{lem:MHB:no complete overlap}, \(s_i-\Delta_i<s_j-\Delta_j\).
	Hence, either (i) \(t_{k_1}\leq s_i-\Delta_i\) and therefore \(t_{k_1}<s_j-\Delta_j<t_{k_1}\), or (ii) \(t_{k_1}\geq s_i+\Delta_i\) and therefore \(t_{k_1}>t_{k_2}\); in both cases, we end up with a clear contradiction.
	A similar argument shows that \(\max\mathcal{L}_{k_1}\leq\max\mathcal{L}_{k_2}\).

	To prove Eqn.~\eqref{lem:MHB:at most m}, we observe that, by Eqns.~\eqref{eqn:proof of MHB:intermed1} and \eqref{eqn:proof of MHB:intermed2},
	\begin{align*}
		\MoveEqLeft
		\st[\big]{k\in\st{2, \dots, n}\colon \pr{\not\exists\ell\in\st{1, \dots, m}}~\st{t_{k-1}, t_k}\subseteq \ooi{s_\ell-\Delta_\ell}{s_\ell+\Delta_\ell}} \\
		&= \st[\big]{k\in\st{2, \dots, n}\colon \max\mathcal{L}_{k-1}<\min\mathcal{L}_k} \\
		&\subseteq \st[\big]{k\in\st{2, \dots, n}\colon \max\mathcal{L}_{k-1}<\max\mathcal{L}_k}.
	\end{align*}
	Since \(1\leq\max\mathcal{L}_1\leq \cdots \leq\max\mathcal{L}_n\leq m\) by Eqn.~\eqref{eqn:proof of MHB:intermed2}, we infer from this that
	\begin{equation*}
		\card[\big]{\st[\big]{k\in\st{2, \dots, n}\colon \pr{\not\exists\ell\in\st{1, \dots, m}}~\st{t_{k-1}, t_k}\subseteq \ooi{s_\ell-\Delta_\ell}{s_\ell+\Delta_\ell}}}
		\leq m-1,
	\end{equation*}
	as required.

	Finally, we prove the final part of the statement, so we assume that \(\altpoints\) is convex.
	Assume \emph{ex absurdo} that there is an index~\(\ell^\star\in \st{2, \dots, m}\) such that
	\begin{equation*}
		s_{\ell^\star-1}+\Delta_{\ell^\star-1}
		< s_{\ell^\star}-\Delta_{\ell^\star}.
	\end{equation*}
	Fix any \(t\in \ooi{s_{\ell^\star-1}+\Delta_{\ell^\star-1}}{s_{\ell^\star}-\Delta_{\ell^\star}}\).
	Note that \(\altpoints\) includes \(\cci{s_{\ell^\star-1}}{s_{\ell^\star}}\) since \(\altpoints\) is closed and convex, so  \(t\) belongs to~\(\altpoints\).
	However, it follows from \ref{lem:MHB:no complete overlap} that \(s_i+\Delta_i\leq s_{\ell^\star-1}-\Delta_{\ell^\star-1}\) for all \(i\in \st{1, \dots, \ell^\star-1}\) and \(s_{\ell^\star}-\Delta_{\ell^\star}\leq s_j-\Delta_j\) for all \(j\in \st{\ell^\star, \dots, m}\).
	From this, we infer that \(t\) does not belong to~\(\bigcup_{\ell=1}^m C_{s_\ell}\), which is a contradiction because \(C_{s_1}\), \dots, \(C_{s_m}\) is a cover of~\(\altpoints\).
\end{proof}
The second one is an obvious observation about the expected number of jumps.
\begin{lemma}
\label{lem:exp jumps as sum for dist}
	Consider a consistent collection~\(\dist_{\noarg}\) of finite-dimensional charges.
	Then for all \(u,v\in\setoftseq\) such that \(u\sqsubseteq v\),
	\begin{equation*}
		\prev_{\dist_u}\pr{\tupjumps_u}
		\leq \prev_{\dist_v}\pr{\tupjumps_v}.
	\end{equation*}
\end{lemma}
\begin{proof}
	Let us enumerate the time points in~\(u\) as \(\pr{r_1, \dots, r_n}\) and in~\(v\) as \(\pr{s_1, \dots, s_m}\).
	It follows from Eqn.~\eqref{eqn:prevdist jump} and the consistency of~\(\dist_{\noarg}\) that
	\begin{equation*}
		\prev_{\dist_u}\pr{\tupjumps_u}
		= \sum_{k=2}^n \dist_{u}\pr{\st{x_u\in\stsp_u\colon x_{r_{k-1}}\neq x_{r_k}}}
		= \sum_{k=2}^n \dist_{v}\pr{\st{y_v\in\stsp_v\colon y_{r_{k-1}}\neq y_{r_k}}}.
	\end{equation*}
	Let \(\ell_1, \dots, \ell_n\) be as defined in the proof of Lemma~\ref{lem:jumps is monotone}.
	Then as explained there, for all \(k\in\st{2, \dots, n}\),
	\begin{equation*}
		\st{y_v\in\stsp_v\colon y_{r_{k-1}}\neq y_{r_k}}
		\subseteq \bigcup_{\ell=\ell_{k-1}+1}^{\ell_k} \st{y_v\in\stsp_v\colon y_{s_{\ell-1}}\neq y_{s_\ell}}.
	\end{equation*}
	Together with the (sub-)additivity of~\(\dist_v\), this implies that
	\begin{align*}
		\prev_{\dist_u}\pr{\tupjumps_u}
		&\leq \sum_{k=2}^n \sum_{\ell=\ell_{k-1}+1}^{\ell_k} \dist_v\pr{\st{y_v\in\stsp_v\colon y_{s_{\ell-1}}\neq y_{s_\ell}}}\\
		&\leq \sum_{\ell=2}^m\dist_v\pr{\st{y_v\in\stsp_v\colon y_{s_{\ell-1}}\neq y_{s_\ell}}}.
	\end{align*}
	The inequality in the statement follows immediately due to Eqn.~\eqref{eqn:prevdist jump}.
\end{proof}
\begin{proof}[Proof of Proposition~\ref{prop:stronger I}]
	To verify \ref{def:dist regular:cont}, we fix any \(t\in\timepoints\cap\rlims\pr{\timepoints}\) and \(\epsilon\in\posreals\).
	Fix any \(n\in\nats\) such that \(\abs{t}<n\).
	Then by the condition in the statement, there is some \(\delta'\in\posreals\) such that for all $r\in\timepoints\cap\ooi{t}{t+\delta'}$
	\begin{equation*}
		\dist_{\pr{t,r}}\pr[\big]{\stsp^2_{\neq}}
		< \pr{r-t}\pr{\lambda_n+\epsilon}.
	\end{equation*}
	Let \(\delta\coloneqq\min\st{\delta', \epsilon/\pr{\lambda_n+\epsilon}}\).
	Then for all $r\in\timepoints\cap\ooi{t}{t+\delta}$
	\begin{align*}
		\dist_{\pr{t,r}}\pr[\big]{\stsp^2_{=}}
		= 1 - \dist_{\pr{t,r}}\pr[\big]{\stsp^2_{\neq}}
		&> 1-(r-t)(\lambda_n+\epsilon)\\
		&\geq1-\delta(\lambda_n+\epsilon)
		\geq1-\epsilon.
	\end{align*}
	Since \(\epsilon\) was an arbitrary strictly positive real number, we conclude that \(\lim_{\timepoints\ni r\searrow t}\dist_{\pr{t,r}}\pr[\big]{\stsp^2_{=}}=1\), as required.

	To verify \ref{def:dist regular:bound on exp of jumps}, we fix any \(n\in\nats\) such that \(\cci{-n}{n}\cap\timepoints\neq\emptyset\).
	Furthermore, we fix any \(\epsilon\in\posreals\).
	Then by the condition in the statement, for all \(t\in\timepoints\cap\cci{-n}{n}\) there is some \(\delta_t\in\posreals\) such that
	\begin{equation}
	\label{eqn:proof of ex of jumps bounded:bound}
		\pr[\big]{\forall s,r\in\timepoints\cap\ooi{t-\delta_t}{t+\delta_t}\colon s<r}~
		\dist_{\pr{s,r}}\pr{\stsp^2_{\neq}}
		<\pr{r-s}\pr{\lambda_n+\epsilon}.
	\end{equation}
	Because \(\timepoints\cap\cci{-n}{n}\) is closed by assumption and clearly bounded, it follows from Lemma~\ref{lem:Modified Heine Borel} that there is some tuple of time points~\(\pr{s_1, \dots, s_m}\in\setoftseq_{\cci{-n}{n}\cap\timepoints}\) such that
	\begin{equation}
	\label{eqn:proof of ex of jumps bounded:cover}
		\timepoints\cap\cci{-n}{n}
		\subseteq \bigcup_{\ell=1}^m \ooi{s_\ell-\delta_{s_\ell}}{s_\ell+\delta_{s_\ell}}
	\end{equation}
	and, for all \(u=\pr{t_1, \dots, t_k}\in\setoftseq_{\cci{-n}{n}\cap\timepoints}\) with \(k\geq2\),
	\begin{equation}
	\label{eqn:proof of ex of jumps bounded:cover2}
		\card[\big]{\st[\big]{i\in\st{2, \dots, k}\colon \pr{\not\exists\ell\in\st{1, \dots, m}}~\st{t_{i-1}, t_i}\subseteq \ooi{s_\ell-\delta_{s_\ell}}{s_\ell+\delta_{s_\ell}}}}
		< m.
	\end{equation}

	Fix any grid~\(u=\pr{t_1, \dots, t_k}\in \setoftseq_{\cci{-n}{n}\cap\timepoints}\), and recall from Eqn.~\eqref{eqn:prevdist jump as sum over differences} that
	\begin{equation*}
		\prev_{\dist_u}\pr{\tupjumps_u}
		= \sum_{i=2}^k \dist_{\pr{t_{i-1}, t_i}}\pr[\big]{\stsp^2_{\neq}}.
	\end{equation*}
	Let us split the index set~\(\st{2, \dots, k}\) into two parts:
	\begin{equation*}
		\mathcal{I}_1
		\coloneqq \st[\big]{i\in\st{2, \dots, k}\colon\pr{\exists\ell\in\st{1, \dots, m}}~\st{t_{i-1}, t_i}\subseteq\ooi{s_\ell-\delta_{s_\ell}}{s_\ell+\delta_{s_\ell}}}
	\end{equation*}
	and \(\mathcal{I}_2\coloneqq\st{2, \dots, k}\setminus\mathcal{I}_2\).

	On the one hand, it follows from Eqn.~\eqref{eqn:proof of ex of jumps bounded:bound} that for all \(i\in \mathcal{I}_1\),
	\begin{equation*}
		\dist_{\pr{t_{i-1}, t_i}}\pr[\big]{\stsp^2_{\neq}}
		< \pr{t_i-t_{i-1}}\pr{\lambda_n + \epsilon}.
	\end{equation*}
	On the other hand, we know from Eqn.~\eqref{eqn:proof of ex of jumps bounded:cover2} that there are at most \(m-1\) indices in~\(\mathcal{I}_2\), and for these indices~\(i\) clearly
	\begin{equation*}
		\dist_{\pr{t_{i-1}, t_i}}\pr[\big]{\stsp^2_{\neq}}
		\leq 1
		< \pr{t_i-t_{i-1}}\pr{\lambda_n + \epsilon} + 1.
	\end{equation*}
	Thus, we find that
	\begin{align*}
		\prev_{\dist_u}\pr{\tupjumps_u}
		&= \sum_{i=2}^k \dist_{\pr{t_{i-1}, t_i}}\pr[\big]{\stsp^2_{\neq}} \\
		&= \sum_{i\in\mathcal{I}_1} \dist_{\pr{t_{i-1}, t_i}}\pr[\big]{\stsp^2_{\neq}} + \sum_{i\in\mathcal{I}_2} \dist_{\pr{t_{i-1}, t_i}}\pr[\big]{\stsp^2_{\neq}} \\
		&\leq m-1+\sum_{i=2}^k \pr{t_i-t_{i-1}}\pr{\lambda_n+\epsilon} \\
		&\leq m-1+\pr{t_k-t_1}\pr{\lambda_n+\epsilon}.
	\end{align*}
	This inequality holds for all \(u=\pr{t_1, \dots, t_k}\in\setoftseq_{\cci{-n}{n}\cap\timepoints}\) -- so with \(-n\leq t_1\leq t_k\leq n\) -- and therefore
	\begin{equation*}
		\sup\st[\big]{\prev_{\dist_u}\pr{\tupjumps_u}\colon u\in\setoftseq_{\cci{-n}{n}\cap\timepoints}}
		\leq m-1+2n\pr{\lambda_n+\epsilon}
		<+\infty.
	\end{equation*}
	A similar argument as the one in the second part of the proof of Proposition~\ref{prop:only expected number of jumps for regularity} shows that this implies \ref{def:dist regular:bound on exp of jumps}.
\end{proof}

For the third sufficient condition, we assume that for all \(n\in\nats\), \(\timepoints\cap\cci{-n}{n}\) is some finite union of intervals -- bounded, closed and convex subsets of~\(\timeaxis\) -- and we impose a condition on the limit superior of the rate of change of the probability of changing states.
The argument in the second part of our proof generalises the proof of Lemma~5.49 in \cite{2021Erreygers-Phd}.
\begin{proposition}
\label{prop:stronger II}
	Suppose that for all \(n\in\nats\) such that \(\timepoints\cap\cci{-n}{n}\neq\emptyset\), there are time points~\(s_{n,1}\leq r_{n,1}<\cdots<s_{n,m_n}\leq r_{n,m_n}\) such that \(\timepoints\cap\cci{-n}{n}=\bigcup_{k=1}^{m_n} \cci{s_{n,k}}{r_{n,k}}\).
	Consider a consistent collection~\(\dist_{\noarg}\) of finite-dimensional distributions.
	If for all \(n\in\nats\) such that $\timepoints\cap\cci{-n}{n}\neq\emptyset$ there is some~\(\lambda_n\in\nnegreals\) such that for all \(k\in\st{1, \dots, m_n}\),
	\begin{align*}
		\limsup_{r\searrow t}\frac{\dist_{\pr{t,r}}\pr[\big]{\stsp^2_{\neq}}}{r-t}
		&\leq \lambda_n
		\quad\text{for all } t\in\coi{s_{n,k}}{r_{n,k}}
	\shortintertext{and}
		\limsup_{s\nearrow t}\frac{\dist_{\pr{s,t}}\pr[\big]{\stsp^2_{\neq}}}{t-s}
		&\leq \lambda_n
		\quad\text{for all } t\in\oci{s_{n,k}}{r_{n,k}},
	\end{align*}
	then \(\dist_{\noarg}\) is regular.
\end{proposition}
\begin{proof}
	The proof for \ref{def:dist regular:cont} is almost exactly the same as the one in our proof for Proposition~\ref{prop:stronger I}, although it does become a bit simpler.
	Fix any \(t\in\timepoints\cap\rlims\pr{\timepoints}\) and some \(n\in\nats\) such that \(\abs{t}<n\).
	Then it follows from the condition in the statement that there are time points~\(s^\star,r^\star\in\timepoints\) and a non-negative real number~\(\lambda_n\in\nnegreals\) such that \(t\in\coi{s^\star}{r^\star}\subseteq\timepoints\) and
	\begin{equation*}
		\limsup_{r\searrow t} \frac{\dist_{\pr{t,r}}\pr[\big]{\stsp^2_{\neq}}}{r-t}
		\leq \lambda_n.
	\end{equation*}
	Clearly, this can only be the case if \(\lim_{r\searrow t} \dist_{\pr{t,r}}\pr[\big]{\stsp^2_{\neq}}=0\).
	Since \(\dist_{\pr{t,r}}\pr[\big]{\stsp^2_{=}}=1-\dist_{\pr{t,r}}\pr[\big]{\stsp^2_{\neq}}\) for all \(r\in\oci{t}{r^\star}\), we infer from this that
	\begin{equation*}
		\lim_{r\searrow t} \dist_{\pr{t,r}}\pr[\big]{\stsp^2_{=}}
		= \lim_{r\searrow t} 1-\dist_{\pr{t,r}}\pr[\big]{\stsp^2_{\neq}}
		= 1,
	\end{equation*}
	as required.

	To verify \ref{def:dist regular:bound on exp of jumps}, we fix any \(n\in\nats\) such that \(\cci{-n}{n}\cap\timepoints\neq\emptyset\).
	Then by the conditions in the statement there are time points~\(s_{n,1}\leq r_{n,1}<\cdots<s_{n,m_n}\leq r_{n,m_n}\in\timepoints\) such that \(\cci{-n}{n}\cap\timepoints=\bigcup_{k=1}^{m_n}\cci{s_{n,k}}{r_{n,k}}\).
	Let \(u=\pr{t_1, \dots, t_m}\) be an arbitrary tuple of time points in~\(\setoftseq_{\cci{-n}{n}\cap\timepoints}\), and recall from Eqn.~\eqref{eqn:prevdist jump as sum over differences} that
	\begin{equation*}
		\dist_{u}\pr{\tupjumps_u}
		= \sum_{\ell=2}^m \dist_{\pr{t_{\ell-1}, t_\ell}}\pr{\stsp^2_{\neq}}.
	\end{equation*}
	For all \(\ell\in\st{1, \dots, m}\), there is a unique index~\(k_\ell\in\st{1, \dots, m_n}\) such that \(t_\ell\in\cci{s_{n,k_\ell}}{r_{n,k_\ell}}\).
	It is obvious that \(k_1\leq\cdots\leq k_m\) because \(t_1<\cdots<t_m\), so the index set
	\begin{equation*}
		\mathcal{L}_1
		\coloneqq \st[\big]{\ell\in\st{2, \dots, m}\colon k_{\ell-1}<k_\ell}
	\end{equation*}
	contains at most \(m_n-1\) indices.
	For these \(m_n-1\) indices~\(\ell\in\mathcal{L}_1\), the most we can say is that
	\begin{equation*}
		\dist_{\pr{t_{\ell-1}, t_\ell}}\pr{\stsp^2_{\neq}}
		\leq 1.
	\end{equation*}
	However, we can say more for those indices~\(\ell\) for which \(t_{\ell-1}\) and \(t_\ell\) belong to the same interval~\(\cci{s_{n, k}}{r_{n, k}}\); that is, for those indices in the index set
	\begin{equation*}
		\mathcal{L}_2
		\coloneqq \st{2, \dots, m}\setminus\mathcal{L}_1
		= \st[\big]{\ell\in\st{2, \dots, m}\colon k_{\ell-1}=k_\ell}.
	\end{equation*}

	To ease our notation, we fix two time points~\(s,r\in\timepoints\) such that \(s<r\) and \(\cci{s}{r}\subseteq\cci{s_{n, k}}{r_{n, k}}\) for some \(k\in\st{1, \dots, m_n}\).
	We set out to show that
	\begin{equation}
	\label{eqn:proof of stronger II:inequality for jumps}
		\dist_{\pr{s,r}}\pr{\stsp^2_{\neq}}
		\leq \pr{r-s} \lambda_n.
	\end{equation}
	To this end, we recall from Eqn.~\eqref{eqn:prevdist jump as sum over differences} that
	\begin{equation}
	\label{eqn:proof of stronger II:prob to jumps}
		\dist_{\pr{s,r}}\pr{\stsp^2_{\neq}}
		= \prev_{\dist_{\pr{s,r}}}\pr{\tupjumps_{\pr{s,r}}}.
	\end{equation}
	By the assumptions in the statement, there is some \(\lambda_n\in\nnegreals\) such that for all \(t\in\cci{s}{r}\) and \(\epsilon\in\posreals\), there are positive real numbers~\(\delta_t^+, \delta_t^-\in\posreals\) such that
	\begin{align}
		\pr[\big]{\forall r'\in\cci{s}{r}\cap\ooi{t}{t+\delta_t^+}}~
		\frac{\dist_{\pr{t,r'}}\pr{\stsp^2_{\neq}}}{r'-t}
		&< \lambda_n +\epsilon
		\label{eqn:proof of stronger II:lambda bound right}
	\shortintertext{and}
		\pr[\big]{\forall s'\in\cci{s}{r}\cap\ooi{t-\delta_t^-}{t}}~
		\frac{\dist_{\pr{s',t}}\pr{\stsp^2_{\neq}}}{t-s'}
		&< \lambda_n +\epsilon.
		\label{eqn:proof of stronger II:lambda bound left}
	\end{align}
	Fix any \(\epsilon\in\posreals\), and for all \(t\in\cci{s}{r}\), let \(\delta_t\coloneqq\min\st{\delta_t^+, \delta_t^-}\).
	Since~\(\cci{s}{r}\) is clearly a closed, bounded and convex subset of~\(\timepoints\), it follows from Lemma~\ref{lem:Modified Heine Borel} that there is a tuple of time points~\(\pr{s_1, \dots, s_{p}}\in\setoftseq_{\cci{s}{r}}\) such that, with \(\Delta_i\coloneqq\delta_{s_i}\) for all \(i\in\st{1, \dots, p}\),
	\begin{enumerate}
		\item \(\cci{s}{r}\subseteq\bigcup_{i=1}^{p} \ooi{s_i-\Delta_i}{s_i+\Delta_i}\),
		\item \(s_i-\Delta_i<s_j-\Delta_j\) and \(s_i+\Delta_i<s_j+\Delta_j\) for all \(i,j\in\st{1, \dots, p}\) such that \(i<j\), and
		\item \(s_{i-1}+\Delta_{i-1}>s_i-\Delta_i\) for all \(i\in\st{2, \dots, p}\).
	\end{enumerate}
	Let~\(v=\pr{r_1, \dots, r_{q}}\in\setoftseq_{\cci{s}{r}}\) be a tuple of time points that (i) starts in~\(r_1=s\), (ii) ends in~\(r_q=r\), (iii) contains all the time points~\(s_1\), \dots, \(s_{p}\), and (iv) for all \(i\in\st{2, \dots, p}\), contains one time point in~\(\ooi{s_{i-1}}{s_i}\cap\ooi{s_i-\Delta_i}{s_{i-1}+\Delta_{i-1}}\) -- to see that this intersection is non-empty, recall the three properties above.
	Note that \(q=2p+1\) if \(s<s_1\) and \(s_p<r\), \(q=2p\) if either \(s=s_1\) or \(s_p=r\), and \(q=2p-1\) if \(s=s_1\) and \(s_p=r\).
	Then by construction \(\pr{s,r}\sqsubseteq v\) and, due to Eqns.~\eqref{eqn:proof of stronger II:lambda bound right} and \eqref{eqn:proof of stronger II:lambda bound left},
	\begin{equation*}
		\dist_{\pr{r_{i-1}, r_i}}\pr{\stsp^2_{\neq}}
		< \pr{r_i-r_{i-1}}\pr{\lambda_n+\epsilon}
		\quad\text{for all } i\in\st{2, \dots, q}.
	\end{equation*}
	From this and Eqn.~\eqref{eqn:prevdist jump as sum over differences}, we infer that
	\begin{equation*}
		\prev_{\dist_v}\pr{\tupjumps_v}
		= \sum_{i=2}^{q} \dist_{\pr{r_{i-1}, r_i}}\pr{\stsp^2_{\neq}} \\
		< \sum_{i=2}^{q} \pr{r_i-r_{i-1}}\pr{\lambda_n+\epsilon} \\
		= \pr{r-s} \pr{\lambda_n+\epsilon}.
	\end{equation*}
	From this, Eqn.~\eqref{eqn:proof of stronger II:prob to jumps} and Lemma~\ref{lem:exp jumps as sum for dist}, it now follows that
	\begin{equation*}
		\dist_{\pr{s,r}}\pr{\stsp^2_{\neq}}
		= \prev_{\dist_{\pr{s,r}}}\pr{\tupjumps_{\pr{s,r}}}
		\leq \prev_{\dist_v}\pr{\tupjumps_v}
		< \pr{r-s} \pr{\lambda_n+\epsilon}.
	\end{equation*}
	Since this inequality holds for arbitrary \(\epsilon\in\posreals\), we have proven Eqn.~\eqref{eqn:proof of stronger II:inequality for jumps}.

	Let us return to where we were before.
	Due to Eqn.~\eqref{eqn:proof of stronger II:inequality for jumps}, we know that for all \(\ell\in\mathcal{L}_2\),
	\begin{equation*}
		\dist_{\pr{t_{\ell-1},t_\ell}}\pr{\stsp^2_{\neq}}
		\leq \pr{t_\ell-t_{\ell-1}} \lambda_n.
	\end{equation*}
	For this reason, and because \(\dist_{\pr{t_{\ell-1},t_\ell}}\pr{\stsp^2_{\neq}}\leq 1\) for all \(\ell\in\mathcal{L}_1=\st{2, \dots, m}\setminus\mathcal{L}_2\),
	\begin{align*}
		\dist_u\pr{\tupjumps_u}
		&= \sum_{\ell=2}^m \dist_{\pr{t_{\ell-1},t_\ell}}\pr{\stsp^2_{\neq}} \\
		&= \sum_{\ell\in\mathcal{L}_1} \dist_{\pr{t_{\ell-1},t_\ell}}\pr{\stsp^2_{\neq}} + \sum_{\ell\in\mathcal{L}_2} \dist_{\pr{t_{\ell-1},t_\ell}}\pr{\stsp^2_{\neq}} \\
		&\leq \sum_{\ell\in\mathcal{L}_1} 1 + \sum_{\ell\in\mathcal{L}_2} \pr{t_\ell-t_{\ell-1}}\lambda_n \\
		&\leq m_n-1 + \pr{t_m-t_1}\lambda_n \\
		&\leq m_n-1+2n\lambda_n,
	\end{align*}
	where for the second inequality we used that \(\card{\mathcal{L}_1}\leq m_n-1\) and for the final inequality we used that \(u\in\setoftseq_{\cci{-n}{n}\cap\timepoints}\).
	This inequality holds for any \(u=\pr{t_1, \dots, t_m}\in\setoftseq_{\cci{-n}{n}\cap\timepoints}\), so we conclude that
	\begin{equation*}
		\sup\st[\big]{\prev_{\dist_u}\pr{\tupjumps_u}\colon u\in\setoftseq_{\cci{-n}{n}\cap\timepoints}}
		\leq m_n-1+2n\lambda_n
		< +\infty.
	\end{equation*}
	Here too, \ref{def:dist regular:bound on exp of jumps} follows from the same argument as that in the second part of the proof of Proposition~\ref{prop:only expected number of jumps for regularity}.
\end{proof}

\end{document}